\documentclass[leqno,final]{siamltex}

\usepackage{amsmath}
\usepackage{amssymb}
\usepackage{amsfonts}
\usepackage{graphicx}

\newcommand{\ex}{\mathrm{e}}
\newcommand{\E}{\rm{E}}
\newcommand{\I}{\mathrm{i}}

\DeclareMathOperator{\meas}{meas}

\newtheorem{remark}{Remark}[section]
\newtheorem{example}{Example}[section]
{\catcode `\@=11 \global\let\AddToReset=\@addtoreset}
\AddToReset{equation}{section}
\renewcommand{\theequation}{\thesection.\arabic{equation}}
\newcommand{\N}{{\mathbb N }}
\newcommand{\Z}{{\mathbb Z }}
\newcommand{\R}{{\mathbb R}}

\newcommand{\e}{{\varepsilon }}

\newcommand{\ie}{{\sl i.e.\/ }}
\newcommand{\cf}{{\sl cf.\/ }}
\newcommand{\eg}{{\sl e.g.\/}}

\newcommand{\D}{{\mathrm {d}}}

\renewcommand{\H}{{\mathbf H}}

\newcommand{\p}{\partial}

\newcommand{\f}[2]{\frac{#1}{#2}}

\newcommand{\Ld}{\Lambda}

\newcommand{\Gm}{\Gamma}
\newcommand{\gm}{\gamma}
\newcommand{\vp}{\varphi}

\newcommand{\ift}{\infty}

\newcommand{\fa}{\forall}

\newcommand{\tg}{\triangle}

\newcommand{\wt}{\widetilde}

\renewcommand{\theequation}{\arabic{section}.\arabic{equation}}
\newcommand{\be}{\begin{equation}}
\newcommand{\ee}{\end{equation}}
\newcommand{\ba}{\begin{array}}
\newcommand{\ea}{\end{array}}
\newcommand{\bea}{\begin{eqnarray}}
\newcommand{\eea}{\end{eqnarray}}
\newcommand{\beas}{\begin{eqnarray*}}
\newcommand{\eeas}{\end{eqnarray*}}
\newcommand{\dpm}{\displaystyle}

\def\({\left(}
\def\){\right)}
\def\<{\left\langle}
\def\>{\right\rangle}
\def\O{\mathcal O}
\newcommand{\newpar}{\par}\parindent =15pt\parskip=3pt\textheight = 615pt
\renewcommand{\L}{{L}}
\newcommand{\Id}[1]{{\rm I\kern-2pt I_{#1}}}
\renewcommand{\hbar}{{\displaystyle\bar{\phantom{x}}\kern-6pt h}}

\numberwithin{equation}{section}

\title
{A Bloch decomposition based split-step pseudo spectral method for
quantum dynamics with periodic potentials\thanks{This work was
partially supported by the Wittgenstein Award 2000 of P.~A.~M.,
NSF grant No.~DMS-0305080, the NSFC Projects No.~10301017 and
10228101, the National Basic Research Program of China under the
grant 2005CB321701, SRF for ROCS, SEM and the Austrian-Chinese
Technical-Scientific Cooperation Agreement. C.~S. has been
supported by the APART grant of the Austrian Academy of Science.}}

\author{Zhongyi Huang\thanks{Department of Mathematical Sciences,
Tsinghua University, Beijing 100084,
China, Phone: (+8610) 62796893, Fax: (+8610) 62773400,
({\tt zhuang@math.tsinghua.edu.cn})} \and
Shi Jin\thanks{Department of
Mathematics, University of Wisconsin, Madison, WI 53706, USA and
Department of Mathematical Sciences, Tsinghua University, Beijing
100084, China, Phone: (608)263-3302, Fax: (608)263-8891
({\tt jin@math.wisc.edu})}
\and Peter A.
Markowich\thanks{Wolfgang Pauli Institute Vienna \& Faculty of
Mathematics, University of Vienna, Nordbergstra\ss e 15, A-1090
Vienna, Austria, Phone: (+43) 1427750611, Fax: (+43) 142779506,
({\tt peter.markowich@univie.ac.at})} \and Christof
Sparber\thanks{Wolfgang Pauli Institute Vienna \& Faculty of
Mathematics, University of Vienna, Nordbergstra\ss e 15, A-1090
Vienna, Austria, Phone: (+43) 1427750716, Fax: (+43) 1427750650,
({\tt christof.sparber@univie.ac.at})}
}

\begin{document}

\maketitle

\begin{abstract}
We present a new numerical method for accurate computations of
solutions to (linear) one dimensional Schr\"odinger equations with
periodic potentials. This is a prominent model in solid state
physics where we also allow for perturbations by non-periodic
potentials describing external electric fields. Our approach is
based on the classical Bloch decomposition method which allows to
diagonalize the periodic part of the Hamiltonian operator. Hence,
the dominant effects from dispersion and periodic lattice
potential are computed together, while the non-periodic potential
acts only as a perturbation. Because the split-step communicator
error between the periodic and non-periodic parts is relatively
small, the step size can be chosen substantially larger than for
the traditional splitting of the dispersion and potential
operators. Indeed it is shown by the given examples, that our
method is unconditionally stable and more efficient than the
traditional split-step pseudo spectral schemes. To this end a
particular focus is on the semiclassical regime, where the new
algorithm naturally incorporates the adiabatic splitting of slow
and fast degrees of freedom.
\end{abstract}



\pagestyle{myheadings} \thispagestyle{plain} \markboth{Z. Huang,
S. Jin, P. A. Markowich, C. Sparber}{Bloch decomposition based
numerical method}


\section{Introduction}\label{sec:intro}

One of the main problems in solid state physics is to describe the
motion of electrons within the periodic potentials generated by
the ionic cores. This problem has been studied from a physical, as
well as from a mathematical point of view in, \eg, \cite{AsKn, Bl,
Lu, PST, Za}, resulting in a profound theoretical understanding of
the novel dynamical features. Indeed one of the most striking
effect, known as \emph{Peirl's substitution}, is a modification of
the dispersion relation for Schr\"odinger's equation, where the
classical energy relation $E_{\rm free}(k) = \frac{1}{2}|k|^2$ has
to be replaced by the $E_m(k)$, $m\in \N$, the  energy
corresponding to the $m$th \emph{Bloch band} \cite{B}. The basic
idea behind this replacement is a separation of scales which is
present in this context. More precisely one recognizes that
experimentally imposed, and thus called external, electromagnetic
fields typically vary on much \emph{larger} spatial scales than
the periodic potential generated by the cores. Moreover this
external fields can be considered \emph{weak} in comparison to the
periodic fields of the cores \cite{AsMe}.

To study this problem, consider the  Schr\"odinger equation for
the electrons in a \emph{semiclassical} asymptotic scaling
\cite{CMS, PST, Te}, \ie in $d=1$ dimensions
\begin{equation}\label{s0}
\left \{
\begin{aligned}
& \I \e   \partial_t  \psi   = -\f{\e^2}{2} \, \p_{xx}  \psi+
V_\Gm\left(\f{x}{\e}\right)\psi+U(x)\psi ,\qquad x \in \R, \ t\in \R,\\
& \psi\big |_{t=0}=\psi_{\rm in}(x),
\end{aligned}
\right.
\end{equation}
where $0<\e \ll 1$, denotes the small \emph{semiclassical
parameter} describing the microscopic/macroscopic scale ratio. The
(dimensionless) equation \eqref{s0} consequently describes the
motion of the electrons on the macroscopic scales induced by the
external potential $U(x)\in \R$. The highly oscillating
\emph{lattice-potential} $V_\Gm(y)\in \R$ is assumed to be
\emph{periodic} with respect to some \emph{regular lattice} $\Gm
$. For definiteness we shall assume that \be\label{lp0}
V_\Gm(y+2\pi )=V_\Gm(y) \quad \fa y\in \R, \quad \ee \ie $\Gamma =
2 \pi \Z$. In the following we shall assume $\psi_{\rm in}\in
L^2(\R)$, such that the \emph{total mass} is $\mathrm M_{\rm
in}\equiv {\| \, \psi_{\rm in} \|}_{L^2}=1$, a normalization which
is henceforth preserved by the evolution.

The mathematically precise asymptotic description of $\psi(t)$,
solution to \eqref{s0}, as $\e \to 0$, has been intensively
studied in, \eg, \cite{BLP, GMMP, GRT, PST}, relying on different
analytical tools. On the other hand the numerical literature on
these issues is not so abundant \cite{Go, GoMa, GoMau}. Here we
shall present a novel approach to  the numerical treatment of
\eqref{s0} relying on the classical \emph{Bloch decomposition}
method, as explained in more detail below. The main idea is to
treat in one step the purely dispersive part $\propto
\partial_{xx}$ of the Schr\"odinger equation together with the
periodic potential $V_\Gamma$, since this combined operator allows
for some sort of ``diagonalization'' via the Bloch transformation.
The corresponding numerics is mainly concerned with the case $\e
\ll 1$ but we shall also show examples for a rather large $\e
=\frac{1}{2}$. Our numerical experiments show that the new method
converges with $\Delta x=\O(\e)$ and $\Delta t=\O(1)$, the latter
being a huge advantage in comparison with a more standard
time-splitting method used in \cite{Go, GoMa, GoMau}, and which
usually requires $\Delta t=\O(\e)$. Moreover we find that the use
of only a few Bloch bands is mostly enough to achieve very high
accuracy, even in cases where $U(x)$ is no longer smooth. We note
that our method is unconditionally stable and comprises spectral
convergence for the space discretization as well as second order
convergence in time. The only drawback of the method is that we
first have to compute the energy bands for a given periodic
potential, although this is needed only in a preprocessing step
rather than during the time marching. On the other hand, this
preprocessing also handles a possible lack of regularity in
$V_\Gamma$, which consequently does not lead to numerical problems
during the time-evolution. In any case the numerical cost of this
preliminary step is much smaller than the costs spend in computing
the time-evolution and this holds true for whatever method we
choose.

We remark that linear and nonlinear evolutionary PDEs with
periodic coefficients also arise in the study of photonic
crystals, laser optics, and Bose-Einstein condensates in optical
lattices, \cf \cite{Bu, CMS, HFKW} and the references given
therein. We expect that our algorithm can adapted to these kind of
problems too. Also note, that in the case of a so-called
\emph{stratified medium}, see, \eg, \cite{BLP, BePo}, an
adaptation of our code to higher dimensions is very likely.
Finally, the use of the Bloch transformation in problems of
\emph{homogenization} has been discussed in \cite{COV, CoVa} and
numerically studied in \cite{CSV} for elliptic problems. Our
algorithm might be useful in similar time-dependent numerical
homogenization problems.
\newpar
The paper is organized as follows: In Section \ref{BD}, we recall
in detail the Bloch-decomposition method and we show how to
numerically calculate the corresponding energy bands. Then, in
Section \ref{BDTSS} we present our new algorithm, as well as the
usual time-splitting spectral method for Schr\"odinger equations.
In section \ref{NUM}, we show several numerical experiments, and
compare both methods. Different examples of $U$ and $V_\Gamma$ are
considered, including the non-smooth cases. Finally we shall also
study a WKB type semiclassical approximation in Section
\ref{sec:asym} and compare its numerical solution to solution of
the full problem. This section is mainly included since it gives a
more transparent description of the Bloch transformation, at least
in cases where a semiclassical approximation is justified.

\section{The emergence of Bloch bands }\label{BD}

First, let us introduce some notation used throughout this paper,
respectively recall some basic definitions used when dealing with
periodic Schr\"odinger operators \cite{AsMe, BLP, Te, Wi}.

With $V_\Gamma$ obeying \eqref{lp0} we have:
\begin{itemize}
    \item The fundamental domain of our lattice $\Gamma = 2\pi \mathbb Z$,
is ${\mathcal C} =(0,2\pi)$.
    \item The \emph{dual lattice} $\Gm^*$ can then be defined as the set of
all wave numbers $k \in \R$, for which plane waves of the form
$\exp(\I k x)$ have the same periodicity as the potential $V_\Gm$.
This yields $\Gm^*=\Z$ in our case.
    \item The fundamental domain of the dual lattice, \ie the (first)
\emph{Brillouin zone}, $\mathcal B=\mathcal C^*$ is the set of all
$k\in \R$ closer to zero than to any other dual lattice point.
In our case, that is $\mathcal B=\left(-\f{1}{2},\f{1}{2}\right)$.
\end{itemize}
\subsection{Recapitulation of Bloch's decomposition method}\label{bloch}

One of our main points in all what follows is that the dynamical
behavior of \eqref{s0} is mainly governed by the periodic part of
the Hamiltonian, in particular for $\e \ll 1$. Thus it will be
important to study its spectral properties. To this end consider
the periodic \emph{Hamiltonian} (where for the moment we set $y =
x/ \e $ for simplicity) \be\label{Hap} H=  -\f{1}{2}\,
\p_{yy}+V_\Gm\left(y \right), \ee which we will regard here only
on $L^2(\mathcal C)$. This is possible since due to the
periodicity of $V_\Gamma$ which allows to then to cover all of
$\R$ by simple translations. More precisely, for $k \in
\overline{{\mathcal B}}=\left[-\f{1}{2},\f{1}{2}\right]$ we equip
the operator $H$ with the following \emph{quasi-periodic} boundary
conditions \be\label{perpsi1} \left\{
\begin{aligned}
\psi(t,y+2\pi)= & \, \ex^{ 2 \I k\pi}\psi(t,y) \quad \fa \, y\in \R,
\ k \in \overline{\mathcal B},\\
\p_y\psi(t,y+2\pi)= & \, \ex^{2 \I k\pi}\p_y\psi(t,y) \quad \fa \,
y\in \R,\ k\in \overline{\mathcal B}.
\end{aligned}
\right. \ee It is well known \cite{Wi} that under very mild
conditions on $V_\Gamma$, the operator $H$ admits a complete set
of eigenfunctions $\vp_m(y,k), m\in \N$, providing, for each fixed
$k \in \overline{{\mathcal B}}$, an orthonormal basis in
$L^2(\mathcal C)$. Correspondingly there exists a countable family
of real-valued eigenvalues which can be ordered according to
$E_1(k)\le E_2(k)\leq \cdots \le E_m(k)\le \cdots$, $ m\in \N,$
including the respective multiplicity. The set $\label{Em}
\{E_m(k)\,|\, k\in \mathcal B\}\subset \R$ is called the $m$th
\emph{energy band} of the operator $H$ and the eigenfunctions
$\vp_m(\cdot,k)$ is usually called \emph{Bloch function}. (In the
following the index $m \in \N$ will \emph{always} denote the
\emph{band index}.) Concerning the dependence on $k\in \mathcal
B$, it has been shown \cite{Wi} that for any $m\in \N$ there
exists a closed subset $\mathcal A\subset \mathcal B$ such that:
$E_m(k)$ is analytic and $\vp_m(\cdot, k)$ can be chosen to be
real analytic function for all $k \in \overline{\mathcal B}
\backslash \mathcal A$. Moreover
\begin{equation}
\label{iso} E_{m-1} < E_m(k) < E_{m+1}(k)\quad \forall \, k \in
\overline{\mathcal B} \backslash \mathcal A.
\end{equation}
If this condition indeed holds for all $k\in \mathcal B$ then
$E_m(k)$ is called an \emph{isolated Bloch band} \cite{Te}.
Moreover, it is known that \be \meas \mathcal A = \meas \, \{ k\in
\overline{\mathcal B} \ | \ E_n(k)=E_{m}(k), \ n\not = m \}=0. \ee
In this set of measure zero one encounters so called \emph{band
crossings}. Note that due to \eqref{perpsi1} we can rewrite
$\vp_m(y,k)$ as \be \vp_m(y,k) = \ex^{\I k y} \chi_m(y,k)\quad
\forall \, m \in \N, \ee for some $2\pi$-periodic function
$\chi_m(\cdot,k)$. In terms of $\chi_m(y,k)$ the \emph{Bloch
eigenvalue problem} reads \be \label{chi1} \left \{
\begin{aligned}
H(k) \chi_m(y,k)= & \, E_m(k)\chi_m(y,k),\\
\chi_m(y+2\pi ,k)=& \, \chi_m(y,k) \quad \fa  \, k\in \mathcal B,
\end{aligned}
\right. \ee where $H(k)$ denotes the shifted Hamiltonian \be
H(k):= \f{1}{2}(-\I\p_y+k)^2+V_\Gm(y). \ee

Let us know introduce the so-called \emph{Bloch transform}
$\mathcal T$ of some function $\psi(t,\cdot) \in L^2(\R)$, for any
fixed $t\in \R$, as can be found in, \eg, \cite{PST, Te}. (Some
other variants of this transformation can also be found in the
literature.) The Bloch transformation $\mathcal T$ is just the
regular Fourier transform $\mathcal F$ on the factor $\ell
^2(\Gamma)$ followed by a multiplication with ${\rm e}^{-\I y
k}$, \ie \be \label{btr} (\mathcal T \psi)(t,k,y):= \sum_{\gamma
\in \Z}\psi(t, y+ 2 \pi \gamma)\, \ex^{-\I k ( 2\pi \gamma +y) }
,\quad
 y\in  {\mathcal C},\ k\in   {\mathcal B}.
\ee It is then easy to see that \be \mathcal T H {\mathcal T}^{-1}
= H(k). \ee which provides a link between the eigenvalue problem
\eqref{chi1} and the periodic part of our Schr\"odinger equation
acting on $\psi(t,\cdot)$.

Most importantly though the Bloch transformation allows to
decompose our original Hilbert space $\mathcal H=L^2(\R)$ into a
direct sum of, so called, \emph{band spaces}, \ie \be L^2(\R)=
\bigoplus_{m=1}^\infty \mathcal H_m, \quad \mathcal H_m:=\left
\{\,  \psi_m(t,y)= \int_\mathcal B f(t,k) \, \vp_m(y,k) \, \D k, \
f(t, \cdot)\in L^2(\mathcal B) \right\}, \ee for any fixed $t\in
\R$. This is the well known \emph{Bloch decomposition method},
which implies that \be\label{sum} \forall \, \psi(t, \cdot)\in
L^2(\R):\quad \psi(t,y)=\sum_{m\in \N} \psi_m(t,y), \quad \psi_m
\in \mathcal H_m. \ee The corresponding projection of $\psi(t)$
onto the $m$th band space is thereby given as
\be\label{projection} \psi_m(t,y)\equiv (\mathbb P_m \psi)(t,y)=
\int_{\mathcal B} \left(\int_{\R} \psi(t,\zeta)
\overline{\vp}_m\left(\zeta ,k\right) \D \zeta \right)
\vp_m\left(y,k\right) \D k \ee and we consequently denote by \be
\label{coeff} C_m(t,k):=\int_{\R} \psi(t,\zeta)
\overline{\vp}_m\left(\zeta ,k\right) \D \zeta \ee the
coefficients of the Bloch decomposition. For a complete
description and a rigorous mathematical proof of this
decomposition we refer to, \eg, \cite{ReSi}, chapter XI. Here it
is only important to note that the Bloch transformation allows to
obtain a spectral decomposition of our periodic Hamiltonians $H$,
upon solving the eigenvalue problem \eqref{chi1}. Roughly speaking
$\mathcal T$ can be seen as some sort of Fourier transform adapted
to the inclusion of periodic coefficients (potentials).

This consequently implies that, if $U\equiv 0$, we can indeed
Bloch transform the whole evolution problem \eqref{s0} and
decompose it into the corresponding band spaces $\mathcal H_m$,
\ie we gain some sort of ``diagonalization'' for our evolution
problem. In this case each $\psi_m(t,\cdot)\in \mathcal H_m$ then
evolves according to the newly obtained PDE
\begin{equation}\label{evotr}
\left \{
\begin{aligned}
& \I \e \partial_t \psi_m = E_m(-\I \partial_y) \psi_m , \qquad
y \in \R, \ t\in \R,\\
& \psi_m\big |_{t=0}= (\mathbb P_m \psi_{\rm in})(y).
\end{aligned}
\right.
\end{equation}
Here $E_m(-\I \partial_y)$ denotes the pseudo-differential
operator corresponding to the (Fourier-) symbol $E_m(k)$, \cf
\cite{GMMP, PST, Te}. The above given evolution equation comprises
a rigorous justification of Peirl's substitution. Moreover
\eqref{evotr} is easily solved invoking the standard Fourier
transformation $\mathcal F$ on $L^2(\R)$, which yields \be
\label{blsol} \psi_m (t,y) = \mathcal F^{-1}\left(\ex^{- \I E_m(k)
t /\e} (\mathcal F (\mathbb P^\e_m \psi_{\rm in}))(k)\right). \ee
Here the energy band $E_m(k)$ is understood to be periodically
extended on all of $\R$. To this end, note that the following
relation holds \be \label{Fcoeff} \mathcal F (\psi_{m})(t,k) =
\ex^{- \I E_m(k) t /\e}C_m(0,k) (\mathcal F\chi_m)(0,k), \ee as
can be shown by a lengthy but straightforward calculation.

Of course if $U \not \equiv 0$ (the non-periodic part of the
potential) the time evolution \eqref{s0} in general \emph{mixes}
all band spaces $\mathcal H_m$, \ie we can no longer hope to be
able to diagonalize the whole Hamiltonian operator (which now
involves also non-periodic coefficients). On the other hand, since
$U(x) = U(\e y)$ varies only slowly on the fast (periodic) scale
$y = x/\e$, one might hope that even if $U\not \equiv 0$, the
\emph{effective Schr\"odinger type equation}
\begin{equation}\label{eff}
\left \{
\begin{aligned}
& \I \e \partial_t \psi^{\rm eff}_m = \, E_m(-\I \partial_y)
\psi^{\rm eff}_m + U(\e y ) \psi^{\rm eff}_m, \qquad y \in \R,
\ t\in \R,\\
& \psi^{\rm eff}_m\big |_{t=0}= (\mathbb P_m \psi_{\rm in})(y),
\end{aligned}
\right.
\end{equation}
holds true, at least approximately for small $\e\ll 1$. In other
words, we expect the slowly varying external potential to be
almost constant on the lattice scale and thus yielding only a
small perturbation of the band structure determined via
\eqref{Hap}. Indeed this is the case as has been rigorously proved
in \cite{CMS, GRT, PST}, using different analytical approaches,
(for a broader overview, see \cite{Te} and the references given
therein), where it is shown that \be \sup_{t\in I} {\big
\|(\mathbb P_m \psi)(t) - \psi^{\rm eff}_m (t) \big \|}_{\L^2(\R)}
\leq \O(\e), \ee holds true for any finite time-interval $I\subset
\R$. Here $\psi(t)$ is the solution of the full Schr\"odinger
equation and $\psi^{\rm eff}_m(t)$ is the solution of the
effective model \eqref{eff}. To this end one hast to assume that
the m'th energy band is \emph{isolated} from the rest of the
spectrum though. If this is not the case, energy transfer of order
$\O(1)$ can occur at band crossings, the so-called Landau-Zener
phenomena.

\subsection{Numerical computation of the Bloch bands}\label{numbands}

As a preparatory step for our algorithm we shall first calculate
Bloch's energy bands $E_m(k)$ numerically as follows. Analogously
to  \cite{GoMa, Ko}, we consider the potential $V_\Gm\in C^1(\R)$
and expand it in its Fourier series, \ie \be\label{Vgm}
V_\Gm(y)=\sum_{\lambda \in\Z}\widehat V(\lambda) \, \ex^{\I
\lambda y},\quad \widehat V(\lambda)
=\f{1}{2\pi}\int_0^{2\pi}V_\Gm(y)\, \ex^{-\I \lambda y }\, \D y.
\ee Likewise, we expand any Bloch eigenfunctions
$\chi_m(\cdot,k)$, in its respective Fourier series
\be\label{chiFourier} \chi_m(y,k)=\sum_{\lambda \in\Z}\widehat
\chi_m(\lambda, k)\, \ex^{\I \lambda y},\quad \widehat
\chi_m(\lambda, k)=\f{1}{2\pi}\int_0^{2\pi}\chi_m(y,k)\, \ex^{-\I
\lambda y}\, \D y. \ee (The latter should \emph{not} be confused
with the so-called \emph{Wannier functions} which are given as the
Fourier transformation of $\varphi_m$ w.r.t to $k \in \mathcal
B$.) Clearly the Fourier approximation of $V_\Gamma$, and thus
also the one of $\chi_m$, depends on the regularity of $V_\Gamma$.
If $V_\Gamma \in C^\infty(\R)$ the corresponding Fourier
coefficients $\widehat V(\lambda)$ decay faster than any power, as
$\lambda \to \pm \infty$, and thus we only need to take into
account a few coefficients in this case.

For $\lambda \in \{-\Lambda,\cdots,\Lambda-1 \}\subset \mathbb Z$,
we consequently aim to approximate the Sturm-Liouville problem
\eqref{chi1}, by the following algebraic eigenvalue problem
\bea\label{chieig}
\H(k)\left(%
\begin{array}{c}
  \widehat \chi_m(-\Lambda) \\
  \widehat \chi_m(1-\Lambda) \\
  \vdots \\
  \widehat \chi_m(\Lambda-1) \\
\end{array}%
\right) =E_m(k)
\left(%
\begin{array}{c}
  \widehat \chi_m(-\Lambda) \\
  \widehat \chi_m(1-\Lambda) \\
  \vdots \\
  \widehat \chi_m(\Lambda-1) \\
\end{array}%
\right) \eea where the $2\Lambda\times 2\Lambda$ matrix $\H(k)$ is
given by \be \label{Hk}
\H(k)=\left(%
\begin{array}{cccc}
  \widehat V(0)+\f{1}{2}(k-\Lambda)^2 & \widehat V(-1) & \cdots &
  \widehat V(1-2\Lambda) \\
  \widehat V(1) & \widehat V(0)+\f{1}{2}(k-\Lambda+1)^2  & \cdots &
  \widehat V(2-2\Lambda) \\
  \vdots & \vdots & \ddots & \vdots \\
  \widehat V(2\Lambda-1) & \widehat V(2\Lambda-2) & \cdots & \widehat
  V(0)+\f{1}{2}(k+\Lambda-1)^2  \\
\end{array}%
\right) \ee The above given matrix $\H(k)$ comprises $2\Lambda$
eigenvalues. Clearly, this number has to be large enough such that
all the eigenvalues $E_m(k)$ which we need to use in our
simulations below are counted, \ie we need $m \leq 2\Lambda$. The
numerical cost for this algebraic problem is about
$\O(\Lambda^3)$, \cf \cite{HJ}. Note however that this is the most
expensive case, which becomes considerably smaller if one exploits
possible symmetries within the potential $V_\Gamma$, \cf Example
\ref{nuex1} below (see also \cite{Bu, Ko, HFKW, Ze}). In any case
the number $\Lambda$ is \emph{independent} of the spatial grid,
thus the numerical costs of this eigenvalue problem are almost
negligible compared to those spend in the evolutionary algorithms
below. The approximate numerical computations of the Bloch bands
$E_m(k)$ can be seen as a preprocessing, to be done only once and
remain unchanged as time evolves.

\begin{remark}\label{rem}
Accurate computations of the energy bands needed in practical
applications, \ie in more than one spatial dimensions and for
different kind of (composite) material, becomes a highly
nontrivial task. Nowadays though, there already exists a huge
amount of numerical data comprising the energy band structure of
the most important materials used in, \eg, the design of
semiconductor devices, \cf \cite{FL, JC, LVF}. We note that some
of these data is available online via the URL {\tt
http://www.research.ibm.com/DAMOCLES/home.html}, or {\tt
http://cmt.dur.ac.uk/sjc}, and also {\tt
http://cms.mpi.univie.ac.at/vasp/vasp/vasp.html}. In the context
of photonic crystals the situation is similar \cite{HFKW}. Thus,
relying on such data one can in principle avoid the above given
eigenvalue-computations (and its generalizations to more
dimensions) completely. To this end, one should also note that,
given the energy bands $E_m(k)$, we do not need any knowledge
about $V_\Gamma$ in order to solve \eqref{s0} numerically, \cf the
algorithm described below.
\end{remark}

\section{Bloch decomposition based algorithm vs. time-splitting
spectral methods}\label{BDTSS}
For the convenience of computations, we shall consider the
equation (\ref{s0}) on a bounded domain $\mathcal D$, say on the
interval $\mathcal D= [-\kappa_1, \kappa_2 ]$, for some large
enough $\kappa_1, \kappa_2 >0$. Moreover we shall equip $\mathcal
D$ with \emph{periodic boundary conditions}. However, this
periodic computational domain $\mathcal D$ should \emph{not} be
confused with the periodic structure induced by the lattice
potential. Without loss of any generality, we assume that
$\mathcal D=[0,2\pi]$.

For practical reasons we shall now introduce, for any fixed $t\in
\R$, a new unitary transformation of $\psi(t, \cdot)\in L^2(\R)$
\be\label{blochtrans} \wt\psi(t,y,k):=\sum_{\gamma \in \Z}\psi(t,
\e(y+ 2 \pi \gamma))\, \ex^{-\I 2\pi k \gamma},\quad
 y\in  {\mathcal C},\ k\in   {\mathcal B},
\ee which has the properties that $\wt \psi$ is quasi-periodic
w.r.t $y \in \Gamma$ and periodic w.r.t. $k \in \Gamma^*$, \ie
\be\label{trans-prop} \wt\psi(t,y+2\pi,k)=  \ex^{\I2 \pi k}\,
\wt\psi(t,y,k), \quad \wt\psi(t,y,k+1 )=  \wt\psi(t,y,k). \ee One
should note that $\wt \psi$ is \emph{not} the standard \emph{Bloch
transformation} $\mathcal T $, as defined in \eqref{btr}, but it
is indeed closely related to it via \be (\mathcal T \psi)(t,y,k) =
\wt\psi(t,y,k) \ex^{- \I y k}, \quad k \in \mathcal B, \ee for
$\e=1$. Furthermore, we have the following inversion formula
\be\label{inverse-BT} \psi(t,\e(y+2\pi \gm))=\int_{\mathcal B}
\wt\psi(t,y,k) e^{\I2\pi k \gm} dk, \ee which is again very
similar to the one of the standard Bloch transformation \cite{Te}.
The main advantage in using $\wt \psi$, instead of $\mathcal T
\psi$ itself, is that we can rely on a standard fast Fourier
transform  (FFT) in the numerical algorithm below. If one aims to
use $\mathcal T \psi$ directly one would be forced to modify a
given FFT code accordingly. A straightforward computation then
shows that \be \label{coeff1} C_m(t,k)= \int_{\mathcal C} \wt
\psi(t,\zeta,k) \overline{\vp}_m\left(\zeta ,k\right) \D \zeta,
\ee where $C_m(t,k)$ is the Bloch coefficient, defined in
\eqref{coeff}.

In what follows, let the time step be $\Delta t=T/N$, for some
$N\in \N$, $T>0$. Suppose that there are $L\in \N$ lattice cells
within the computational domain $\mathcal D=[0,2\pi]$. In this
domain, the wave function $\psi$ is numerically computed at $L
\times R$ grid points, for some $R \in \N$. In other words we
assume that there are $R$ grid points in each lattice cell, which
yields the following discretization \be \left \{
\begin{aligned}
k_\ell=& \, -\f{1}{2}+\f{\ell-1}{L},\quad  \mbox{where
$\ell=\{1,\cdots,L\}\subset \N$},\\
y_r=&\, \f{2\pi (r-1)}{R},\quad \quad  \ \mbox{where
$r=\{1,\cdots, R\}\subset \N$}, \\
\end{aligned}
\right. \ee and thus we finally we evaluate $\psi^n = \psi(t_n)$
at the grid points $x = \e(2\pi \gm + y)$, \ie \be x_{\ell,r}= \,
\e(2\pi (\ell - 1)+y_r). \ee We remark that in our numerical
computations we can use $R \ll L$, whenever $\e \ll 1$, \ie we
only use a few grid points within each cell. Now we shall describe
precisely the Bloch decomposition based algorithm used to solve
\eqref{s0}.

\subsection{The Bloch decomposition based algorithm (BD)}

Suppose that at the time $t_n$ we are given
$\psi(t_n,x_{\ell,r})\approx \psi^n_{\ell,r}$. Then
$\psi^{n+1}_{\ell,r}$, \ie the solution at the (next) time step
$t_{n+1}=t_n+\Delta t$, is obtained as follows:

\textbf{Step 1.} First, we solve the equation
\begin{equation}
\label{bd1}
\begin{aligned}
&\, \I \e \partial _t\psi = -\frac{\e^2}{2}\, \p_{xx}\psi
+V_\Gm\left(\f{x}{\e}\right)\psi, \\
\end{aligned}
\end{equation}
on a fixed time-interval $\Delta t$. To this end we shall heavily
use the Bloch-decomposition method, see below.

\newpar
\textbf{Step 2.} In a second step, solve the ordinary differential
equation (ODE)
\begin{equation}
\label{bd2}
\begin{aligned}
& \,\I \e \partial _t\psi  = U(x) \psi, \\
\end{aligned}
\end{equation}
on the same time-interval, where the solution obtained in Step 1
serves as initial condition for Step 2. We easily obtain the exact
solution for this linear ODE by \be \label{bd3}
\psi(t,x)=\psi(0,x)\,\ex^{-\I U(x) t/\e}. \ee
\begin{remark} \label{rem1}
Clearly, the algorithm given above is first order in time. But we
could easily obtain also a second order scheme by the Strang
splitting method, which means that we use Step 1 with time-step
$\tg t/2$, then Step 2 with time-step $\tg t$, and finally
integrate Step 1 again with $\tg t/2$. Note that in both cases the
scheme conserves the particle density $\rho(t,x):= |\psi(t,x)|^2$
, also on the fully discrete level.
\end{remark}
Indeed Step 1 consists of several intermediate steps which we
shall present in what follows:

\textbf{Step 1.1.} We first compute $\wt \psi$ at time $t^n$ by
\be\label{blochst1}
\wt\psi^n_{\ell,r}=\sum_{j=1}^{L}\psi^{n}_{j,r}\, \ex^{-\I
k_\ell\cdot x_{j,1}}. \ee
\newpar
\textbf{Step 1.2.} Next, we compute the $m$th band Bloch
coefficient $C_m(t,k)$, at time $t^n$, via \eqref{coeff1}, \ie
\be\label{blochst2}
\begin{aligned}
C_m(t_n,k_\ell)\approx C^n_{m,\ell}= &\,\frac{2\pi}{R}\,
\sum_{r=1}^{R}\wt\psi^n_{\ell,r}\overline{{\chi}_{m}}(y_r,k_\ell)\,
\ex^{-\I  k_\ell  y_r} \\\approx &\, \frac{2\pi}{R}\,
\sum_{r=1}^{R}\wt\psi^n_{\ell,r}\sum_{\lambda=-R/2}^{R/2-1}
\overline{\widehat{\chi}_{m}}(\lambda,k_\ell)\,
\ex^{-\I  (k_\ell+\lambda)  y_r} ,
\end{aligned}
\ee where for the second line we simply inserted the Fourier
expansion of $\chi_m$, given in \eqref{chiFourier}. Note that in
total we have $R$ Fourier coefficients for $\chi_m$. Clearly this
implies that we need $\Ld > R/2$ to hold, where $\Lambda$ is the
number of Fourier modes required in the numerical approximation of
Bloch's eigenvalue problem as discussed in Section \ref{numbands}.
Here we only take the $R$ lowest frequency Fourier coefficients.
\newpar
\textbf{Step 1.3.} The obtained Bloch coefficients are then
evolved up to the time $t^{n+1}$, according to the explicit
solution formula \eqref{blsol}, taking into account
\eqref{Fcoeff}. This yields \be\label{blochst3}
C^{n+1}_{m,\ell}=C^n_{m,\ell}\, \ex^{-\I E_m(k_\ell)\Delta t/\e}.
\ee
\newpar
\textbf{Step 1.4.} From here, we consequently compute $\wt \psi$
at the new time $t^{n+1}$ by summing up all band contributions and
using the analytical formulas \eqref{projection} and
\eqref{coeff}, \ie \be\label{blochst4} \wt\psi^{n+1}_{\ell,r}=
\sum_{m=1}^M (\mathbb P_m \wt \psi)_{\ell, r}^{n+1} \approx
\sum_{m=1}^MC^{n+1}_{m,\ell}\sum_{\lambda=-R/2}^{R/2-1}
\widehat{\chi}_{m}(\lambda,
k_\ell)\, \ex^{\I  (k_\ell+\lambda)  y_r}. \ee
\newpar
\textbf{Step 1.5.} Finally we numerically perform the inverse
transformation to \eqref{blochtrans}, \ie we compute
$\psi^{n+1}_{\ell,r}$ from $\wt \psi^{n+1}_{\ell,r}$. Thus from
\eqref{inverse-BT}, we get \be\label{blochst5}
\psi^{n+1}_{\ell,r}=\f{1}{L}\sum_{j=1}^{L}\wt\psi^{n+1}_{j,r}\,
\ex^{\I k_j  x_{\ell,1}}. \ee

Note that in the BD algorithm, the main numerical costs are
introduced via the FFT in Steps 1.1 and 1.5. This also implies
that on the same spatial grid, the numerical costs of our Bloch
transform based algorithm is of the same order as the classical
time-splitting spectral method below. Moreover, we want to stress
the fact that if there is no external potential, \ie $U(x)\equiv
0$, then the above given algorithm numerically computes the
\emph{exact} solution of the evolutionary problem \eqref{s0},
which can be seen analogous to a standard spectral method, adapted
to periodic potentials. In particular this fact allows us to solve
the Schr\"odinger equation \eqref{s0} for very long time steps,
even if $\e$ is small (see the results given below). Moreover, one
should note that a possible lack of regularity in $V_\Gamma$ only
requires numerical care when approximating \eqref{chi1} by the
algebraic problem \eqref{chieig}. In particular, $V_\Gamma$ itself
does not enter in the time-evolution but only $E_m(k)$.

\subsection{A simple time-splitting spectral method (TS)}\label{tsm}
Ignoring for a moment the additional structure provided by the
periodic potential $V_\Gamma$, one might wish to solve \eqref{s0}
by using a classical time-splitting spectral scheme. Such schemes
already proved to be successful in similar circumstances, see,
\eg, \cite{BJM, BJM1, GoMa, HJMSZ}. For the purpose of a detailed
comparison, we present this method here:

\textbf{Step 1.} In the first step we solve the equation
\begin{equation}
\label{st1}
\begin{aligned}
&\, \I \e \partial _t\psi = -\frac{\e^2}{2}\, \p_{xx}\psi, \\
\end{aligned}
\end{equation}
on a fixed time interval $\Delta t$, relying on the
pseudo-spectral method.
\newpar
\textbf{Step 2.} Then, in a second step, we solve the ordinary
differential equation
\begin{equation}
\label{st2}
\begin{aligned}
& \,\I \e \partial _t\psi  = \left(V_\Gm \left(\f{x}{\e}\right)
+ U(x) \right)\psi, \\
\end{aligned}
\end{equation}
on the same time-interval, where the solution obtained in step 1
serves as initial condition for step 2. Again it is easily seen,
that such a scheme conserves the particle density. It is clear
however that, due to the inclusion of
$V_\Gamma\left(\f{x}{\e}\right)$, the exact solution of
\eqref{st2} \be \psi(t,x)=\psi_{\rm in}(x)\,\ex^{-\I \left(V_\Gm
(x/\e)+ U(x) \right) t/\e}, \ee involves high oscillations on
\emph{different} length- and time-scales as $\e\to 0$ (which one
has to resolve), in contrast to \eqref{bd3}, where only
$t/\e$-oscillations are present.

\begin{remark}\label{rmk_comparison}
In our BD algorithm, we compute the dominant effects from
dispersion and periodic lattice potential in one step, and treat
the non-periodic potential as a perturbation. Because the
split-step communicator error between the periodic and
non-periodic parts is relatively small, the step size can be
chosen considerably larger than for the SP algorithm.
\end{remark}

\begin{remark}\label{rmk_no_lattice}
Clearly, if there is no lattice potential, \ie $V_\Gm(y)\equiv 0$,
the BD algorithm simplifies to the described time-splitting method
TS. Moreover, a second order second order scheme (based on the
Strang splitting algorithm) can be analogously obtained to the one
described above, see Remark \ref{rem1}, and a comparison of these
second order schemes gives similar results as those shown in the
following.
\end{remark}

\begin{remark}\label{rmk_complexity}
For the BD algorithm, the complexities of Step 1.1 and 1.5 are
$\O(RL \log (L))$, the complexities of Step 1.2 and 1.4 are
$\O(MLR\log(R))$, and for Step 1.3 we have $\O(ML)$. Also the
complexity of the eigenvalue problem \eqref{chieig} is
$\O(\Ld^3)$. However, since $\Ld$ (or $R$) is independent of $\e$
and since we only need to solve the eigenvalue problem
\eqref{chieig} once in a preparatory step, the computation costs
for this problem are negligible. On the other hand, for the TS
algorithm, the complexities of Step 1 and 2 are $\O(RL \log (RL))$
and $\O(RL)$ respectively. As $M$ and $R$ are independent of $\e$,
we can use $R \ll L$ and $M\ll L$, whenever $\e \ll 1$. Finally
the complexities of the BD and TS algorithm in each time step are
comparable.
\end{remark}
\section {Numerical experiments}\label{NUM}
\begin{figure} 
\begin{center}\footnotesize
\resizebox{4.8in}{2.2in}{\includegraphics{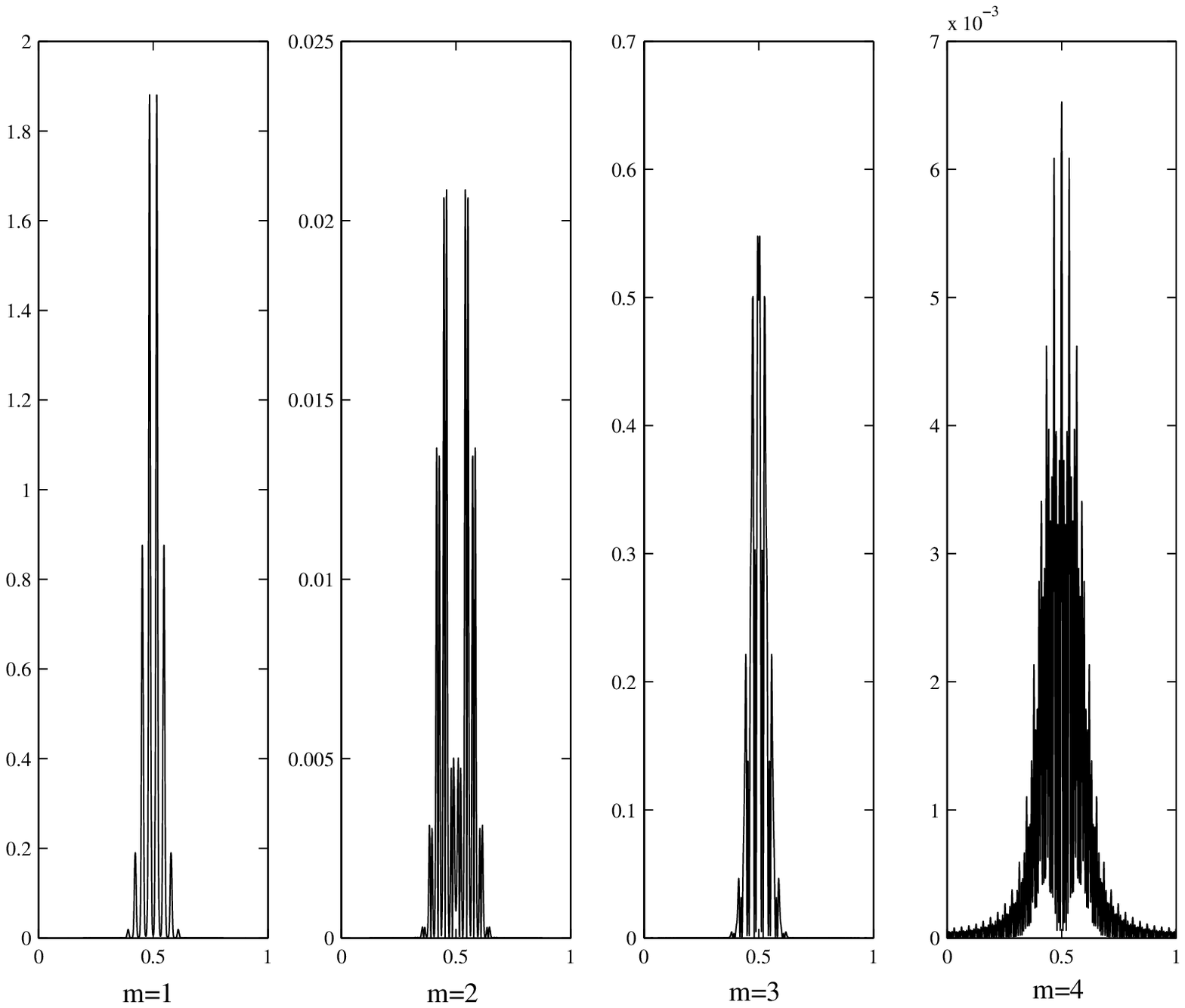}}\vspace{-2mm}
\end{center}
\caption{$|\mathbb P^\e_m \psi_{\rm in}|^2$, $m=1,\cdots,4$ for
$\e=\f{1}{32}$.} \label{fig:energy}\vspace{2mm}
\end{figure}

In this section, we shall use several numerical examples to show
the efficiency of our algorithm. We shall choose for \eqref{s0}
initial data $\psi_{\rm in} \in \mathcal S(\R)$ of the following
form \be\label{ex0} \psi_{\rm in}(x)=
\left(\f{10}{\pi}\right)^{1/4}\ex^{-5(x-\pi)^2}, \ee Let us
perform a decomposition of $\psi_{\rm in}$ in terms of the Bloch
bands, and take a summation of the first $m= 1,\cdots, M_0$ energy
bands, for some finite (cut-off) number $M_0\in \N$. A picture of
the corresponding band densities $\rho_m^\e:= |\mathbb P^\e_m
\psi_{\rm in}|^2$ is given in Figure \ref{fig:energy}, for
$m=1,\cdots,4$. Here $(\mathbb P^\e_m \psi_{\rm in})(x)$ is the
$\e$-scaled projection onto $\mathcal H_m^\e$, obtained from
\eqref{projection} by replacing $\vp_m(y,k)\to \e^{-1/2} \vp_m(x/
\e,k)$.
Since $\psi_{\rm in}$ is smooth we expect that only very few bands
have to be taken into account in the Bloch decomposition. Indeed
we observe that the amount of mass corresponding to $\mathbb
P^\e_m \psi_{\rm in}$, \ie the mass concentration in each Bloch
band, decays rapidly as $m \to \infty$, see Table
\ref{tab_energy}. In other words, the number $M_0$ is essentially
determined by the regularity of $\psi_{\rm in}$ in each cell. Note
that $M_0$ is independent of $\e$.
\begin{table}[htbp]
\caption{The values of $\mathrm M_{m}^\e:={\big \| \, \mathbb
P^\e_m \psi_{\rm in} \big\|}_{L^2(\R)}$, for
$\e=\frac{1}{32}$:}\label{tab_energy}\vspace{-1mm}
\begin{center}\footnotesize
\renewcommand{\arraystretch}{1.3}
\begin{tabular}{ccccc}
\hline $m$ & 1 & 2& 3& 4\\ \hline ${\mathrm M}^\e_m$
&     $7.91\E-1$
&     $1.11\E-1$
&     $5.92\E-1$
&     $8.80\E-2$
\\ \hline $m$ & 5& 6& 7& 8\\ \hline ${\mathrm M}^\e_m$
&     $8.67\E-2$
&     $2.81\E-3$
&     $2.80\E-3$
&     $4.98\E-5$
\\ \hline
\end{tabular}
\end{center}
\end{table}

To compute the evolution of these initial data we shall take into
account $M\geq M_0$ bands. Note that only in cases where
$U(x)\equiv 0$ one can take $M$ to be identical to $M_0$, the
initial band cut-off. The reason is that if $U(x)$ is nonzero Step
2 in the BD algorithm given above mixes all bands. In particular
all the $\psi_m(t)$ are no longer orthogonal to each other.
Roughly speaking however, if $\e$ is very small, all band spaces
$\mathcal H_m$ remain ``almost orthogonal'' and thus the mass
within each Bloch band, \ie $\mathrm M_m^\e(t):={\big \| \,
\mathbb P^\e_m \psi(t) \big\|}^2_{L^2(\R)}$ is ``almost
conserved''. More precisely it is conserved up to errors $\O(\e)$
on time scales $\O(1)$. Thus, by checking mass conservation after
each time step one gets a rather reliable measure on the amount of
mixing of the bands. In other words if the mass conservation after
some time steps gets worse, one has to take into account more
bands to proceed.

We find numerically that the use of $M = M_0 \approx 8$ bands
already yields satisfactory results for $\e = \frac{1}{32}$. In
the following though we shall even compute  $M = 32$ energy bands,
which is by far sufficient for our purposes (even if $\e =
\frac{1}{2}$). Note that the number of required bands $M$ depends
on the regularity properties of $U(x)$, as well as on the
considered time-scales (which might be even longer than $\O(1)$,
the case considered here). This approximation problem is more or
less analogous to the one appearing in spectral schemes for PDEs
with non-smooth coefficients.

Concerning slowly varying, external potentials $U$, we shall
choose, on the one hand, smooth functions which are either of the
form \be \label{eqlin} U(x) = \mathcal E x, \ee modelling a
constant (electric) force field $\mathcal E \in \R$, or given by a
\emph{harmonic oscillator} type potential \be\label{eqhar_p}
U(x)=|x-\pi|^2. \ee On the other hand, we shall also consider the
case of an external (non-smooth) \emph{step potential}, \ie
\be\label{eqstep}
U(x)=\left\{\ba{ll}1, & x\in \left[\f{\pi}{2},\f{3\pi}{2}\right]\\
0, & \mbox{else} . \ea\right. \ee Within the setting described
above, we shall focus on two particular choices for the lattice
potential, namely:

\begin{example}[\textbf{Mathieu's model}]\label{nuex1}
The so-called \emph{Mathieu's model}, \ie \be \label{mathieu}
V_\Gm(x)=\cos(x), \ee as already considered in \cite{GoMa}. (For
applications in solid state physics this is rather unrealistic,
however it fits quite good with experiments on Bose-Einstein
condensates in optical lattices.) In this case all Fourier
coefficients $\widehat V(\lambda)$, appearing in \eqref{Vgm} are
zero, except for $\widehat V(\pm 1 ) = \frac{1}{2}$ and thus
$\H(k)$, given in \eqref{Hk}, simplifies to a tri-diagonal matrix.
\end{example}

\begin{example}[\textbf{Kronig-Penney's model}]\label{nuex2}
The so-called \emph{Kronig-Penney's model}, \ie \be \label{KrP}
V_\Gm(x)=1-\sum_{\gamma \in\Z}{\bf 1}_{x\in\left[\f{\pi}{2}+2\pi
\gamma,\f{3\pi}{2}+2 \pi \gamma\right]}, \ee where ${\bf
1}_{\Omega}$ denotes the characteristic function of a set $\Omega
\subset \R$. In contrast to Mathieu's model this case comprises a
non-smooth lattice potential. The corresponding Bloch eigenvalue
problem is known to be explicitly solvable (see, \eg,
\cite{GoMa}).
\end{example}

In order to compare the different numerical algorithms we denote
by $\psi^{\rm ts}(t,x)$ the solution gained from the
time-splitting spectral method, whereas $\psi^{\rm bd}(t,x)$
denotes the solution obtained via the new method base on Bloch's
decomposition. Both methods will be compared to the ``exact''
solution $\psi^{\rm ex}(t,x)$, which is obtained using a very fine
spatial grid. We consider the following errors \be\label{error1}
\begin{aligned}
\Delta^{\rm bd/ts}_\infty(t):= & \,
{\big \| \, \psi^{\rm ex}(t,\cdot) - \psi^{\rm bd/ts}(t,\cdot)
\big \|}_{L^\infty(\R)},\\
\Delta^{\rm bd/ts}_2(t):= &\, {\big \| \, \psi^{\rm ex}(t,\cdot) -
\psi^{\rm bd/ts}(t,\cdot) \big \|}_{L^2(\R)}
\end{aligned}
\ee between the ``exact solution'' and the corresponding solutions
obtained via the Bloch decomposition based algorithm resp. the
classical time splitting spectral method. The numerical
experiments are now done in a series of three different settings:
\begin{itemize}
\item First we shall study both cases of $V_\Gamma$, imposing
additionally $U(x)\equiv0$, \ie no external potential. The
obtained results are given in Table \ref{tb1_0}, where $\e=
\f{1}{2}$, $\f{1}{32}$, and $\f{1}{1024}$, respectively. In the
last case the oscillations are extremely spurious. As discussed
before, we can use only one step in time to obtain the numerical
solution, because the Bloch-decomposition method indeed is
``exact'' in this case (independently of $\e$). Thus, even if we
would refine the time steps in the BD algorithm we would not get
more accurate approximations. On the other hand, by using the
usual time-splitting method, one has to refine the time steps
(depending on $\e$) as well as the mesh size in order to achieve
the same accuracy. More precisely we find that $\Delta t=\O(\e)$,
$\Delta x=\O(\e^\alpha)$, for some $\alpha \geq 1$, is needed when
using TS (see also the computations given in \cite{GoMa}). In
particular $\alpha > 1$ is required for the case of a non-smooth
lattice potential $V_\Gamma$. (Note that if $V_\Gamma = 0$ it is
well known that $\Delta x=\O(\e)$, is sufficient, \cf \cite{BJM,
BJM1, HJMSZ}).

\item In a second series of numerical experiments we shall
consider only Example \ref{nuex1} for the periodic potential but
taking into account all three cases of the external potentials
$U$, as given above. In Fig. \ref{fig18}--\ref{fig17}, we show the
obtained numerical results for $\e=\f{1}{2}$,  and $\e =
\f{1}{1024}$, respectively. We observe that, if $\e=\O(1)$, the
Bloch-decomposition method gives almost the same results as
time-splitting spectral method. However, if $\e\ll 1$, we can
achieve quite good accuracy by using the Bloch-decomposition
method with $\Delta t=\O(1)$ and $\Delta x=\O(\e)$. On the other
hand, using the standard TS algorithm, we again have to rely on
much finer spatial grids and time steps to achieve the same
accuracy. \item We finally show the numerical results obtained by
combining external fields and a non-smooth lattice potential given
by Example \ref{nuex2}. As before we include all three cases for
the external potential $U$. The cases $\e=\f{1}{2}$, and
$\f{1}{1024}$ are studied and the obtained results are given in
Fig. \ref{fig27}--\ref{fig26}, respectively. We observe that the
results of the Bloch-decomposition are much better than the
time-splitting spectral method, even if $\e = \f{1}{2}$. Moreover,
as $\e$ gets smaller, the advantages of the Bloch-decomposition
method are even better visible.
\end{itemize}
To convince ourselves that only a few Bloch bands contribute to
$\|\psi \|_{L^2(\R)}$, even after time steps $\O(1)$, we show in
the following table the numerical values of $\mathrm
M^\e_m(t)={\big \| \, \mathbb P^\e_m \psi(t) \big\|}^2_{L^2(\R)}$,
for $m=1,\cdots,8$, corresponding to the solution of Example
\ref{nuex1} with $U$ given by \eqref{eqhar_p}.
\begin{table}[htbp]
\caption{The mass of $\psi_{\rm }(t,x)$, solution to Example
\ref{nuex1} with external potential \eqref{eqhar_p}, decomposed
into the Bloch bands for $\e=\frac{1}{32}$ at time
$t=1$:}\label{tab_energy1}
\begin{center}\footnotesize
\renewcommand{\arraystretch}{1.3}
\begin{tabular}{ccccc}
\hline $m$ & 1 & 2& 3& 4\\ \hline $\mathrm M_m^\e$
&     $7.89\E-1$
&     $1.10\E-2$
&     $5.92\E-1$
&     $9.38\E-2$
\\ \hline $m$ & 5& 6& 7& 8\\ \hline $\mathrm M_m^\e$
&     $7.15\E-2$
&     $3.50\E-3$
&     $1.80\E-3$
&     $5.63\E-5$
\\ \hline
\end{tabular}
\end{center}
\end{table}

We also check the conservation of the total (discrete) mass, \ie
${\|\psi(t) \|}_{l^2(\mathcal D)}$. We find that numerically it is
of the order $10^{-6}$ for the smooth lattice potential
\eqref{mathieu} and $10^{-3}$ for the non-smooth case \eqref{KrP}.
The latter however can be improved by using a refined spatial grid
and more time steps.

In summary we find (at least for our one dimensional computations)
that, relying on the new Bloch-decomposition based algorithm, one
can use much larger time steps, and sometimes even a coarser
spatial grid, to achieve the same accuracy as for the usual
time-splitting spectral method. This is particularly visible in
cases, where the lattice potential is non longer smooth and $\e\ll
1$. Indeed in these cases the BD algorithm turns out to be
\emph{considerably faster} than the TS method.
\begin{remark}
In view of our results the earlier numerical studies based on TS
methods \cite{GoMa, Go, GoMau}, should be taken with some care, in
particular when comparing the full Schr\"odinger solution to the
semiclassical approximation beyond caustics.
\end{remark}
\begin{table}[htbp]
\caption{The results of Example \ref{nuex1} with
$U(x)=0$:}\label{tb1_0}
\begin{center}\footnotesize
\renewcommand{\arraystretch}{1.3}
{Spatial discretization error test at time $t=1.0$ for $\e=1/2$.\\
For TS $\tg t=0.0001$ and for BD $\tg t=1$. }
\begin{tabular}{c|cccc}\hline
mesh size $\tg x/\e$& $1/2$ & $1/4$ & $1/8$ &   $1/16$ \\ \hline
$\ba{c}\vspace*{-4mm} \\
\left\|\psi^{\rm ts}_{\tg x,\tg t}(t,\cdot)-\psi^{\rm
ex}(t,\cdot)\right\|_{l^2}\\ \vspace*{-4mm} \ea$
 &  4.33E-1 & 2.53E-1 &  2.80E-2 & 6.42E-6 \\ \hline
convergence order & & 0.8 & 3.2& 12.1 \\ \hline \hline
$\ba{c}\vspace*{-4mm} \\
\left\|\psi^{\rm bd}_{\tg x,\tg t}(t,\cdot)-\psi^{\rm
ex}(t,\cdot)\right\|_{l^2}\\ \vspace*{-4mm} \ea$
 &  3.01E-1& 1.95E-1 &  1.39E-2 & 1.17E-6 \\ \hline
convergence order & & 0.6& 3.8& 13.5  \\ \hline 
\end{tabular}\vspace{5mm}

{Spatial discretization error test at time $t=0.1$ for $\e=1/32$.\\
For TS $\tg t=0.00001$ and for BD $\tg t=0.1$. }
\begin{tabular}{c|cccc}\hline
mesh size $\tg x/\e$& $1/2$ & $1/4$ & $1/8$ &   $1/16$ \\ \hline
$\ba{c}\vspace*{-4mm} \\
\left\|\psi^{\rm ts}_{\tg x,\tg t}(t,\cdot)-\psi^{\rm
ex}(t,\cdot)\right\|_{l^2}\\ \vspace*{-4mm} \ea$
 &  2.88E-1 & 1.08E-1 &  9.63E-4 & 1.33E-7 \\ \hline
convergence order & & 1.4 & 6.8& 12.8 \\ \hline \hline
$\ba{c}\vspace*{-4mm} \\
\left\|\psi^{\rm bd}_{\tg x,\tg t}(t,\cdot)-\psi^{\rm
ex}(t,\cdot)\right\|_{l^2}\\ \vspace*{-4mm} \ea$
 &  2.53E-1& 7.34E-2 &  8.97E-4 &  4.95E-10 \\ \hline
convergence order & & 1.8& 6.4& 20.8  \\ \hline 
\end{tabular}\vspace{5mm}

{Spatial discretization error test at time $t=0.01$ for $\e=1/1024$.\\
For TS $\tg t=0.000001$ and for BD  $\tg t=0.01$. }
\begin{tabular}{c|cccc}\hline
mesh size $\tg x/\e$& $1/2$ & $1/4$ & $1/8$ &   $1/16$ \\ \hline
$\ba{c}\vspace*{-4mm} \\
\left\|\psi^{\rm ts}_{\tg x,\tg t}(t,\cdot)-\psi^{\rm
ex}(t,\cdot)\right\|_{l^2}\\ \vspace*{-4mm} \ea$
 &  5.14E-1 & 1.94E-1 &  1.08E-3 & 6.08E-8 \\ \hline
convergence order & & 1.4 & 7.5& 14.1 \\ \hline \hline
$\ba{c}\vspace*{-4mm} \\
\left\|\psi^{\rm bd}_{\tg x,\tg t}(t,\cdot)-\psi^{\rm
ex}(t,\cdot)\right\|_{l^2}\\ \vspace*{-4mm} \ea$
 &  2.64E-1& 6.83E-2 &  2.29E-4 &  1.71E-10 \\ \hline
convergence order & & 2.0& 8.2& 20.4  \\ \hline 
\end{tabular}
\end{center}
\end{table}

\begin{table}[htbp]
\caption{The results of Example \ref{nuex1} with linear external
potential \eqref{eqlin}:}\label{tb1_1}
\begin{center}\footnotesize
\renewcommand{\arraystretch}{1.3}
{Spatial discretization error test at time $t=0.1$ for $\e=1/2$.\\
For TS $\tg t=0.0001$ and for BD $\tg t=0.01$. }
\begin{tabular}{c|cccc}\hline
mesh size $\tg x/\e$& $1/2$ & $1/4$ & $1/8$ &   $1/16$ \\ \hline
$\ba{c}\vspace*{-4mm} \\
\left\|\psi^{\rm ts}_{\tg x,\tg t}(t,\cdot)-\psi^{\rm
ex}(t,\cdot)\right\|_{l^2}\\ \vspace*{-4mm} \ea$
 &  2.73E-1 & 9.22E-2 &  5.78E-3 & 4.73E-6 \\ \hline
convergence order & & 1.6 & 4.0& 10.3 \\ \hline \hline
$\ba{c}\vspace*{-4mm} \\
\left\|\psi^{\rm bd}_{\tg x,\tg t}(t,\cdot)-\psi^{\rm
ex}(t,\cdot)\right\|_{l^2}\\ \vspace*{-4mm} \ea$
 &  3.15E-1& 1.55E-1 &  1.32E-2 & 3.36E-6 \\ \hline
convergence order & & 1.0 & 3.6& 11.9  \\ \hline 
\end{tabular}\vspace{5mm}

{Spatial discretization error test at time $t=0.01$ for $\e=1/1024$.\\
For TS $\tg t=0.00001$ and for BD  $\tg t=0.001$. }
\begin{tabular}{c|cccc}\hline
mesh size $\tg x/\e$& $1/2$ & $1/4$ & $1/8$ &   $1/16$ \\ \hline
$\ba{c}\vspace*{-4mm} \\
\left\|\psi^{\rm ts}_{\tg x,\tg t}(t,\cdot)-\psi^{\rm
ex}(t,\cdot)\right\|_{l^2}\\ \vspace*{-4mm} \ea$
 &  5.22E-1 & 1.98E-1 &  1.53E-2 & 3.19E-5 \\ \hline
convergence order & & 1.4 & 3.7& 8.9 \\ \hline \hline
$\ba{c}\vspace*{-4mm} \\
\left\|\psi^{\rm bd}_{\tg x,\tg t}(t,\cdot)-\psi^{\rm
ex}(t,\cdot)\right\|_{l^2}\\ \vspace*{-4mm} \ea$
 &  4.71E-1& 1.61E-1 &  9.17E-3 &  6.08E-6 \\ \hline
convergence order & & 1.5& 4.1& 10.6  \\ \hline 
\end{tabular}\vspace{5mm}

{Temporal discretization error test at $t=0.1$ for $\e=1/2$ and
$\tg x/\e=1/128$.}
\begin{tabular}{c|cccc}\hline
time step $\tg t$& $1/10$ & $1/20$ & $1/40$ &   $1/80$ \\ \hline
$\ba{c}\vspace*{-4mm} \\
\left\|\psi^{\rm ts}_{\tg x,\tg t}(t,\cdot)-\psi^{\rm
ex}(t,\cdot)\right\|_{l^2}\\ \vspace*{-4mm} \ea$
 &  2.59E-4 & 6.47E-5 &  1.62E-5 & 4.04E-6 \\ \hline
convergence order & & 2.0 & 2.0& 2.0 \\ \hline \hline
$\ba{c}\vspace*{-4mm} \\
\left\|\psi^{\rm bd}_{\tg x,\tg t}(t,\cdot)-\psi^{\rm
ex}(t,\cdot)\right\|_{l^2}\\ \vspace*{-4mm} \ea$
 &  4.86E-5& 1.23E-5 &  3.08E-6 & 7.60E-7 \\ \hline
convergence order & & 2.0 & 2.0& 2.0  \\ \hline 
\end{tabular}\vspace{5mm}

{Temporal discretization error test at $t=0.01$ for $\e=1/1024$
and $\tg x/\e=1/128$.}
\begin{tabular}{c|cccc}\hline
time step $\tg t$& $1/1000$ & $1/2000$ & $1/4000$ &   $1/8000$ \\
\hline
$\ba{c}\vspace*{-4mm} \\
\left\|\psi^{\rm ts}_{\tg x,\tg t}(t,\cdot)-\psi^{\rm
ex}(t,\cdot)\right\|_{l^2}\\ \vspace*{-4mm} \ea$
 &  6.60E-2 & 1.54E-2 &  3.81E-3 & 9.45E-4 \\ \hline
convergence order & & 2.1 & 2.0& 2.0 \\ \hline \hline time step
$\tg t$& $1/100$ & $1/200$ & $1/400$ &   $1/800$ \\ \hline
$\ba{c}\vspace*{-4mm} \\
\left\|\psi^{\rm bd}_{\tg x,\tg t}(t,\cdot)-\psi^{\rm
ex}(t,\cdot)\right\|_{l^2}\\ \vspace*{-4mm} \ea$
 &  3.32E-3& 7.54E-4 &  1.42E-4 & 3.16E-5 \\ \hline
convergence order & & 2.1 & 2.4& 2.2  \\ \hline 
\end{tabular}
\end{center}
\end{table}

\begin{table}[htbp]
\caption{The results of Example \ref{nuex2} with harmonic external
potential \eqref{eqhar_p}:}\label{tb2_1}
\begin{center}\footnotesize
\renewcommand{\arraystretch}{1.3}
{Spatial discretization error test at time $t=0.1$ for $\e=1/2$.\\
For TS $\tg t=0.0001$ and for BD $\tg t=0.01$. }
\begin{tabular}{c|cccc}\hline
mesh size $\tg x/\e$& $1/2$ & $1/4$ & $1/8$ &   $1/16$ \\ \hline
$\ba{c}\vspace*{-4mm} \\
\left\|\psi^{\rm ts}_{\tg x,\tg t}(t,\cdot)-\psi^{\rm
ex}(t,\cdot)\right\|_{l^2}\\ \vspace*{-4mm} \ea$
 &  2.71E-1 & 8.87E-2 &  5.19E-3 & 1.32E-4 \\ \hline
convergence order & & 1.6 & 4.1& 5.3 \\ \hline \hline
$\ba{c}\vspace*{-4mm} \\
\left\|\psi^{\rm bd}_{\tg x,\tg t}(t,\cdot)-\psi^{\rm
ex}(t,\cdot)\right\|_{l^2}\\ \vspace*{-4mm} \ea$
 &  3.23E-1& 9.08E-2 &  7.03E-3 & 1.27E-4 \\ \hline
convergence order & & 1.8 & 3.7& 5.8  \\ \hline 
\end{tabular}\vspace{5mm}

{Spatial discretization error test at time $t=0.01$ for $\e=1/1024$.\\
For TS $\tg t=0.00001$ and for BD  $\tg t=0.001$. }
\begin{tabular}{c|cccc}\hline
mesh size $\tg x/\e$& $1/2$ & $1/4$ & $1/8$ &   $1/16$ \\ \hline
$\ba{c}\vspace*{-4mm} \\
\left\|\psi^{\rm ts}_{\tg x,\tg t}(t,\cdot)-\psi^{\rm
ex}(t,\cdot)\right\|_{l^2}\\ \vspace*{-4mm} \ea$
 &  3.99E-1 & 3.67E-1 &  2.19E-1 & 1.10E-1 \\ \hline
convergence order & & 0.1 & 0.7& 1.0 \\ \hline \hline
$\ba{c}\vspace*{-4mm} \\
\left\|\psi^{\rm bd}_{\tg x,\tg t}(t,\cdot)-\psi^{\rm
ex}(t,\cdot)\right\|_{l^2}\\ \vspace*{-4mm} \ea$
 &  2.06E-1& 5.64E-2 &  8.16E-3 &  6.40E-4 \\ \hline
convergence order & & 1.9& 2.8& 3.7  \\ \hline 
\end{tabular}\vspace{5mm}

{Temporal discretization error test at $t=0.1$ for $\e=1/2$ and
$\tg x/\e=1/128$.}
\begin{tabular}{c|cccc}\hline
time step $\tg t$& $1/10$ & $1/20$ & $1/40$ &   $1/80$ \\ \hline
$\ba{c}\vspace*{-4mm} \\
\left\|\psi^{\rm ts}_{\tg x,\tg t}(t,\cdot)-\psi^{\rm
ex}(t,\cdot)\right\|_{l^2}\\ \vspace*{-4mm} \ea$
 &  1.02E-3 & 6.41E-4 &  3.80E-4 & 2.18E-4 \\ \hline
convergence order & & 0.7 & 0.8& 0.8 \\ \hline \hline
$\ba{c}\vspace*{-4mm} \\
\left\|\psi^{\rm bd}_{\tg x,\tg t}(t,\cdot)-\psi^{\rm
ex}(t,\cdot)\right\|_{l^2}\\ \vspace*{-4mm} \ea$
 &  4.20E-6 &1.02E-6 &  2.22E-7 & 5.56E-8 \\ \hline
convergence order & & 2.0 & 2.2& 2.0  \\ \hline 
\end{tabular}\vspace{5mm}

{Temporal discretization error test at $t=0.01$ for $\e=1/1024$
and $\tg x/\e=1/128$.}
\begin{tabular}{c|cccc}\hline
time step $\tg t$& $1/1000$ & $1/2000$ & $1/4000$ &   $1/8000$ \\
\hline
$\ba{c}\vspace*{-4mm} \\
\left\|\psi^{\rm ts}_{\tg x,\tg t}(t,\cdot)-\psi^{\rm
ex}(t,\cdot)\right\|_{l^2}\\ \vspace*{-4mm} \ea$
 &  1.21E-1 & 1.18E-1 &  1.10E-1 & 1.10E-1 \\ \hline
convergence order & & 0.04 & 0.1& 0.0 \\ \hline \hline time step
$\tg t$& $1/100$ & $1/200$ & $1/400$ &   $1/800$ \\ \hline
$\ba{c}\vspace*{-4mm} \\
\left\|\psi^{\rm bd}_{\tg x,\tg t}(t,\cdot)-\psi^{\rm
ex}(t,\cdot)\right\|_{l^2}\\ \vspace*{-4mm} \ea$
 &  3.30E-5& 5.21E-6 &  1.23E-6 & 3.16E-7 \\ \hline
convergence order & & 2.6 & 2.1& 2.0  \\ \hline
\end{tabular}
\end{center}
\end{table}
\section{Asymptotic analysis in the semiclassical regime}\label{sec:asym}

For completeness we shall also compare the numerical solution of
the Schr\"odinger equation \eqref{s0} with its semiclassical
asymptotic description. To this end we shall rely on a multiple
scales WKB-type expansion methods, even though there are currently
more advanced tools at hand, \cf \cite{GMMP, PST, Te}. The WKB
method however has the advantage of given a rather simple and
transparent description of $\psi(t)$, solution to \eqref{s0}, for
$\e \ll 1$, (at least locally in-time). Since the Bloch
decomposition method itself is rather abstract we include this
approximative description here too, so that the reader gets a
better feeling for the appearing quantities. Moreover this
two-scale WKB method can also be used for \emph{nonlinear}
Schr\"odinger dynamics \cite{CMS}, a problem we shall study
numerically in an upcoming work.

\subsection{The WKB formalism}
To this end let us suppose that the initial condition is of
(two-scale) \emph{WKB-type}. More precisely assume
\bea\label{eq:wkb-IC} \psi_{\rm in}(x)=\sum^M_{m =
1}u_{m}\left(x,\frac{x}{\e}\right) \, \ex^{i \phi(x)/\e}, \eea
with some given \emph{real-valued phase} $\phi\in C^\infty(\R)$
and some given initial (complex-valued) \emph{band-amplitudes}
$u_m(x,y+2\pi)=u_m(x,y)$, each of which admits an asymptotic
description of the following form \be\label{WKBam} u_m(x,y)\sim
{u}^0_{m} \left(x,y\right) +\e {u}^1_{m}\left(x, y\right)+\O(\e^2)
\quad \forall \, m \in \N. \ee Here and in the following we shall
only be concerned with the leading order asymptotic description.
\begin{remark}
Note that we do consider only a \emph{single} initial WKB-phase
$\phi(x)$ for all bands $m\in \N$. We could of course also allow
for more general cases, like one WKB-phase for each band or even a
superposition of WKB-states within each band. However in order to
keep the presentation clean we hesitate to do so. The standard WKB
approximation, for non-periodic problems, involves real-valued
amplitudes $\tilde u^0(x), \tilde u^1(x), \dots$ which only depend
on the slow scale.
\end{remark}
It is well known then, \cf \cite{CMS, GRT}, that the leading order
term ${u}^0_m$, $m\in \N$, can be decomposed as \be \label{u0}
{u}^0_{m} \left(x,\frac{x}{\e}\right) = f_{m}(x)\chi_m
\left(\frac{ x}{\e}, \partial_x \phi(x) \right), \ee where we
assume the m-th energy band to be \emph{non-degenerated} (for
simplicity) and \emph{isolated} from the rest of spectrum. We can
choose an arbitrary $f_m \in \mathcal S(\R)$. In other words,
there is an adiabatic decoupling between the slow scale $x$ and
fast scale $x/\e$. Indeed, a lengthy calculation, invoking the
classical stationary phase argument, \cf chapter 4.7 in
\cite{BLP}, shows that in this case the band projection $\mathbb
P^\e_m \psi$ can be approximated via \be \mathbb P_m^\e \psi(x)
\sim  f_m(x) \, \chi_m\left(\frac{x}{\e},\p_x\phi(x)\right) \ex^{i
\phi(x)/\e}+\O(\e). \ee This approximate formula shows the origin
of the high oscillations induced either by $V_\Gamma$, described
by $\chi_m$, or by the dispersion, described by $\phi(x)$. We note
that in general the higher order terms (in $\e$), such as $u^1_m$
etc., are of a more complicated structure than \eqref{u0}, but we
shall neglect these terms in what follows (see, \eg, \cite{CMS}
for more details). One consequently finds that $\psi(t)$ obeys a
leading order asymptotic description of the form \bea
\label{WKBex} \psi(t,x)\sim \sum_{m=1}^M a_m(t,x) \,
\chi_m\left(\frac{x}{\e},\p_x\phi_m(t,x)\right) \ex^{i
\phi_m(t,x)/\e}+\O(\e), \eea where $\phi_m(t,x)\in
C^\infty([0,t_c)\times \R)$ satisfies the $m$th band
\emph{Hamilton-Jacobi equation} \be\label{eq:H-J_1} \left \{
\begin{aligned}
& \p_t\phi_m(t,x)+E_m(\p_x\phi_m)+U(x)=  \, 0,\quad m\in \N,\\
& \phi_m\big|_{t=0} = \, \phi(x).
\end{aligned}
\right. \ee Also, the (complex-valued) leading order WKB-amplitude
$a_m(t,x)\in C^\infty([0,t_c)\times \R)$ satisfies the following
\emph{semiclassical transport equations} \be\label{eq:trans_1}
\left \{
\begin{aligned}
& \p_t a_m + \p_k E_m(\partial_x\phi_m)\p_x a_m+ \frac{1}{2} \,
\partial_x(\partial_k E_m(\partial_x\phi_m))a_m -
(\beta_m(t,x)\partial_x U(x)) \, a_m=0,\\
& a_m\big|_{t=0}=f_{m}(x).\end{aligned} \right. \ee with
$\beta_m(t,x):={\langle\chi_m(y,k),\p_{k}\chi_m(y,k)\rangle}_{L^2(\mathcal
C)}$, evaluated at $k=\p_x\phi_m$, the so-called \emph{Berry phase
term}.

\begin{remark}
Note that the Berry term is purely imaginary, \ie $\beta_m(t,x)
\in \I \R$, which implies the following conservation law \be \p_t
|a_m|^2  + \partial_x\left(\partial_k E_m(\partial_x\phi_m)
|a_m|^2 \right) =0 \quad \forall \, m \in \N. \ee
\end{remark}
Of course the above given WKB-type expansion method is only valid
up to the (in general finite) time $0\leq t_{\rm c}<\infty$, the
\emph{caustic onset-time} in the solution of \eqref{eq:H-J_1}.
Here we shall simply assume that $t_c >0$ holds, \ie no caustic is
formed at time $t=0$, which is very well possible in general. We
note that in the considered numerical examples below we indeed
have $t_c>0$ and we refer to \cite{Ca} for a broader discussion on
this. For $t\geq t_{\rm c}$ one would need to superimpose several
WKB-type solutions corresponding to the multi-valued solutions of
the flow map $(x,t)\mapsto X_t(x)\equiv X_t(x;\partial_x\phi(x))$,
where
\begin{equation}
\label{semi} \left \{
\begin{aligned}
\dot X_t = & \ \partial_k E_m(\Xi_t),\quad X_{0}=x,\\
\dot \Xi_t = & \ -\partial_x U(X_t) ,\quad \Xi_0=\partial_x
\phi(x).
\end{aligned}
\right.
\end{equation}
Numerically we shall use the relaxation method introduced in
\cite{JiXi} to solve the Hamilton-Jacobi equation
\eqref{eq:H-J_1}. Consequently we can solve the system of
transport equations \eqref{eq:trans_1} by a time-splitting
spectral scheme similar to the ones used above.

\subsection{Numerical examples} We shall finally study the WKB approach,
briefly described above,
by some numerical examples. Denote by \be \psi^{\rm sc} (t,x):=
\sum_{m=1}^M f_m(t,x) \chi_m \left(\f{x}{\e},\partial_x
\phi_m\right)e^{i\phi_m(t,x)/\e}, \ee the approximate
semiclassical solution to the Schr\"odinger equation \eqref{s0}.
In the following examples we only take into account a harmonic
external potential of the form \eqref{eqhar_p}.

\begin{example}[\textbf{Mathieu's model}]\label{nuscex1}
We first consider Mathieu's model \eqref{mathieu} and choose
initial condition in the form \be\label{eq:ini_sc1}
\begin{aligned}
\psi_{\rm in}(x)=
\ex^{-5(x-\pi)^2}\chi_{1}\left(\frac{x}{\e},0\right),
\end{aligned}
\ee \ie we choose $\phi(x)=0$ and restrict ourselves to the case
of only one band with index $m=1$. (Since $E_1(k)$ is an isolated
band the analytical results of \cite{BLP, CMS, GRT}, then imply
that we can neglect the contributions from all other bands $m >1$
up to errors of order $\O(\e)$ in $L^2(\R)\cap L^\infty(\R)$,
uniformly on compact time-intervals.) In this case, we numerically
find that no caustic is formed within the solution of
\eqref{eq:H-J_1} at least up to $t=1$, the largest time in our
computation. Note that \eqref{eq:ini_sc1} concentrates at the
minimum of the first Bloch band, where it is known that \be
E_m(k)\approx \frac{|k|^2}{2m^*} + E_m(0), \ee This is the
so-called \emph{parabolic band approximation}, yielding an
effective mass $m^* \in \R$. In Table
\ref{tb01}, 
we show the results with an additional harmonic external
potential, \cf \eqref{eqhar_p}, for $\e=\f{1}{32}$ and $\e=
\f{1}{1024}$ respectively.
\begin{table}[htbp]
\caption{Difference between the asymptotic solution and the
Schr\"odinger equation for example \ref{nuscex1} ($\tg t=10^{-4}$,
$\tg x=1/32768$):}\label{tb01}
\begin{center}\footnotesize
\renewcommand{\arraystretch}{1.3}
\begin{tabular}{ccc}\hline
$\e$ &         $\f{1}{32}$ & $\f{1}{1024}$ \\ \hline 
$\ba{c}\vspace{-4mm}\\ \dpm\sup_{0\le t\le 1}\left\|\, \psi (t,x)
- \psi^{\rm sc} (t,x) \, \right\|_{L^{2}(\R)}\ea$ & $6.68\E-3$ &
$3.08\E-4$ \vspace*{0.5mm}\\ \hline $\ba{c}\vspace{-4mm}\\
\dpm\sup_{0\le t\le 1}\left\|\, \psi (t,x) - \psi^{\rm sc} (t,x)
\, \right\|_{L^{\ift}(\R)}\ea$ & $5.57\E-2$ & $2.38\E-3$
\vspace*{0.5mm}\\ \hline
\end{tabular}
\end{center}
\end{table}\\
Note that these numerical experiments, together with those given
below, confirm the analytical results given in \cite{BLP, CMS,
GRT}.
\end{example}

\begin{example}[\textbf{Kronig-Penney's model}]\label{nuscex2}
Here, we consider again the Kronig-Penney's model \eqref{KrP}.
First we use the same initial condition as given in
\eqref{eq:ini_sc1} but with $m=2$, which again corresponds to an
isolated energy band. The corresponding numerical results for
$\e=\f{1}{32}$ and $\f{1}{1024}$ are shown in Table \ref{tb02}.
\begin{table}[htbp]
\caption{Difference between the asymptotic solution and the
Schr\"odinger equation for example \ref{nuscex2} for initial
condition \eqref{eq:ini_sc1} ($\tg t=10^{-4}$, $\tg
x=1/32768$):}\label{tb02}
\begin{center}\footnotesize
\renewcommand{\arraystretch}{1.3}
\begin{tabular}{ccc}\hline
$\e$ &         $\f{1}{32}$ & $\f{1}{1024}$ \\ \hline 
$\ba{c}\vspace{-4mm}\\ \dpm\sup_{0\le t\le 0.1}\left\|\, \psi
(t,x) - \psi^{\rm sc} (t,x) \, \right\|_{L^{2}(\R)}\ea$ &
$1.18\E-2$ & $1.08\E-3$ \vspace*{0.5mm}\\ \hline
$\ba{c}\vspace{-4mm}\\ \dpm\sup_{0\le t\le 0.1}\left\|\, \psi
(t,x) - \psi^{\rm sc} (t,x) \, \right\|_{L^{\ift}(\R)}\ea$ &
$9.34\E-2$ & $7.74\E-3$ \vspace*{0.5mm}\\ \hline
\end{tabular}
\end{center}
\end{table}\\
In a second case, we alternatively choose initial data of the form
\be\label{eq:ini_sc2} \psi_{\rm
in}(x)=\ex^{-5(x-\pi)^2}\chi_{2}\left(\frac{x}{\e} ,\sin(x)\right)
\ex^{-\I \cos(x)/\e}, \ee \ie $\phi(x)= - \cos(x)$. Here we find
(numerically) that the caustic onset time is roughly given by
$t_c\approx 0.24$, \cf Fig. \ref{phase_caustic}.
\begin{figure} 
\begin{center}\footnotesize
\resizebox{2in}{!} {\includegraphics{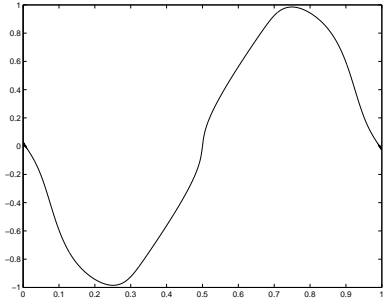}}
\end{center}
\caption{The graph of $\partial_x\phi_2(t,x)$ at
$t=0.24$.}\label{phase_caustic}
\end{figure}
The corresponding numerical results are given in Fig. \ref{fig41}
and Table \ref{tb03}.
\begin{table}[htbp]
\caption{Difference between the asymptotic solution and the
Schr\"odinger equation for example \ref{nuscex2}  for initial
condition \eqref{eq:ini_sc2} ($\tg t=10^{-4}$, $\tg
x=1/32768$):}\label{tb03}
\begin{center}\footnotesize
\renewcommand{\arraystretch}{1.3}
\begin{tabular}{ccc}\hline
$\e$ &         $\f{1}{32}$ & $\f{1}{1024}$ \\ \hline 
$\ba{c}\vspace{-4mm}\\ \dpm\sup_{0\le t\le 0.1}\left\|\, \psi
(t,x) - \psi^{\rm sc} (t,x) \, \right\|_{L^{2}(\R)}\ea$ &
$1.68\E-2$ & $3.19\E-3$ \vspace*{0.5mm}\\ \hline
$\ba{c}\vspace{-4mm}\\ \dpm\sup_{0\le t\le 0.1}\left\|\, \psi
(t,x) - \psi^{\rm sc} (t,x) \, \right\|_{L^{\ift}(\R)}\ea$ &
$2.73\E-1$ & $7.33\E-2$ \vspace*{0.5mm}\\ \hline
\end{tabular}
\end{center}
\end{table}
\end{example}

\section{Acknowledgement}\label{ACK}
The authors are grateful to Prof. Christian Ringhofer for fruitful
discussions on this work.



\begin{figure} \footnotesize
\begin{center}
\resizebox{1.5in}{!}{\includegraphics{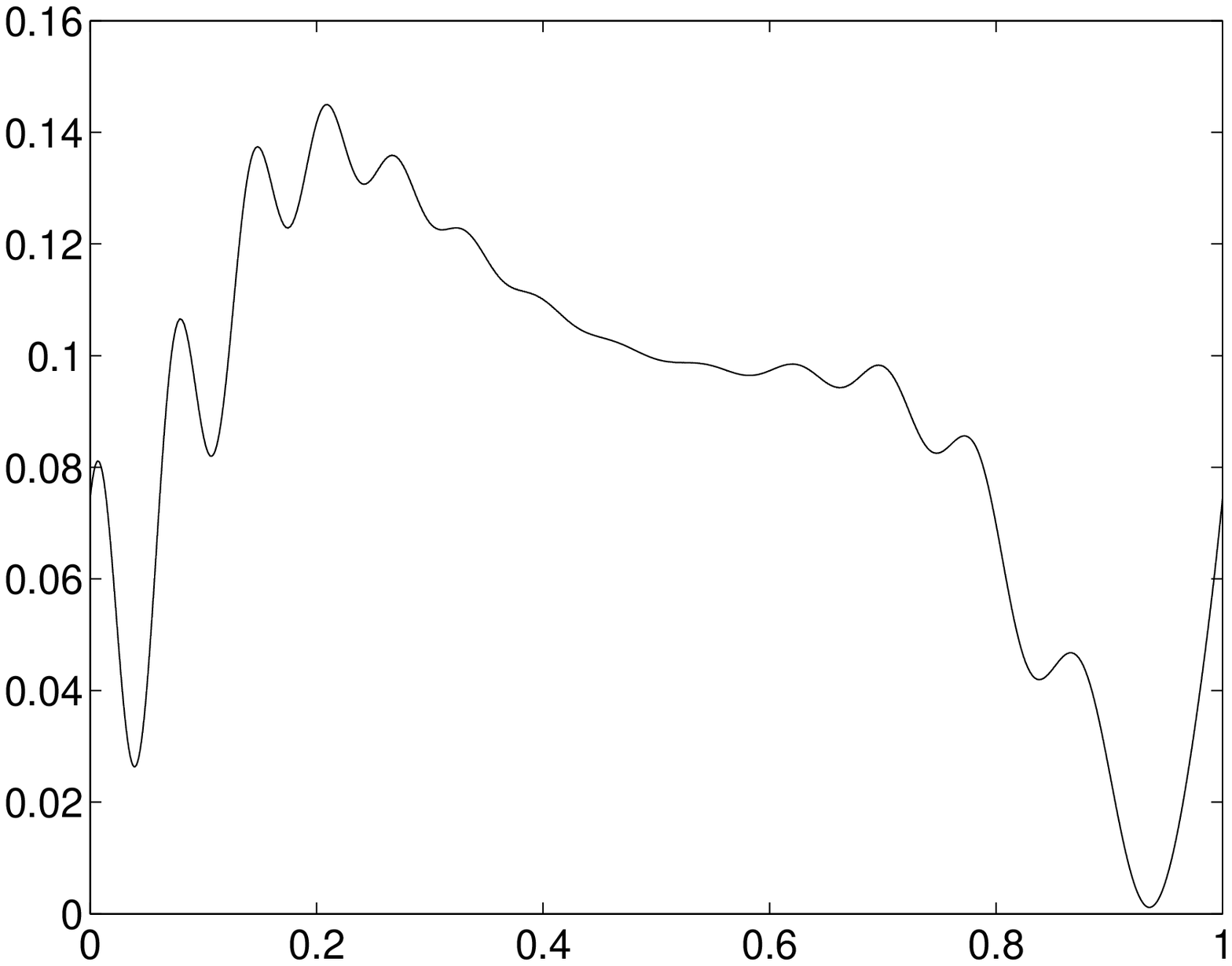}}
\resizebox{1.5in}{!}{\includegraphics{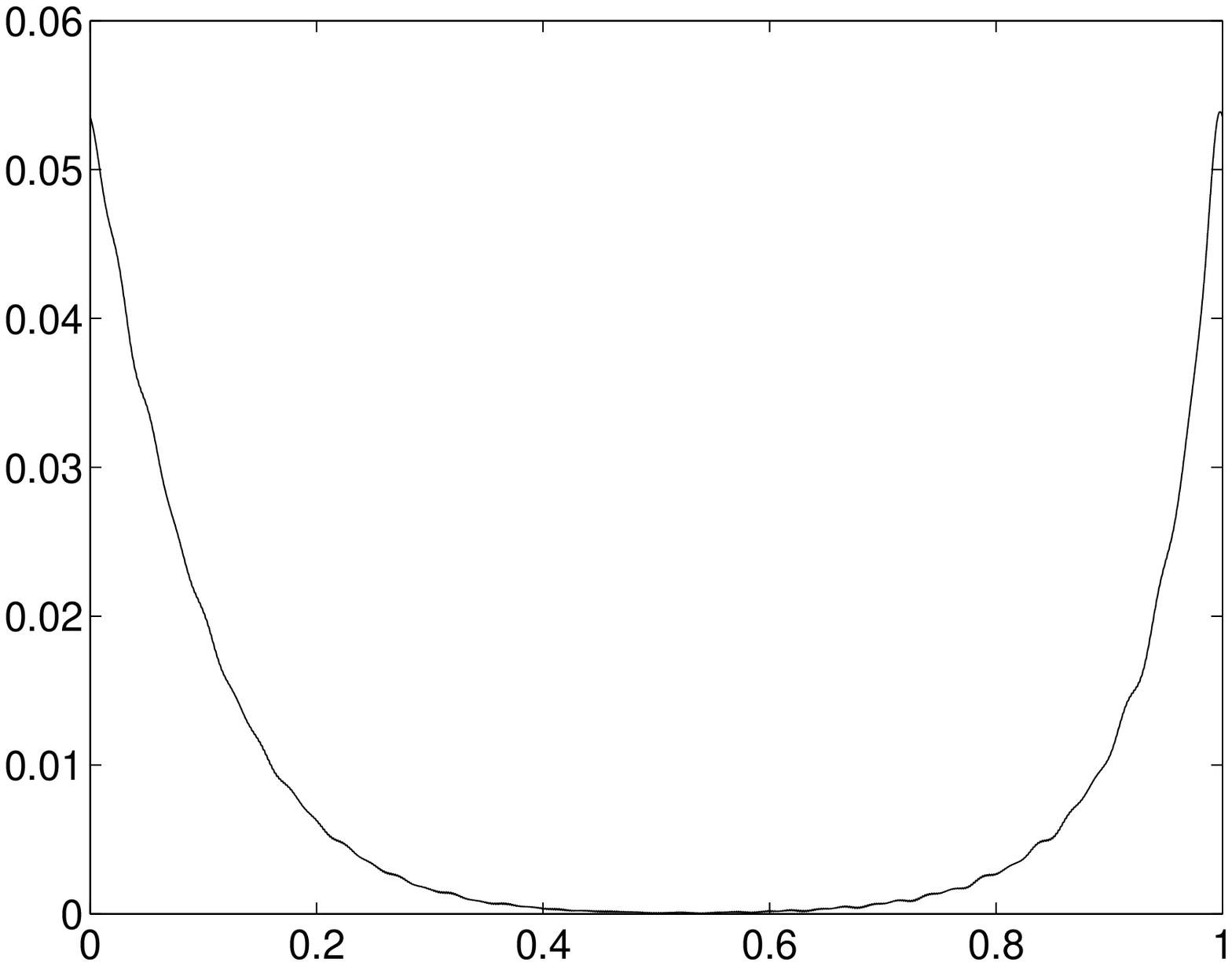}}
\resizebox{1.5in}{!}{\includegraphics{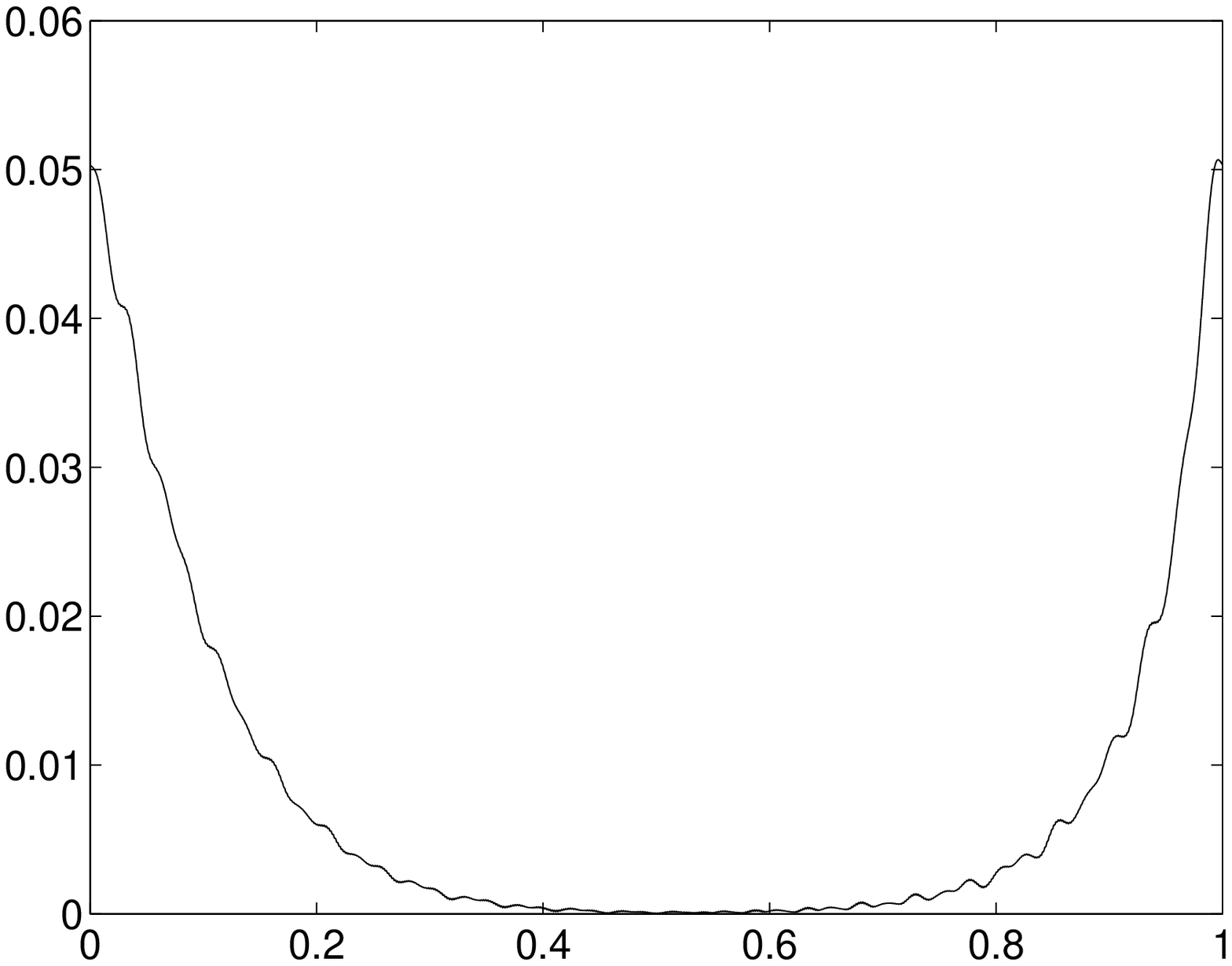}}

$|\psi^{\rm ex}(t,x)|^2$, $|\psi^{\rm ex}(t,x)-\psi^{\rm
ts}(t,x)|$ and $|\psi^{\rm ex}(t,x)-\psi^{\rm bd}(t,x)|$ at
$t=1.0$
\vspace{-1mm} \caption{Numerical results for example
\ref{nuex1} with $U(x)$ given by \eqref{eqlin} and $\e=\f{1}{2}$.
We use $\tg t=\f{1}{100}$, $\tg x=\f{1}{64}$ for the TS and the BD
method, and $\tg t=\f{1}{100000}$, $\tg x=\f{1}{8192}$ for the
``exact'' solution.}\vspace{-3mm}
\[\Delta^{\rm ts}_\infty(t)=5.39\E-2,\
\Delta^{\rm bd}_\infty(t)=5.07\E-2,\
\Delta^{\rm ts}_2(t)=1.56\E-2,\
\Delta^{\rm bd}_2(t)= 1.51\E-2. \]
\label{fig18}

\resizebox{1.5in}{!}{\includegraphics{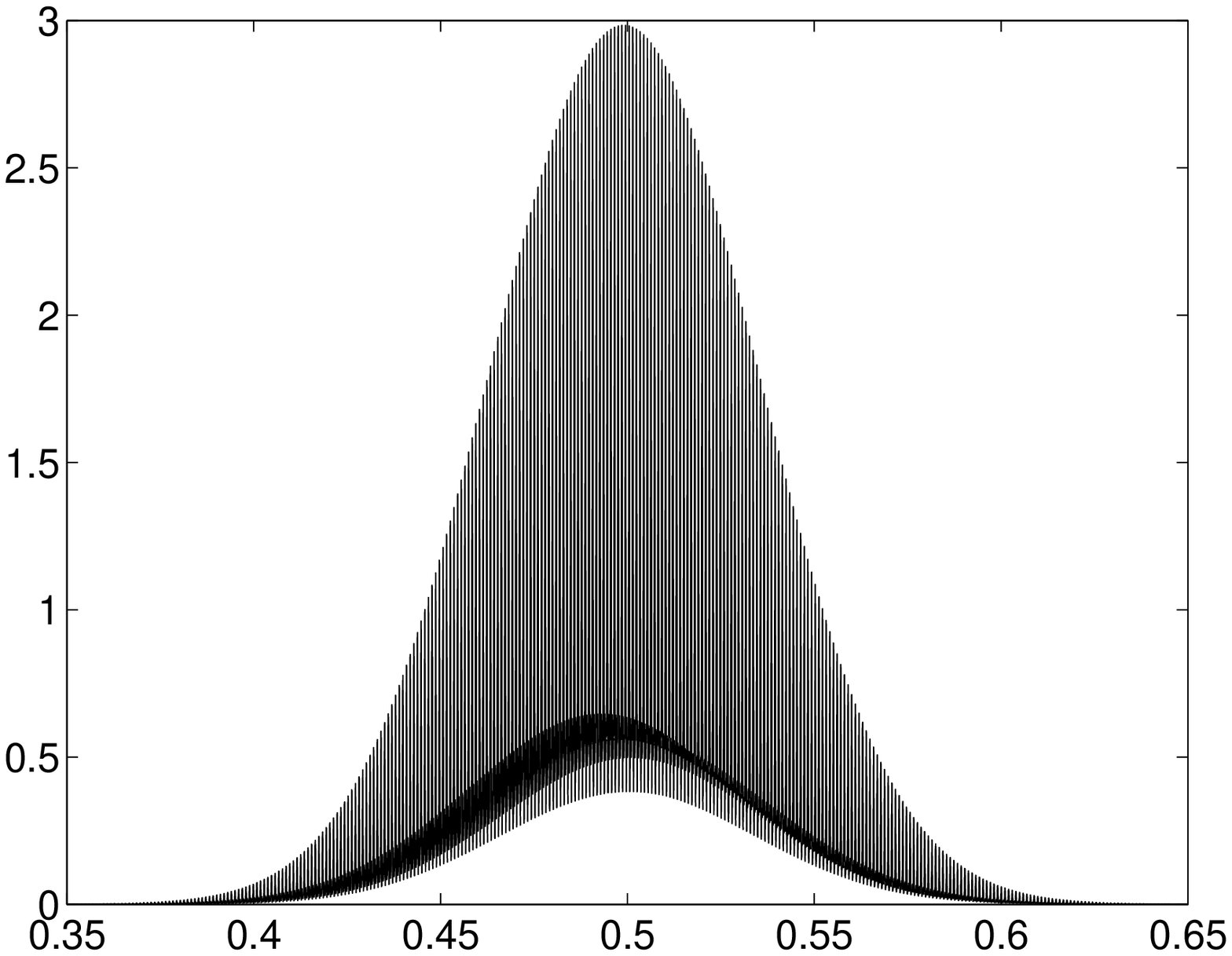}}
\resizebox{1.5in}{!}{\includegraphics{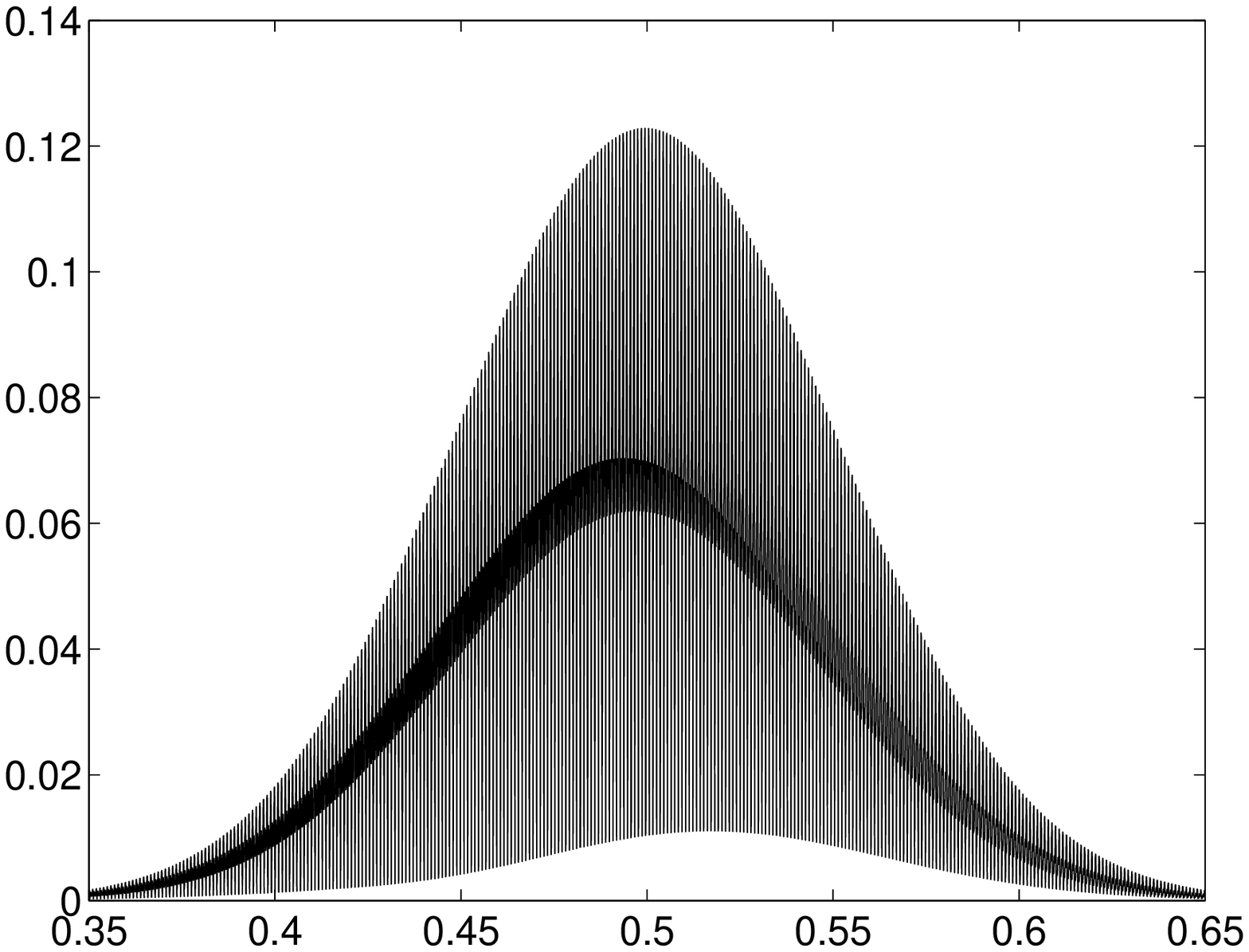}}
\resizebox{1.5in}{!}{\includegraphics{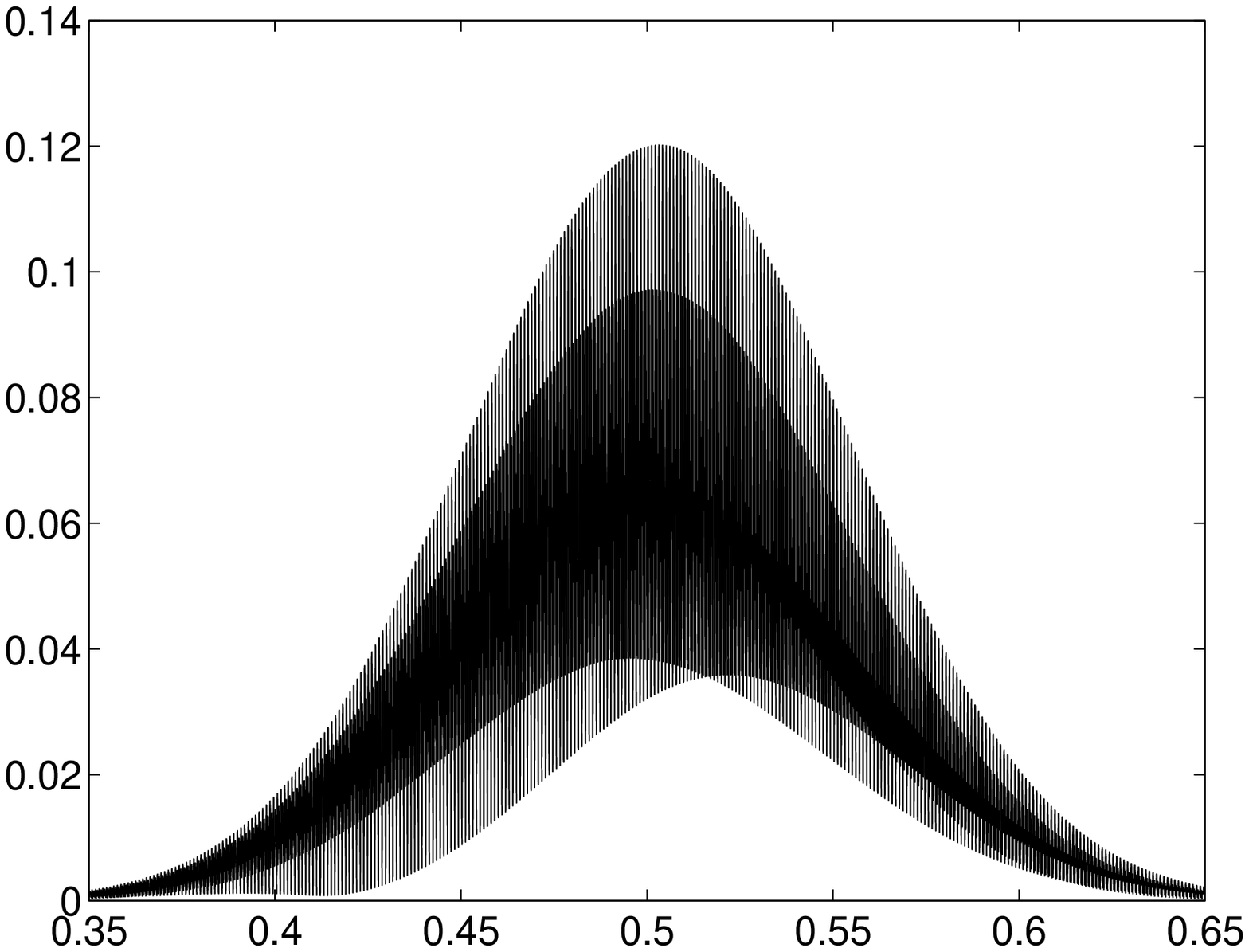}}

$|\psi^{\rm ex}(t,x)|^2$, $|\psi^{\rm ex}(t,x)-\psi^{\rm ts}(t,x)|$
and $|\psi^{\rm ex}(t,x)-\psi^{\rm bd}(t,x)|$ at $t=0.1$\vspace{-1mm}
\end{center}
\caption{Numerical results for example \ref{nuex1} with $U(x)$
given by \eqref{eqlin} and $\e=\f{1}{1024}$. We use $\tg
t=\f{1}{5000}$, $\tg x=\f{1}{16384}$ for the TS, $\tg
t=\f{1}{20}$, $\tg x=\f{1}{8192}$ for the BD method, and $\tg
t=\f{1}{100000}$, $\tg x=\f{1}{131072}$ for the ``exact''
solution. }\label{fig19}\vspace{-3mm}
\[\Delta^{\rm ts}_\infty(t)=1.23\E-1,\
\Delta^{\rm bd}_\infty(t)=1.20\E-1,\
\Delta^{\rm ts}_2(t)=2.29\E-2,\
\Delta^{\rm bd}_2(t)= 2.31\E-2.\]
\end{figure}

\begin{figure} \footnotesize
\begin{center}
\resizebox{1.5in}{!}{\includegraphics{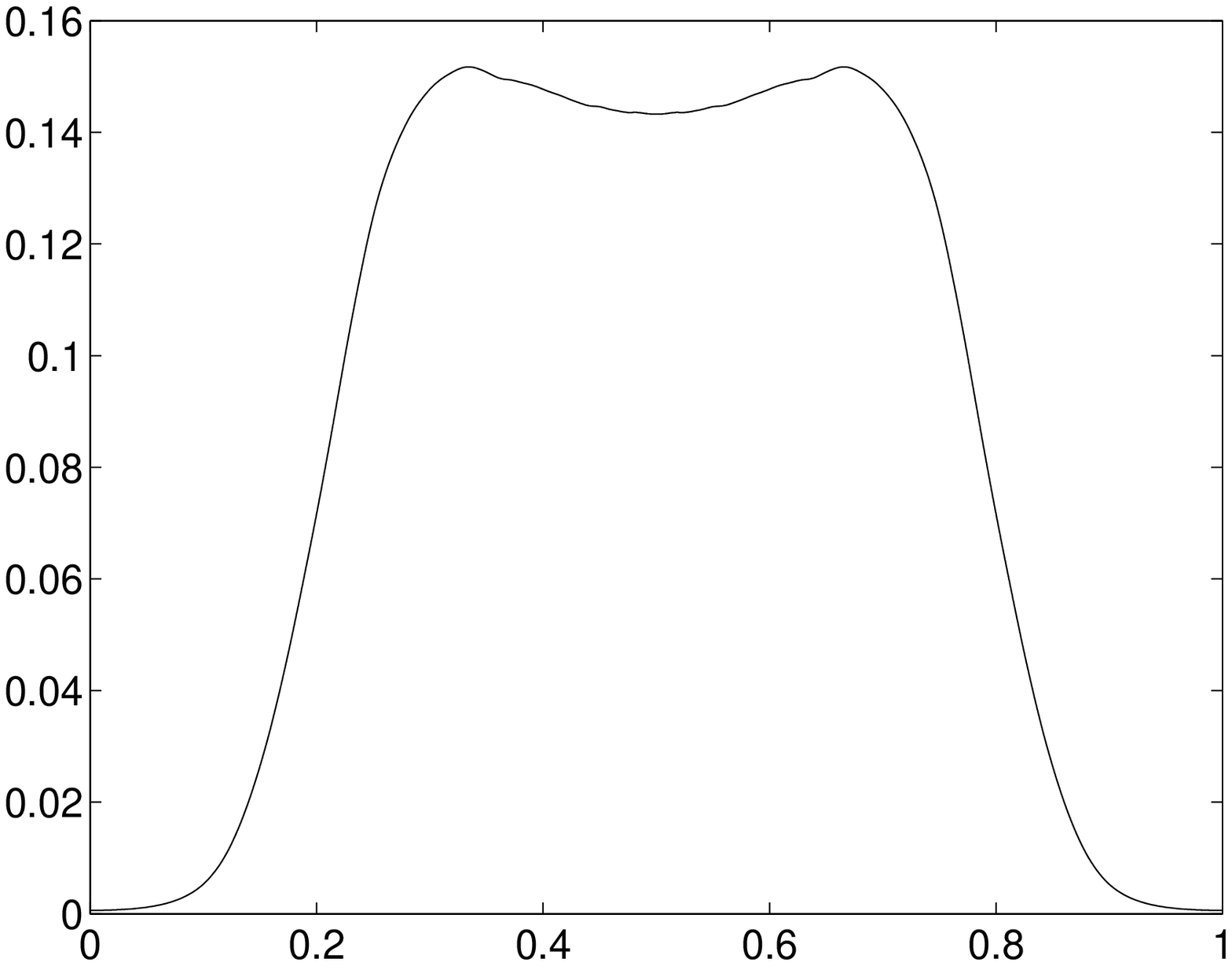}}
\resizebox{1.5in}{!}{\includegraphics{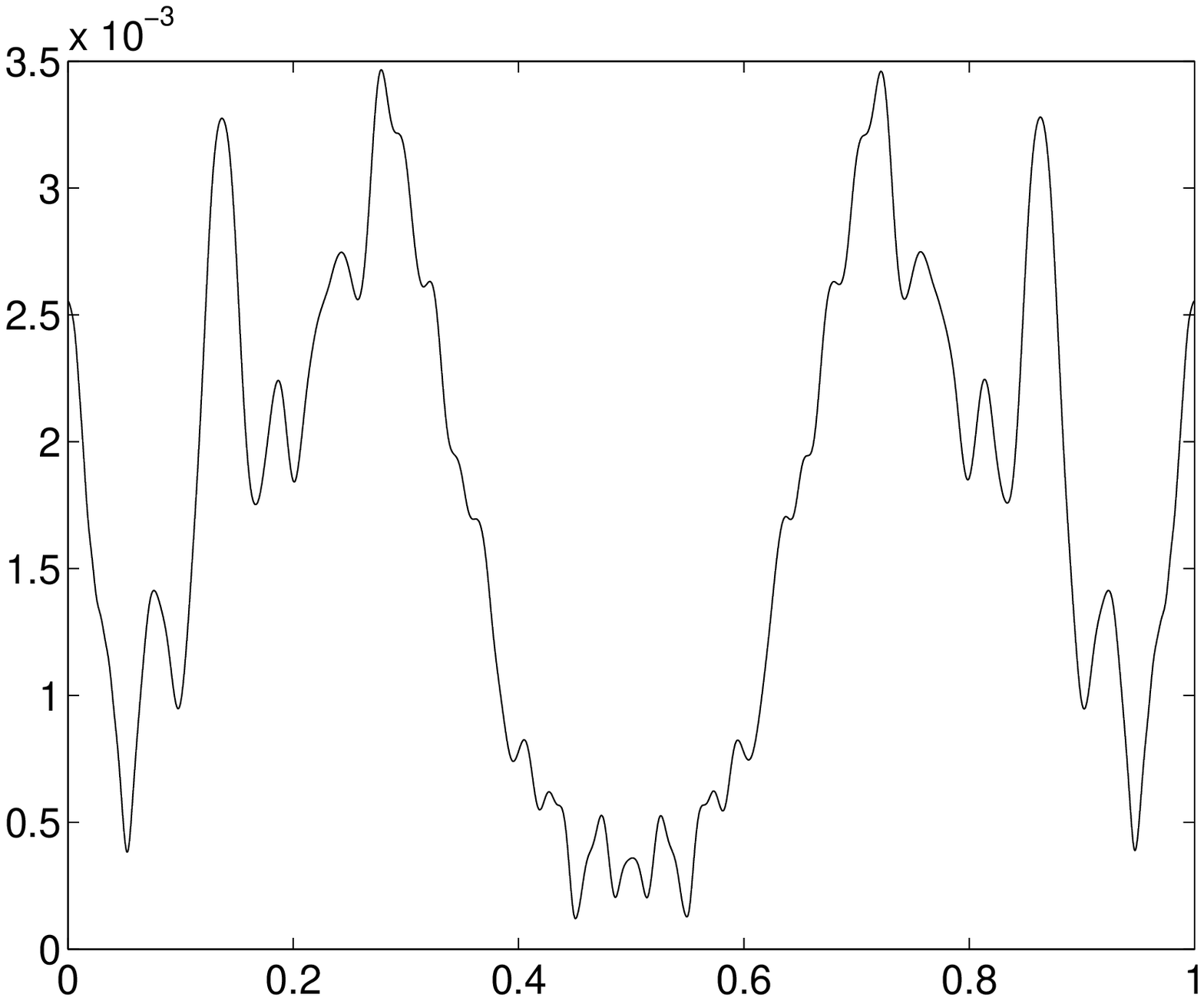}}
\resizebox{1.5in}{!}{\includegraphics{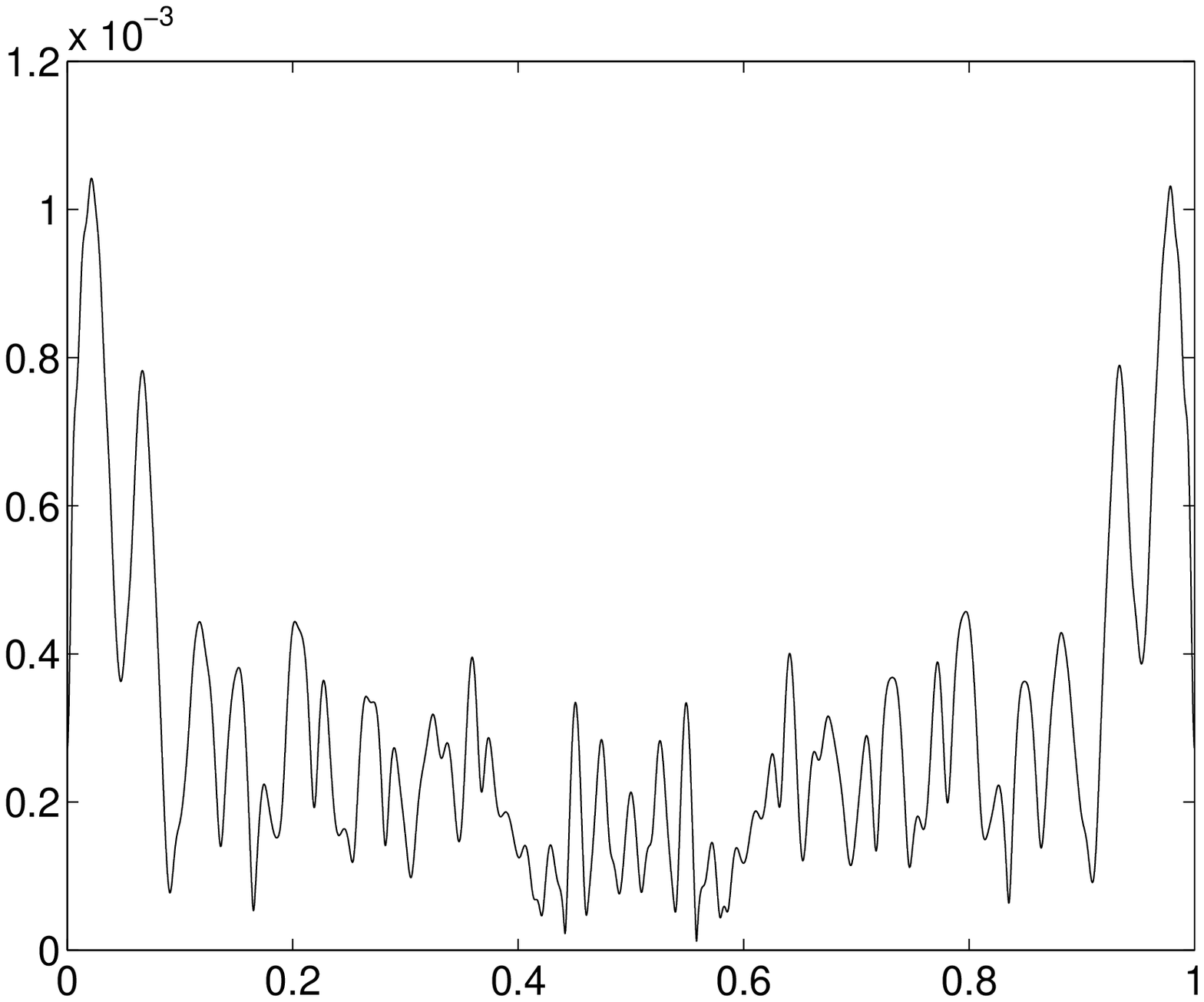}}

$|\psi^{\rm ex}(t,x)|^2$, $|\psi^{\rm ex}(t,x)-\psi^{\rm
ts}(t,x)|$ and $|\psi^{\rm ex}(t,x)-\psi^{\rm bd}(t,x)|$ at
$t=1.0$\vspace{-1mm}
\caption{Numerical results for example
\ref{nuex1} with $U(x)$ given by \eqref{eqhar_p} and
$\e=\f{1}{2}$. We use $\tg t=\f{1}{10}$, $\tg x=\f{1}{32}$ for the
TS and the BD method, and $\tg t=\f{1}{100000}$, $\tg
x=\f{1}{8192}$ for the ``exact'' solution.
}\label{fig14}\vspace{-3mm}
\[\Delta^{\rm ts}_\infty(t)=3.47\E-3,\
\Delta^{\rm bd}_\infty(t)=1.04\E-3,\
\Delta^{\rm ts}_2(t)=1.96\E-3,\
\Delta^{\rm bd}_2(t)= 3.65\E-4.\]

\resizebox{1.5in}{!}{\includegraphics{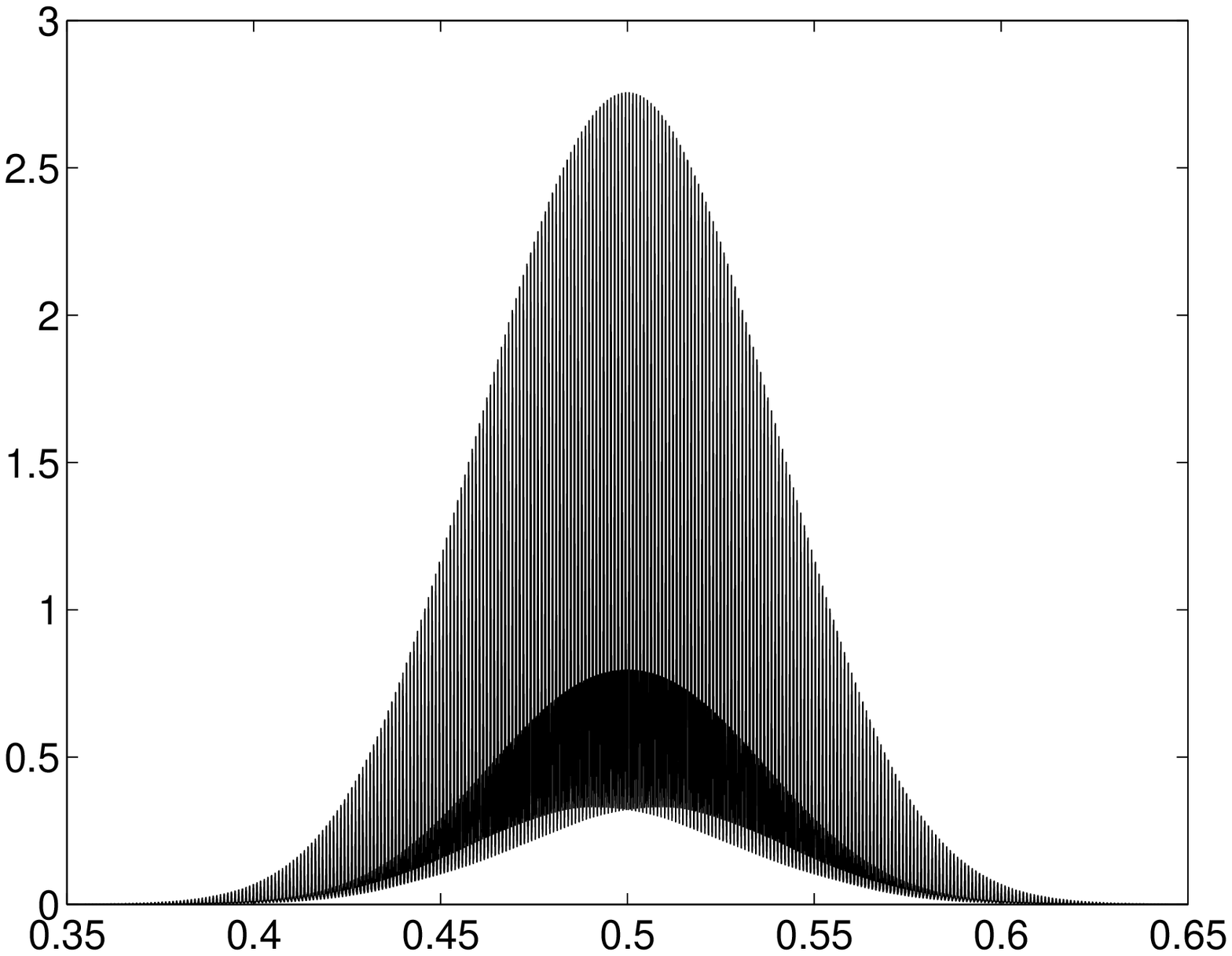}}
\resizebox{1.5in}{!}{\includegraphics{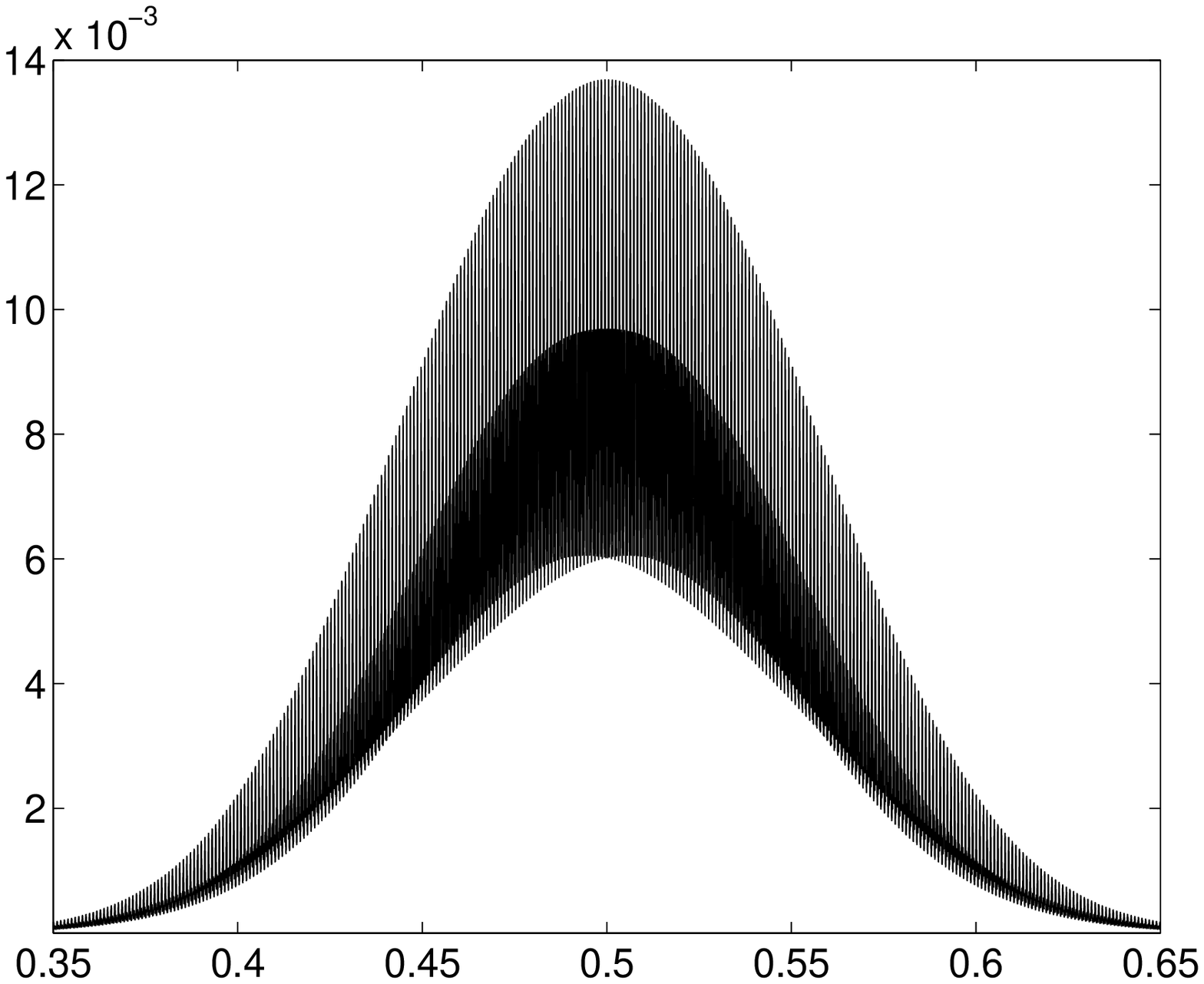}}
\resizebox{1.5in}{!}{\includegraphics{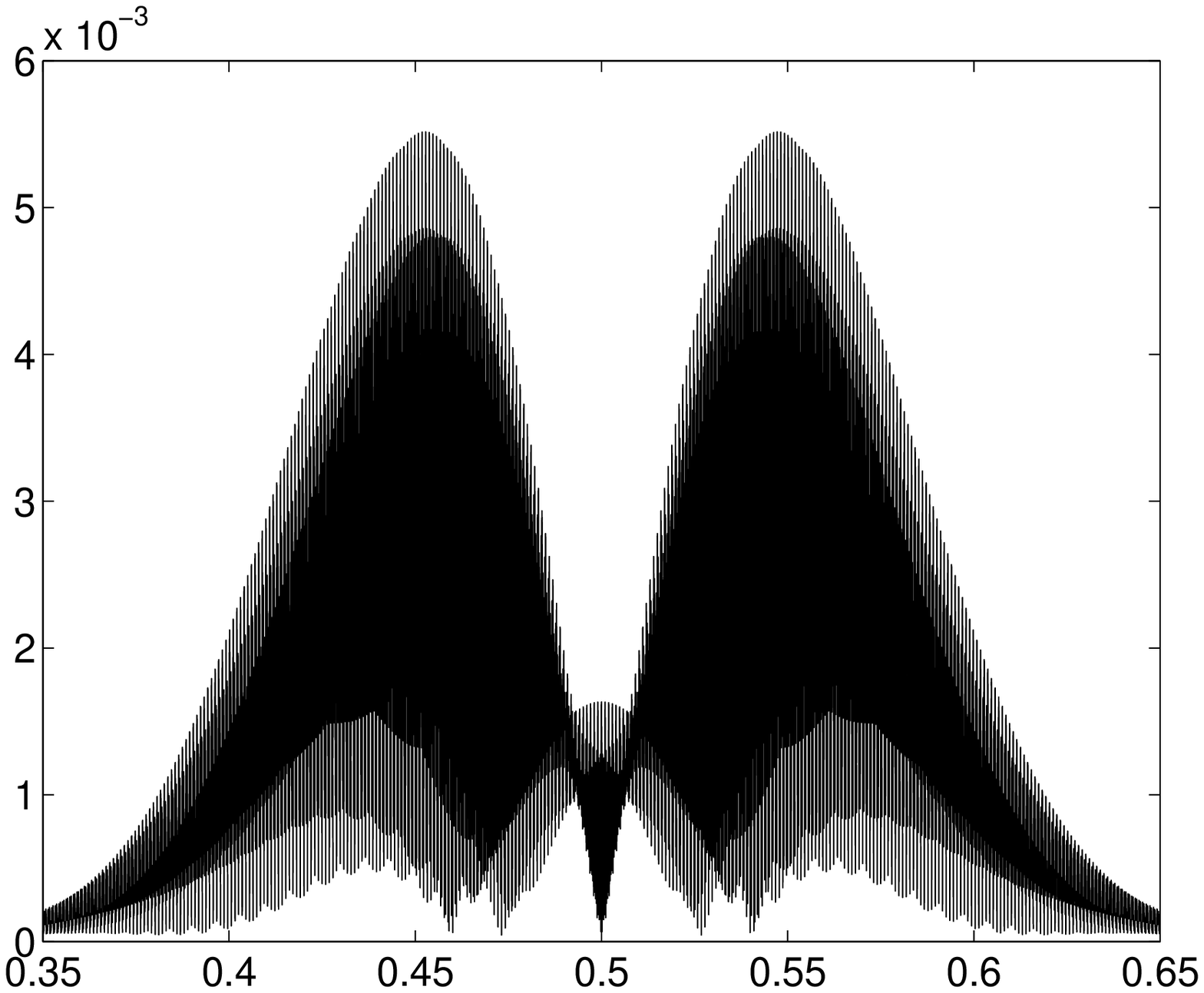}}

$|\psi^{\rm ex}(t,x)|^2$, $|\psi^{\rm ex}(t,x)-\psi^{\rm ts}(t,x)|$
and $|\psi^{\rm ex}(t,x)-\psi^{\rm bd}(t,x)|$ at $t=0.1$
\end{center}\vspace{-1mm}
\caption{Numerical results for example \ref{nuex1} with $U(x)$
given by \eqref{eqhar_p} and $\e=\f{1}{1024}$. We use Here $\tg
t=\f{1}{10000}$, $\tg x=\f{1}{16384}$ for the TS and $\tg
t=\f{1}{100}$, $\tg x=\f{1}{16384}$ for the BD method, and $\tg
t=\f{1}{100000}$, $\tg x=\f{1}{131072}$ for the ``exact''
solution. }\label{fig15}\vspace{-3mm}
\[\Delta^{\rm ts}_\infty(t)=1.37\E-2,\
\Delta^{\rm bd}_\infty(t)=5.52\E-3,\
\Delta^{\rm ts}_2(t)=2.76\E-3,\
\Delta^{\rm bd}_2(t)=1.20\E-3.\]
\end{figure}

\begin{figure} \footnotesize
\begin{center}
\resizebox{1.5in}{!}{\includegraphics{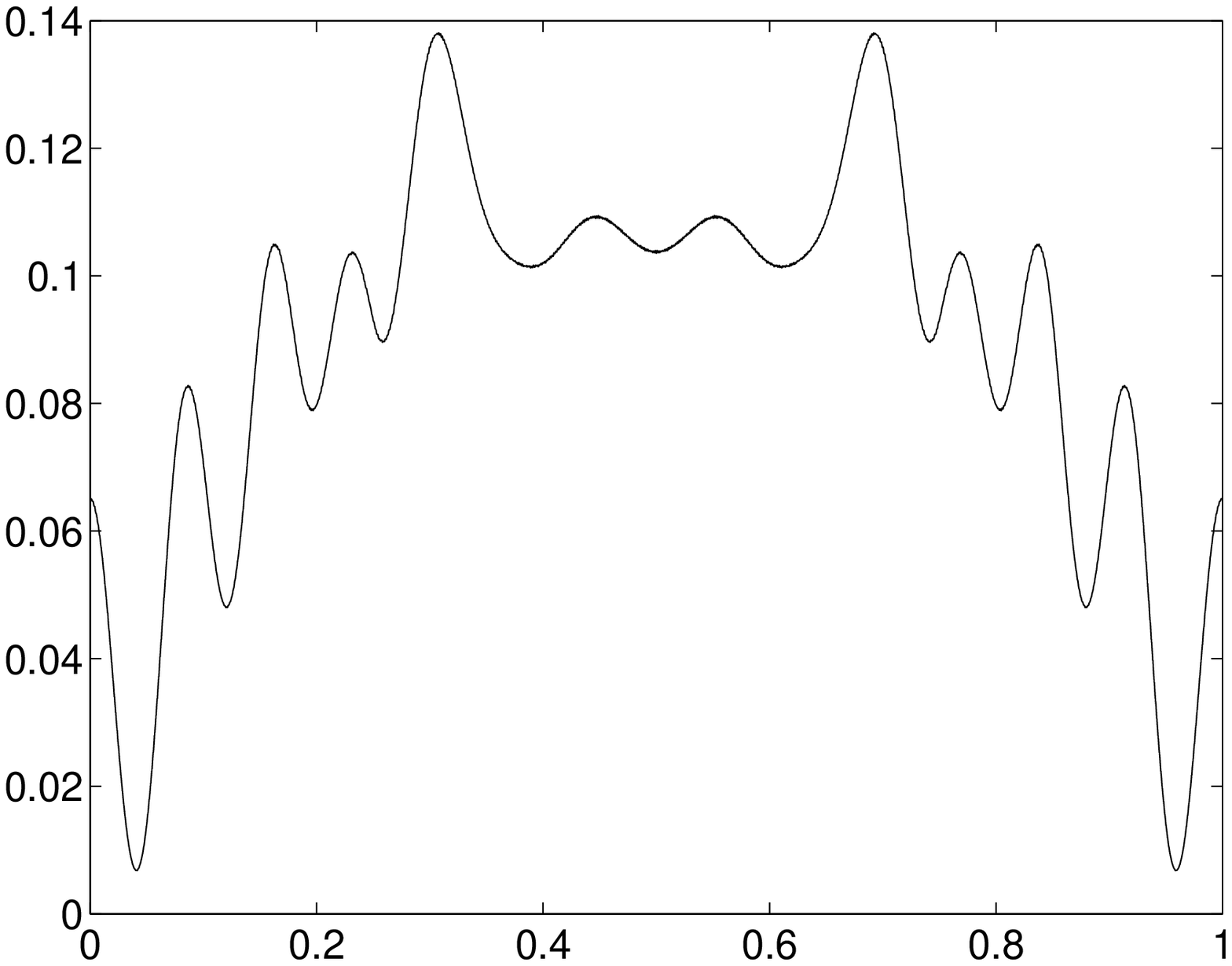}}
\resizebox{1.5in}{!}{\includegraphics{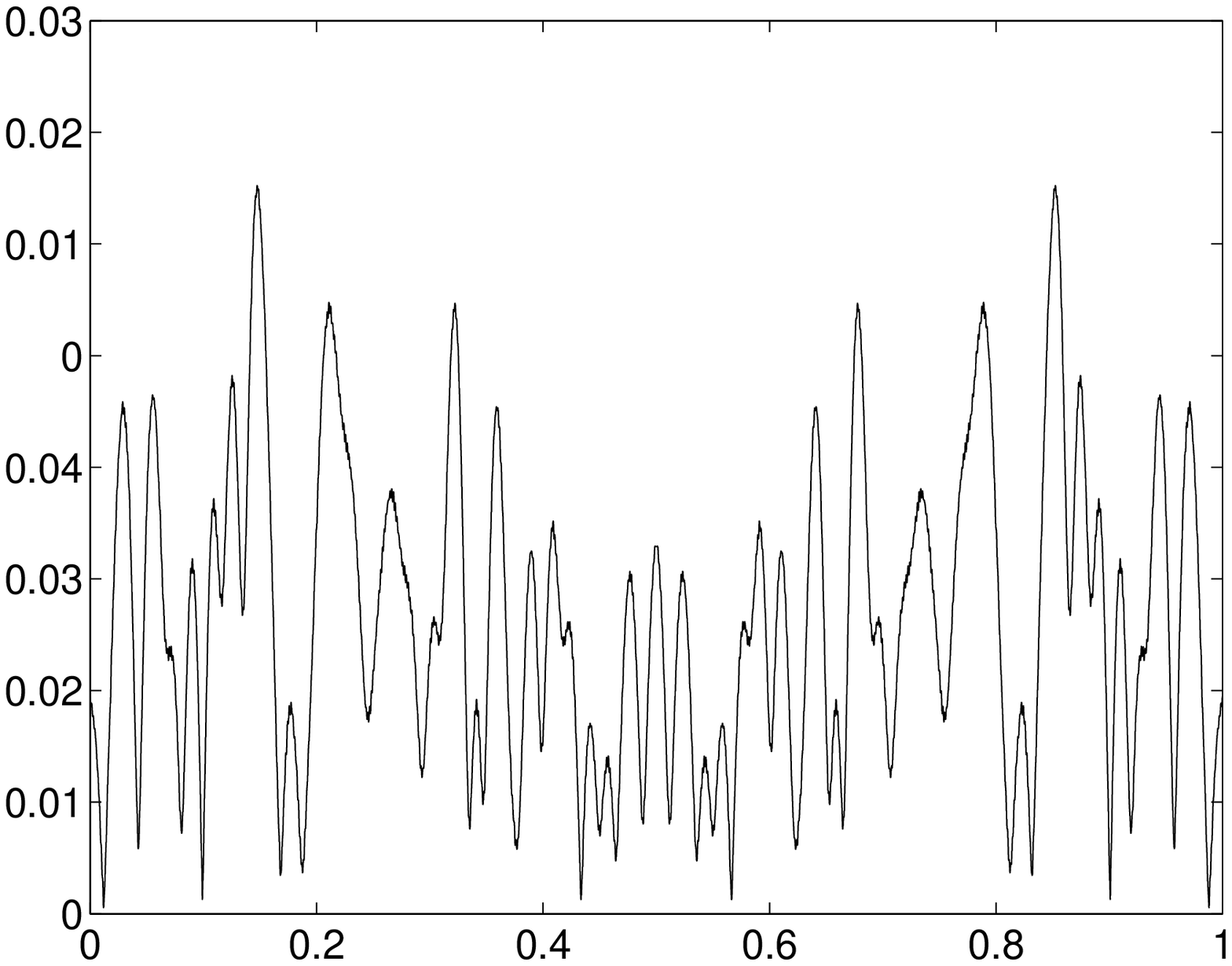}}
\resizebox{1.5in}{!}{\includegraphics{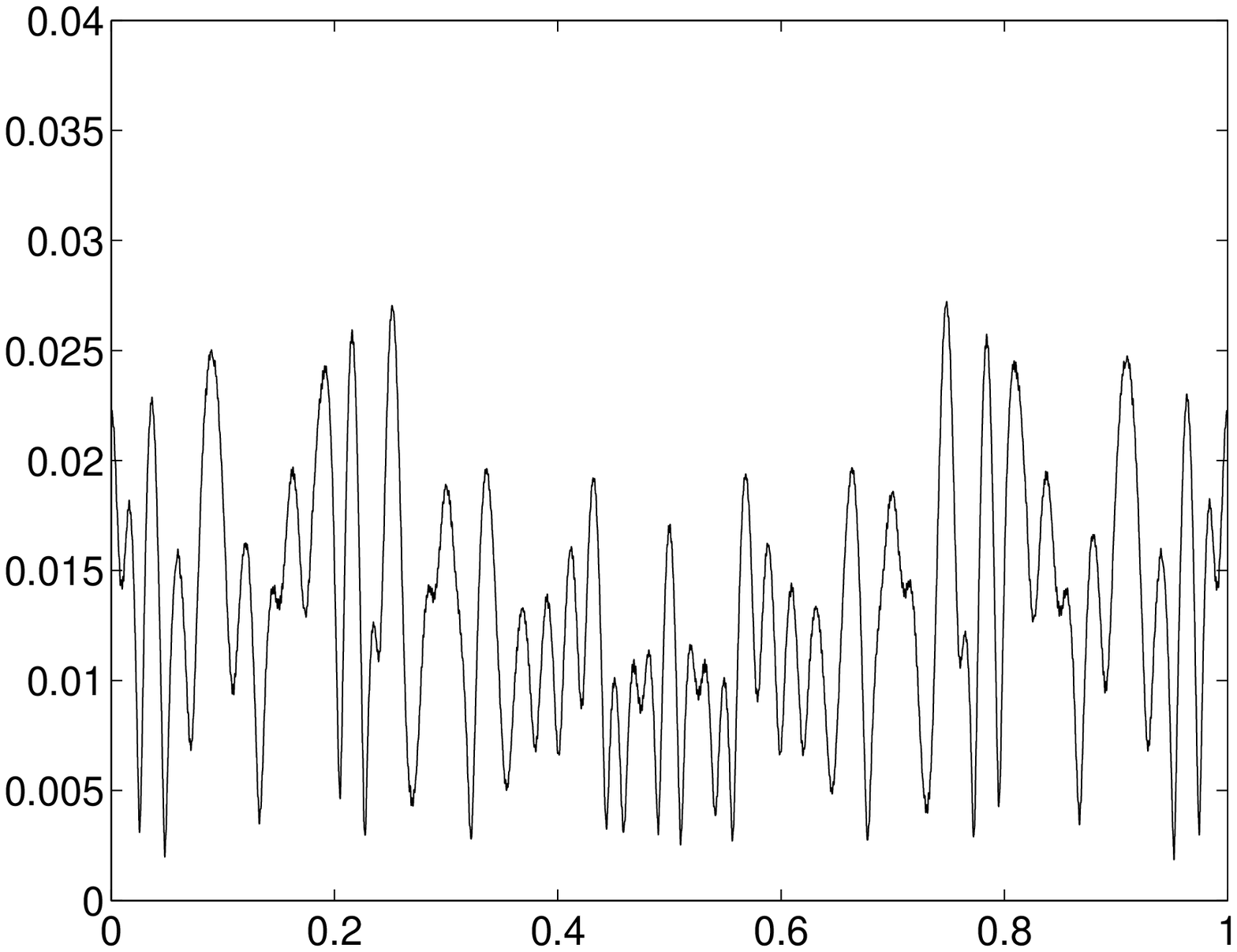}}

$|\psi^{\rm ex}(t,x)|^2$, $|\psi^{\rm ex}(t,x)-\psi^{\rm
ts}(t,x)|$ and $|\psi^{\rm ex}(t,x)-\psi^{\rm bd}(t,x)|$ at
$t=1.0$\vspace{-1mm} \caption{Numerical results for example
\ref{nuex1} with $U(x)$ given by \eqref{eqstep} and $\e=\f{1}{2}$.
We use $\tg t=\f{1}{10}$, $\tg x=\f{1}{32}$ for the TS and the BD
method, and $\tg t=\f{1}{100000}$, $\tg x=\f{1}{8192}$ for the
``exact'' solution. \vspace{-3mm}}\label{fig16}
\[\Delta^{\rm ts}_\infty(t)=3.26\E-2,\
\Delta^{\rm bd}_\infty(t)=2.72\E-2,\
\Delta^{\rm ts}_2(t)=1.51\E-2,\
\Delta^{\rm bd}_2(t)=1.45\E-2.\]\vspace{1mm}

\resizebox{1.5in}{!}{\includegraphics{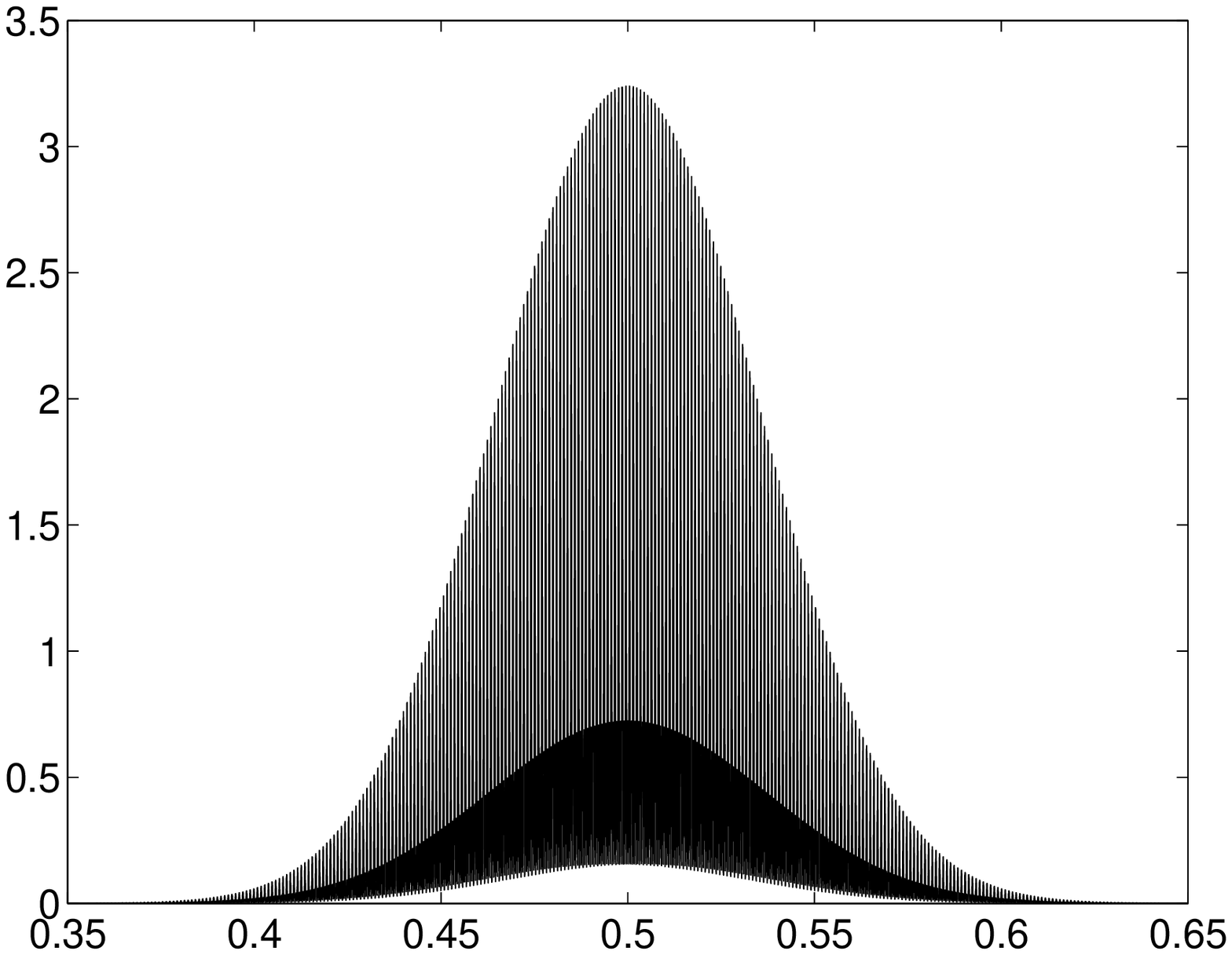}}
\resizebox{1.54in}{!}{\includegraphics{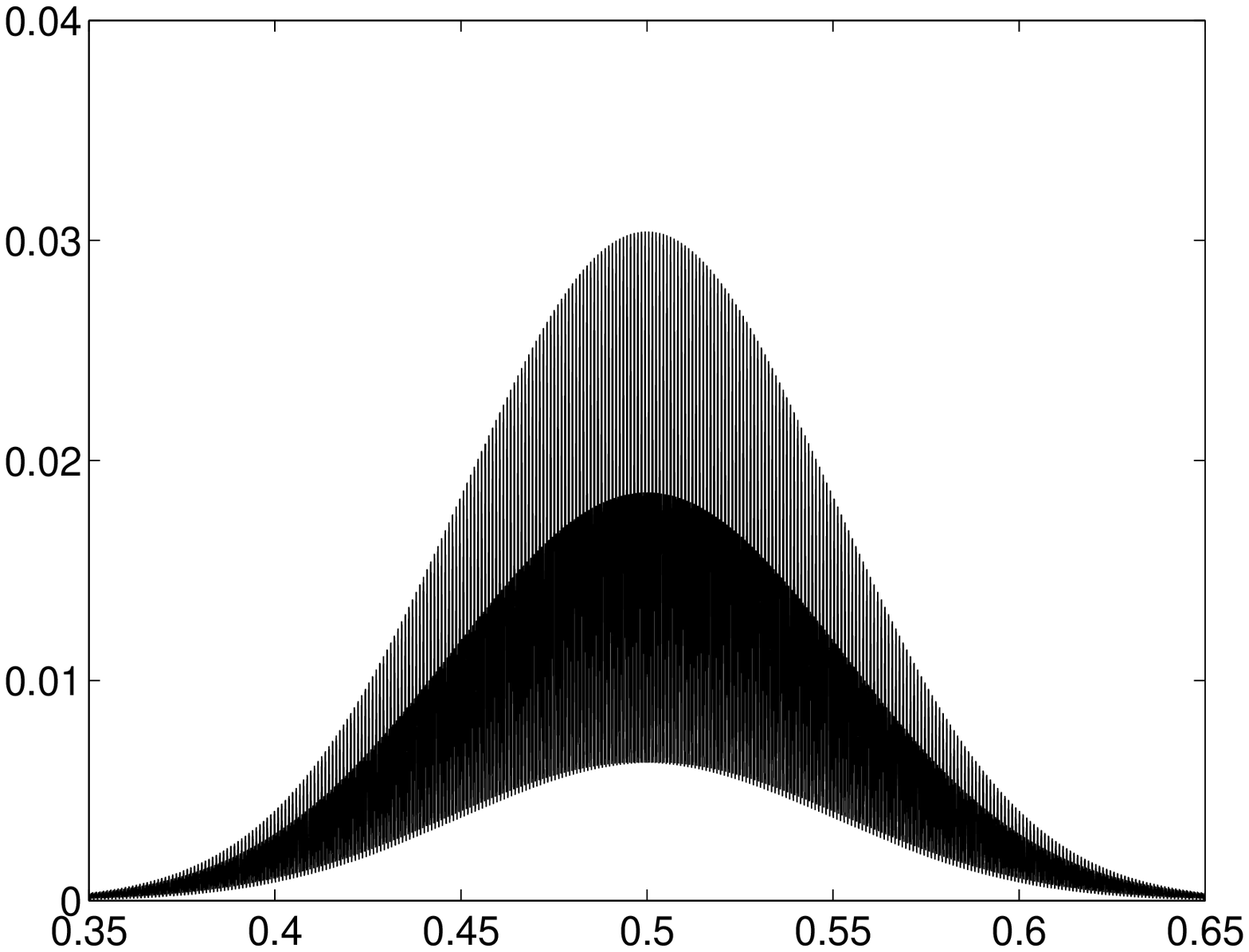}}
\resizebox{1.5in}{!}{\includegraphics{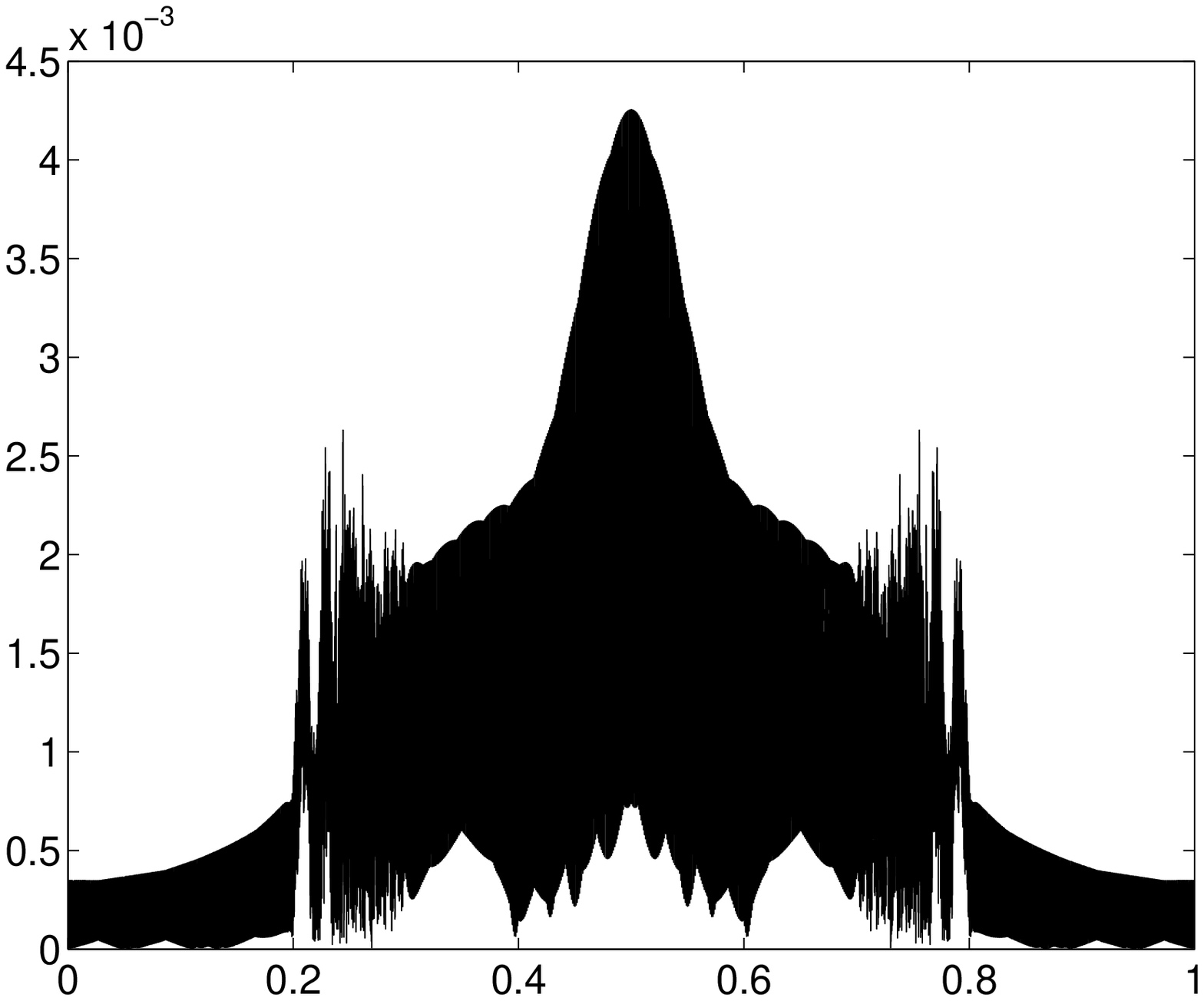}}

$|\psi^{\rm ex}(t,x)|^2$, $|\psi^{\rm ex}(t,x)-\psi^{\rm ts}(t,x)|$
and $|\psi^{\rm ex}(t,x)-\psi^{\rm bd}(t,x)|$ at $t=0.1$\vspace{-1mm}
\end{center}
\caption{Numerical results for example \ref{nuex1} with $U(x)$
given by \eqref{eqstep} and $\e=\f{1}{1024}$. We use $\tg
t=\f{1}{10000}$, $\tg x=\f{1}{16384}$ for the TS and $\tg
t=\f{1}{10}$, $\tg x=\f{1}{8192}$ for the BD method, and $\tg
t=\f{1}{100000}$, $\tg x=\f{1}{131072}$ for the ``exact''
solution. }\label{fig17}\vspace{-3mm}
\[\Delta^{\rm ts}_\infty(t)=3.04\E-2,\
\Delta^{\rm bd}_\infty(t)=4.25\E-3,\
\Delta^{\rm ts}_2(t)=5.35\E-3,\
\Delta^{\rm bd}_2(t)= 1.21\E-3.\]
\end{figure}

\begin{figure} \footnotesize
\begin{center}
\resizebox{1.5in}{!}{\includegraphics{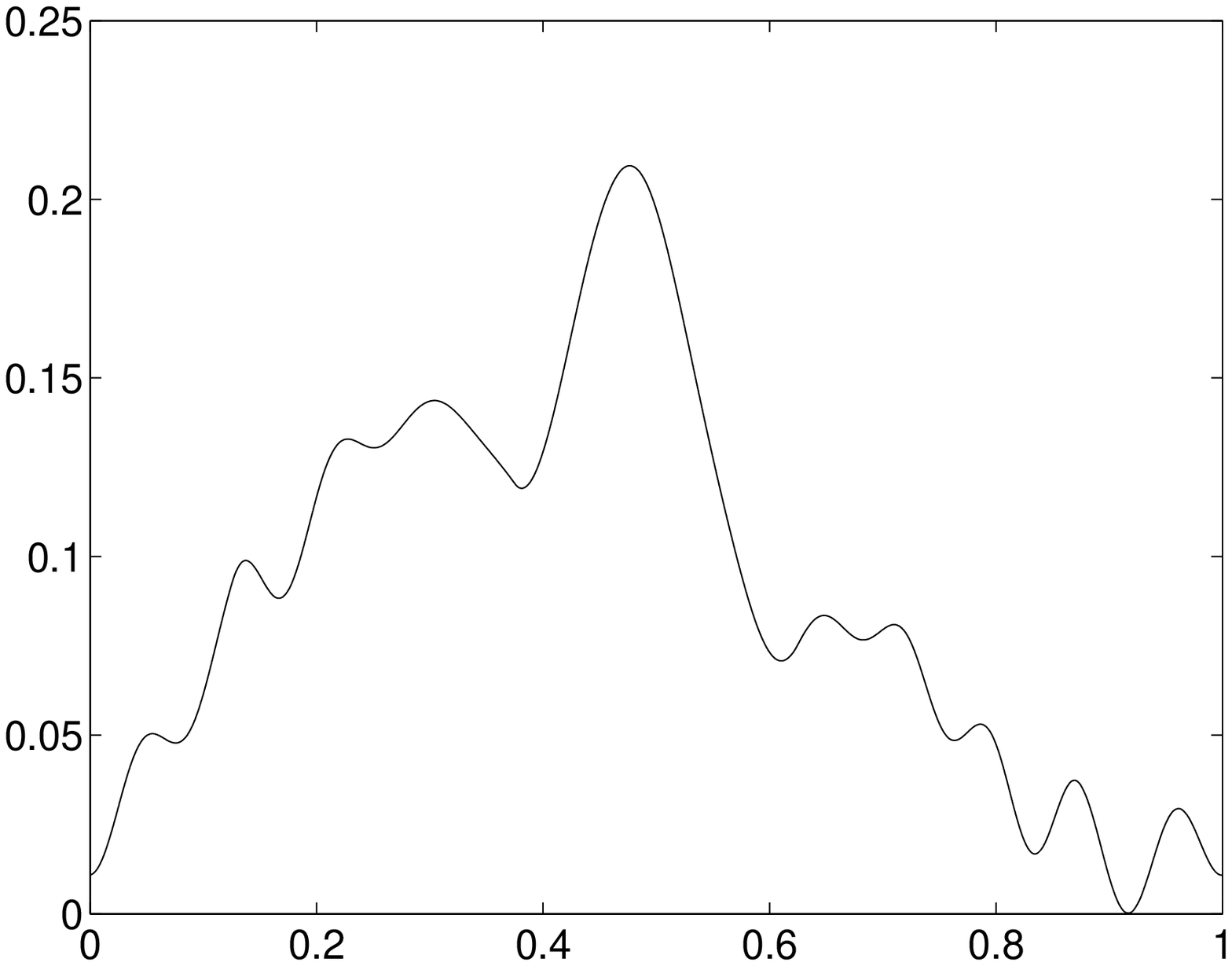}}
\resizebox{1.5in}{!}{\includegraphics{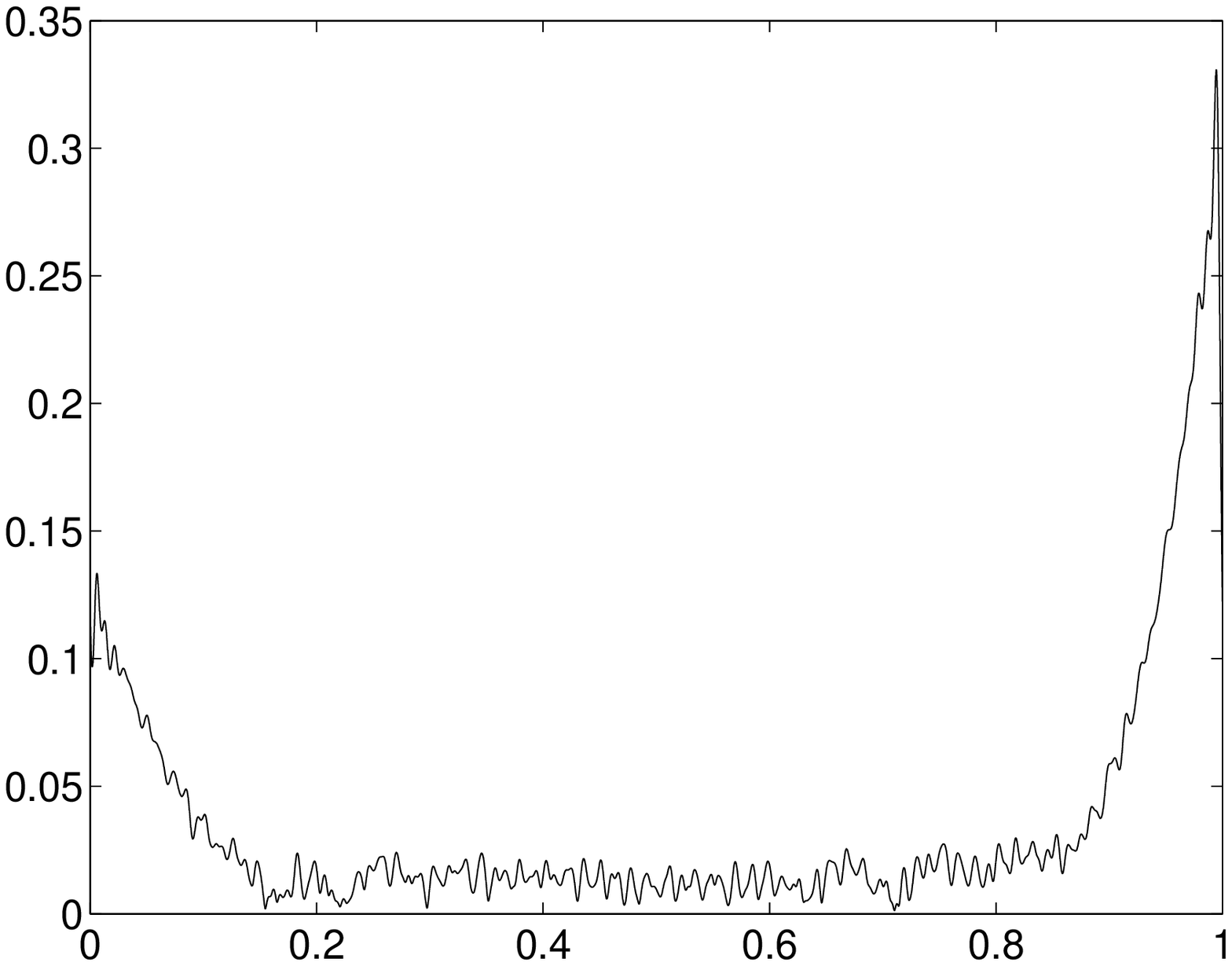}}
\resizebox{1.5in}{!}{\includegraphics{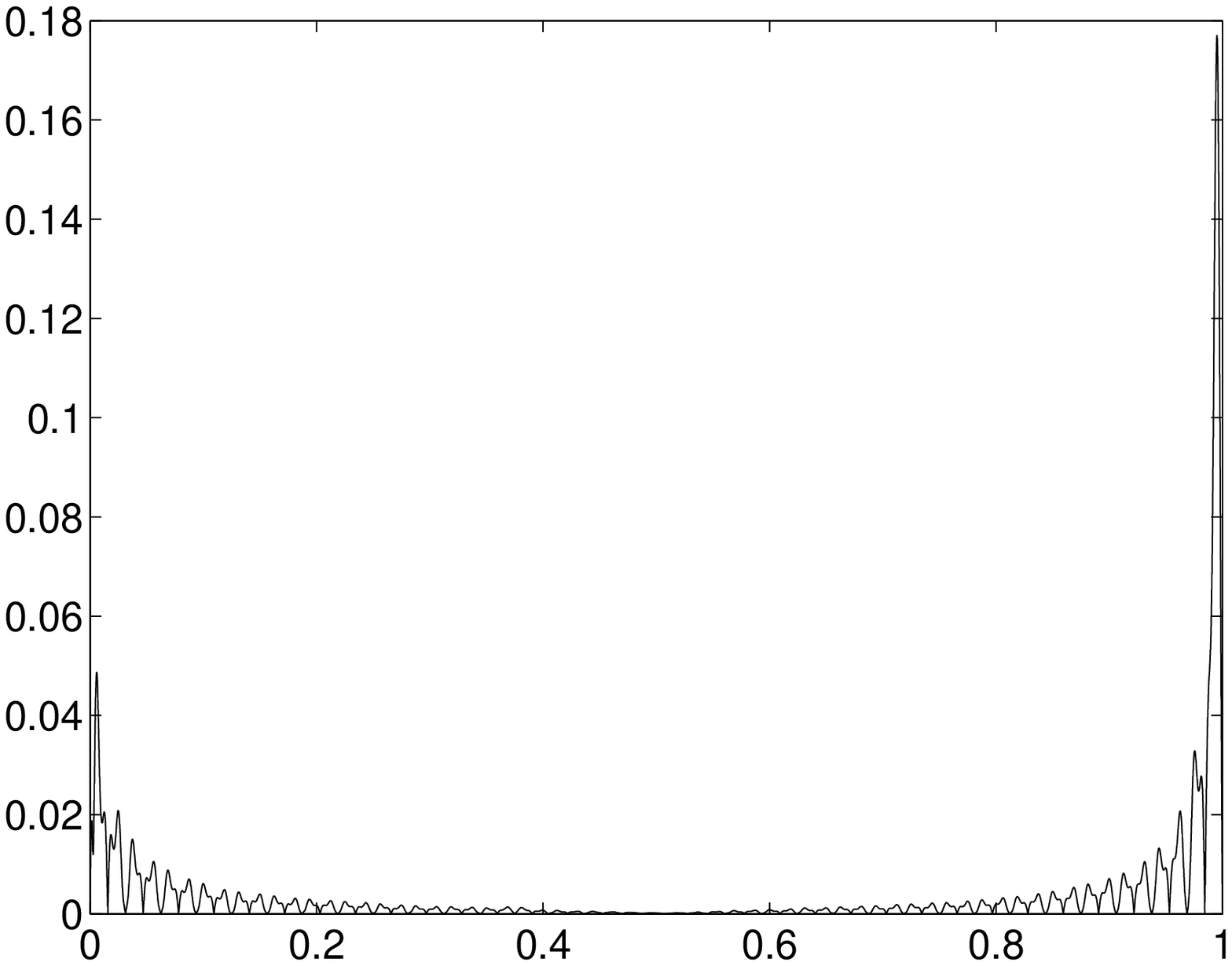}}

$|\psi^{\rm ex}(t,x)|^2$, $|\psi^{\rm ex}(t,x)-\psi^{\rm
ts}(t,x)|$ and $|\psi^{\rm ex}(t,x)-\psi^{\rm bd}(t,x)|$ at
$t=1.0$\vspace{-1mm} \caption{Numerical results for example
\ref{nuex2} with $U(x)$ given by \eqref{eqlin}, $\e=\f{1}{2}$. We
use $\tg t=\f{1}{100}$, $\tg x=\f{1}{64}$ for the TS, $\tg
t=\f{1}{2}$, $\tg x=\f{1}{32}$ for BD method, and $\tg
t=\f{1}{100000}$, $\tg x=\f{1}{8192}$ for the ``exact'' solution.
\vspace{-3mm}}\label{fig27}
\[\Delta^{\rm ts}_\infty(t)=3.31\E-1,\
\Delta^{\rm bd}_\infty(t)=1.77\E-1,\
\Delta^{\rm ts}_2(t)=6.16\E-2,\
\Delta^{\rm bd}_2(t)=1.38\E-2.\]\vspace{1mm}

\resizebox{1.5in}{!}{\includegraphics{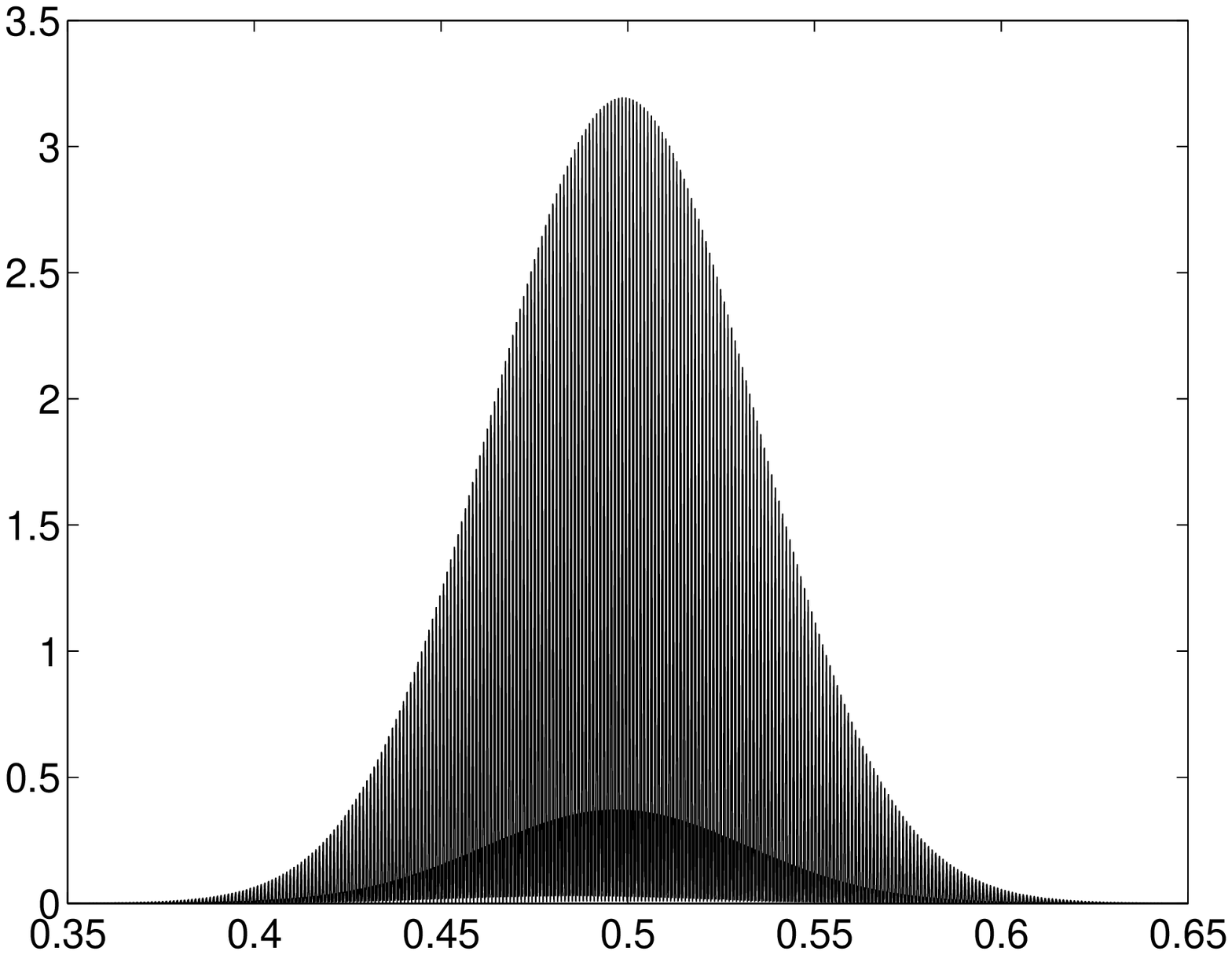}}
\resizebox{1.5in}{!}{\includegraphics{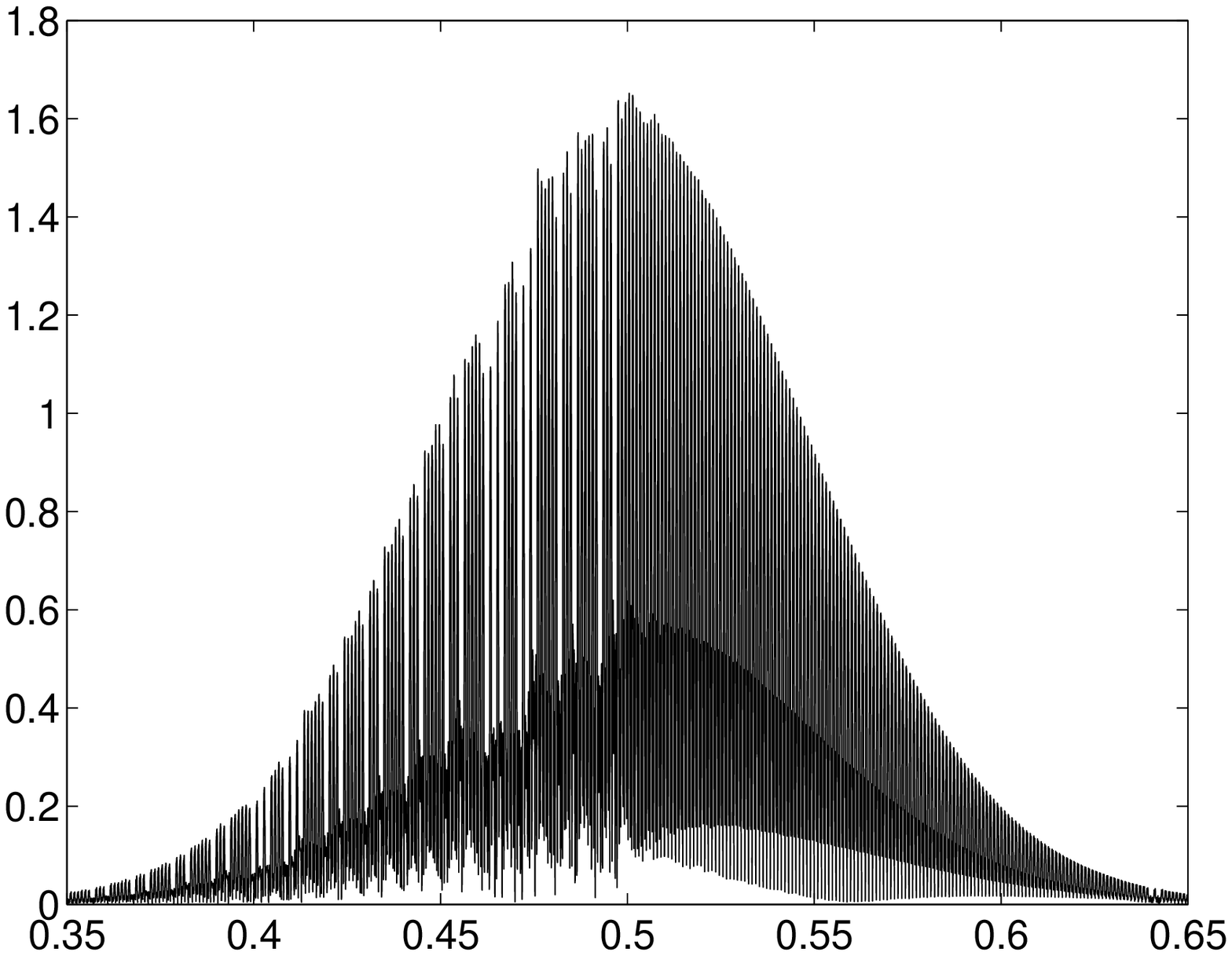}}
\resizebox{1.5in}{!}{\includegraphics{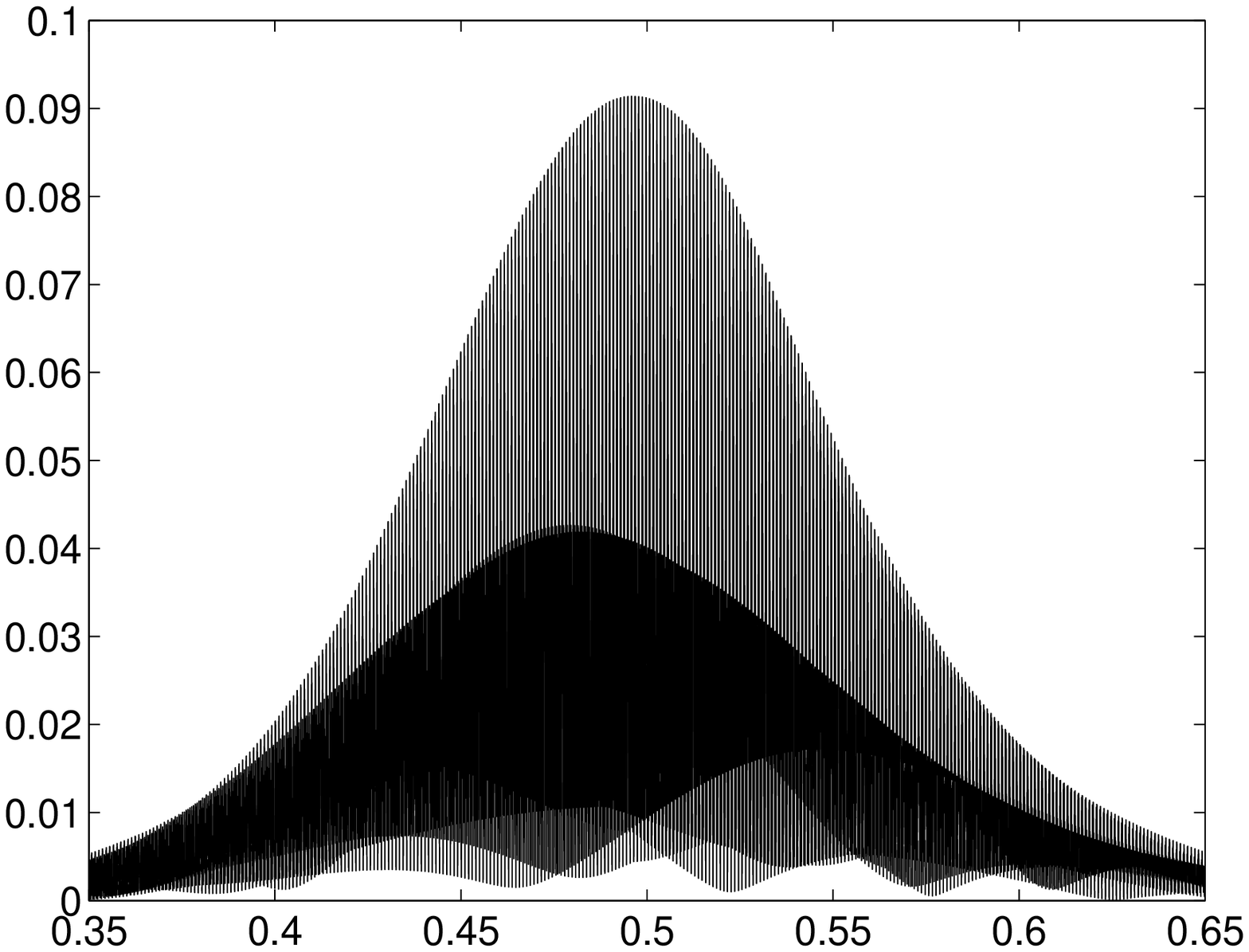}}

$|\psi^{\rm ex}(t,x)|^2$, $|\psi^{\rm ex}(t,x)-\psi^{\rm ts}(t,x)|$
and $|\psi^{\rm ex}(t,x)-\psi^{\rm bd}(t,x)|$ at $t=0.1$\vspace{-1mm}
\end{center}
\caption{Numerical results for example \ref{nuex2} with $U(x)$
given by \eqref{eqlin}, $\e=\f{1}{1024}$. We use $\tg
t=\f{1}{10000}$, $\tg x=\f{1}{65536}$ for the TS, $\tg
t=\f{1}{10}$, $\tg x=\f{1}{8192}$ for the BD method, and $\tg
t=\f{1}{100000}$, $\tg x=\f{1}{131072}$ for the ``exact''
solution. }\label{fig29}\vspace{-3mm}
\[\Delta^{\rm ts}_\infty(t)=1.65,\
\Delta^{\rm bd}_\infty(t)=9.14\E-2,\
\Delta^{\rm ts}_2(t)=2.63\E-1,\
\Delta^{\rm bd}_2(t)=1.39\E-2.\]
\end{figure}

\begin{figure} \footnotesize
\begin{center}
\resizebox{1.5in}{!}{\includegraphics{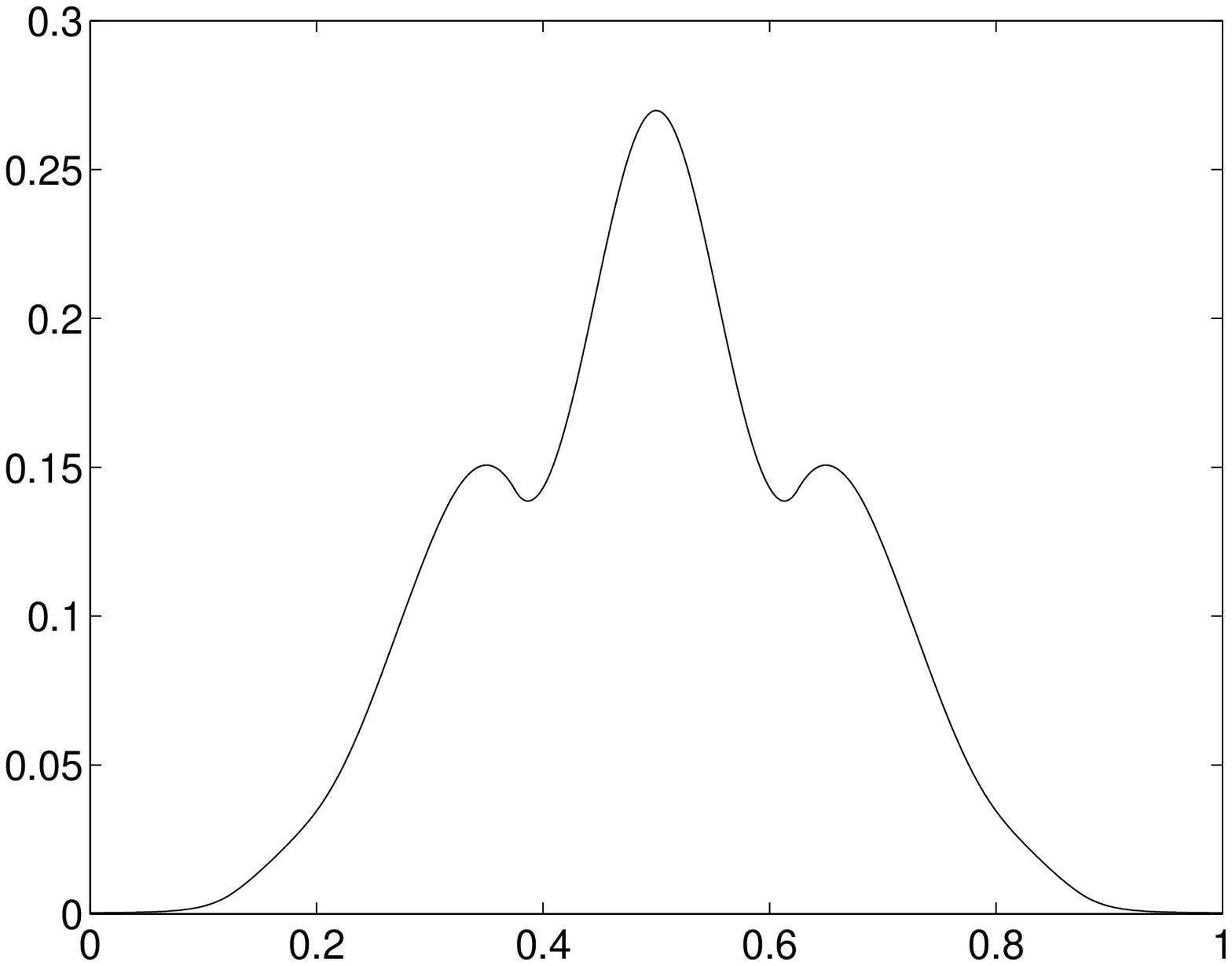}}
\resizebox{1.5in}{!}{\includegraphics{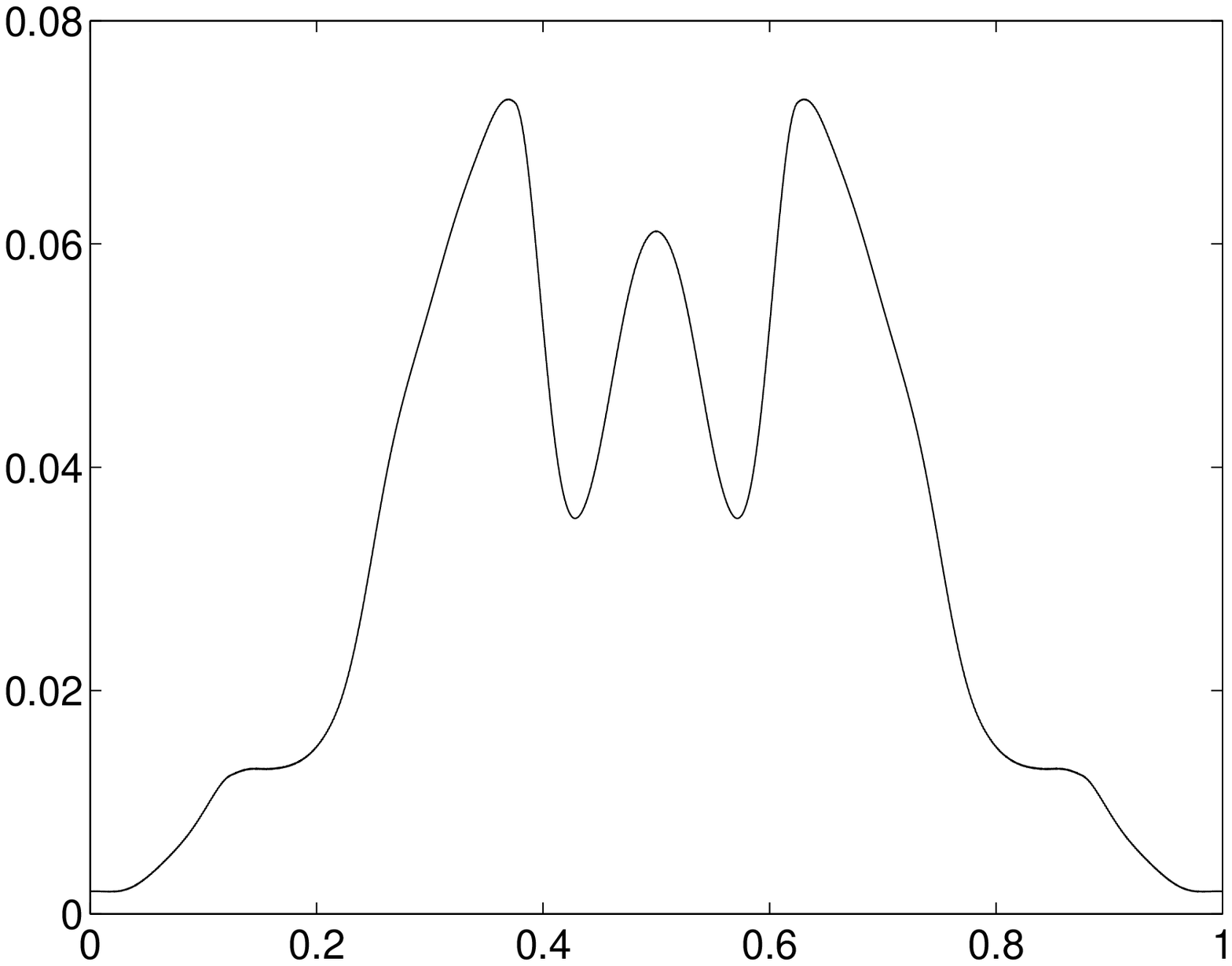}}
\resizebox{1.5in}{!}{\includegraphics{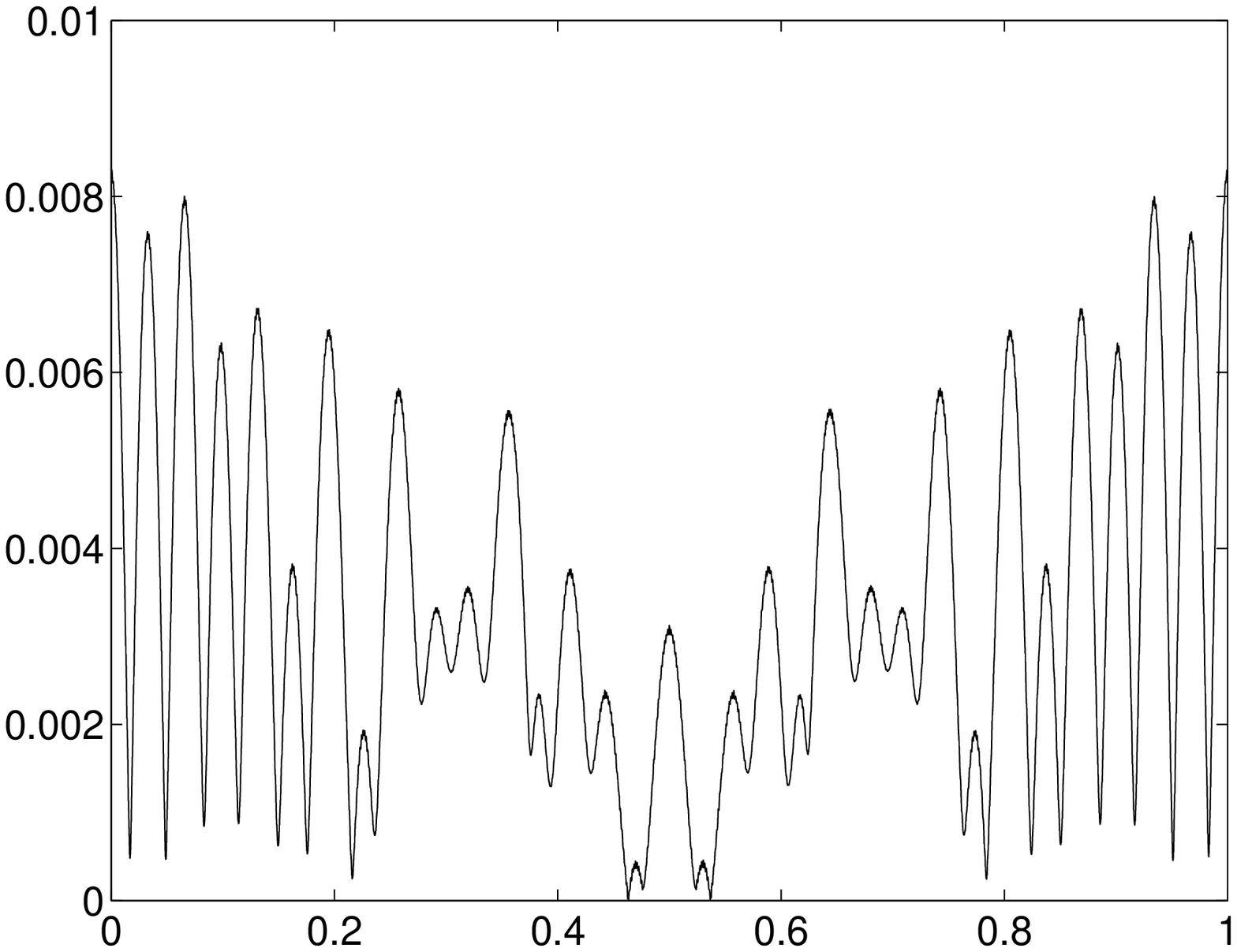}}

$|\psi^{\rm ex}(t,x)|^2$, $|\psi^{\rm ex}(t,x)-\psi^{\rm
ts}(t,x)|$ and $|\psi^{\rm ex}(t,x)-\psi^{\rm bd}(t,x)|$ at
$t=1.0$\vspace{-1mm} \caption{Numerical results for example
\ref{nuex2} with $U(x)$ given by \eqref{eqhar_p} and
$\e=\f{1}{2}$. We use $\tg t=\f{1}{200}$, $\tg x=\f{1}{64}$ for
the TS, $\tg t=\f{1}{5}$, $\tg x=\f{1}{32}$ for the BD method, and
$\tg t=\f{1}{100000}$, $\tg x=\f{1}{8192}$ for the ``exact''
solution. \vspace{-3mm}}\label{fig21}
\[\Delta^{\rm ts}_\infty(t)=7.30\E-2,\
\Delta^{\rm bd}_\infty(t)=8.30\E-3,\
\Delta^{\rm ts}_2(t)= 4.02\E-2,\
\Delta^{\rm bd}_2(t)= 3.89\E-3.\]\vspace{1mm}

\resizebox{1.5in}{!}{\includegraphics{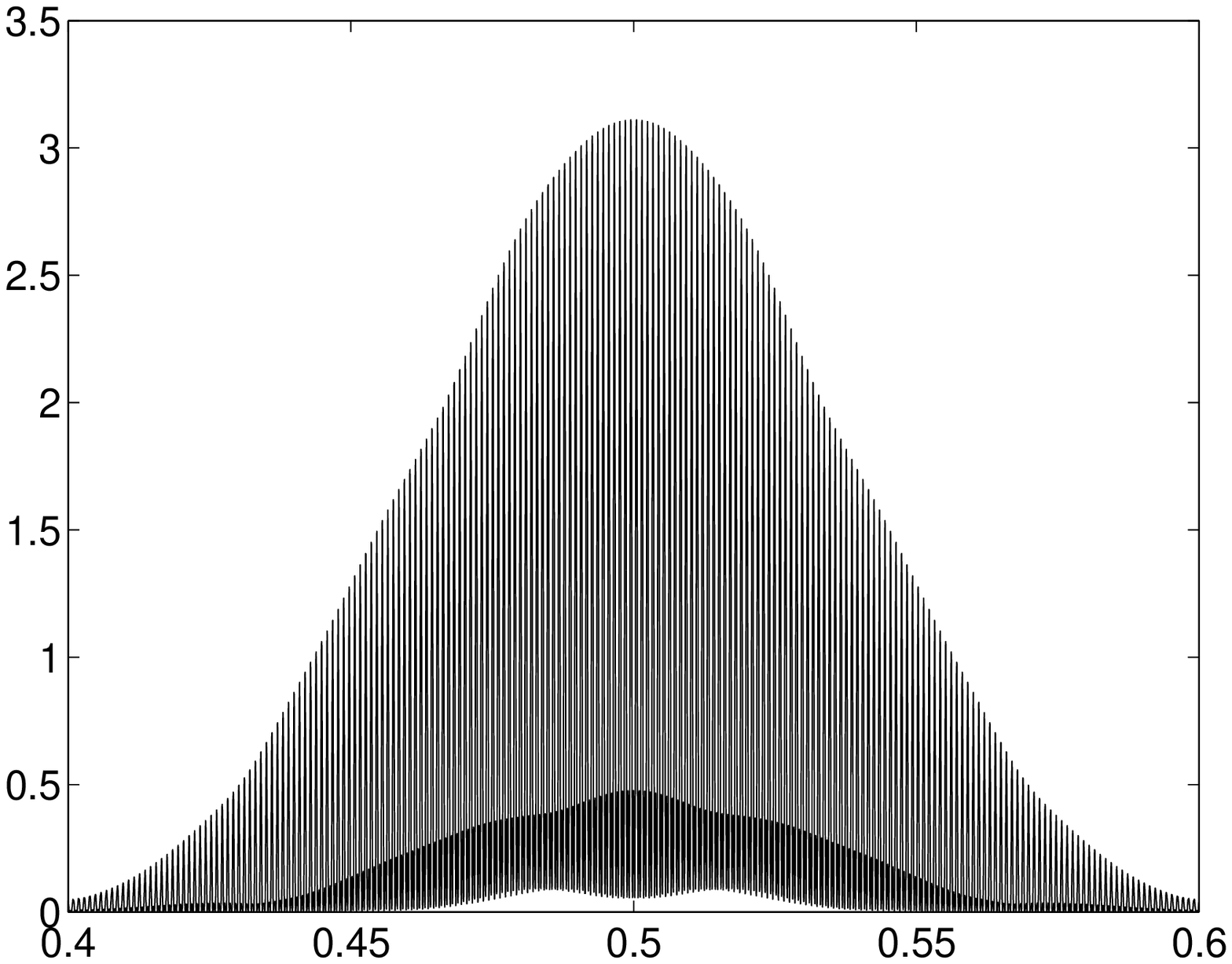}}
\resizebox{1.5in}{!}{\includegraphics{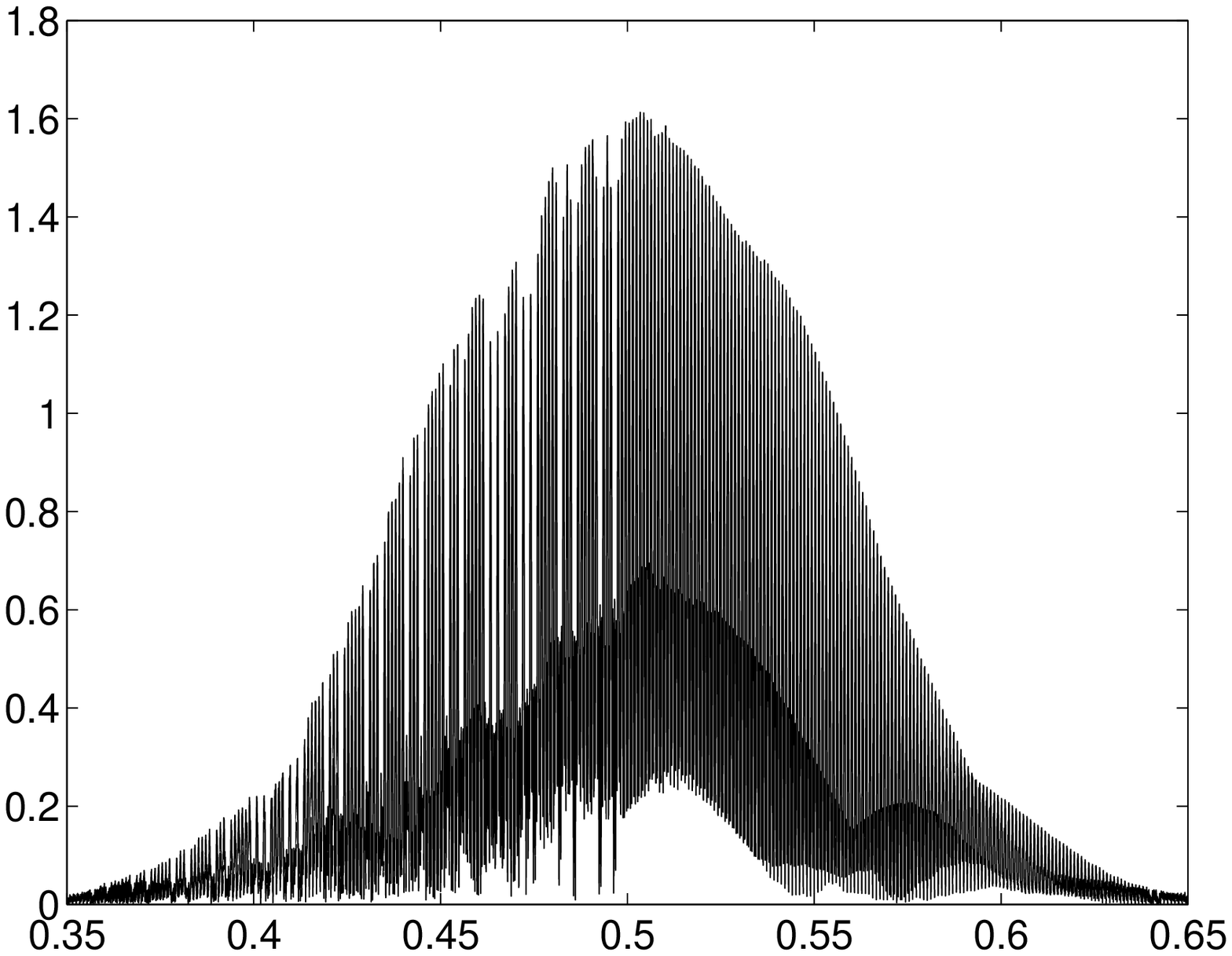}}
\resizebox{1.5in}{!}{\includegraphics{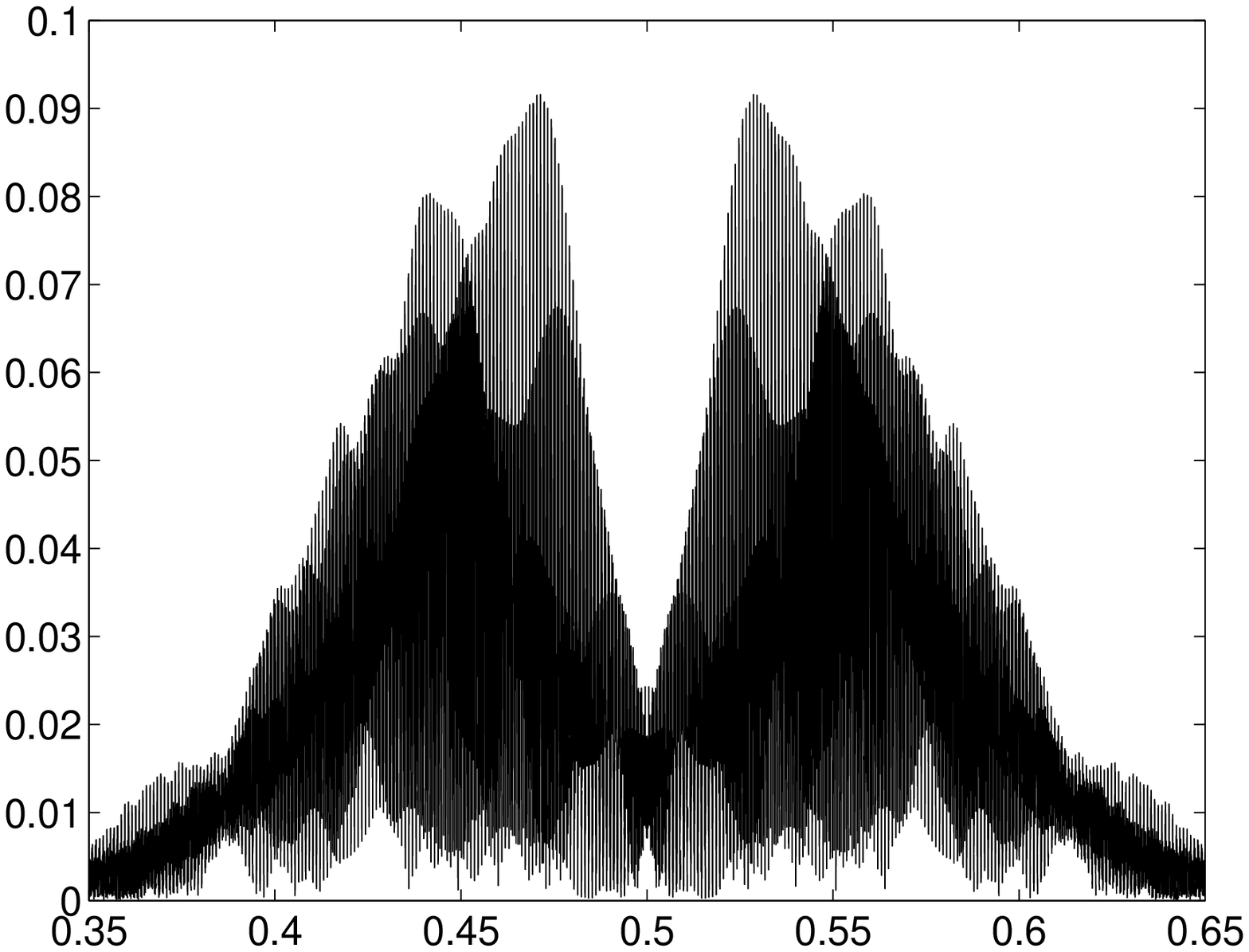}}

$|\psi^{\rm ex}(t,x)|^2$, $|\psi^{\rm ex}(t,x)-\psi^{\rm ts}(t,x)|$
and $|\psi^{\rm ex}(t,x)-\psi^{\rm bd}(t,x)|$ at $t=0.1$\vspace{-1mm}
\end{center}
\caption{Numerical results for example \ref{nuex2} with $U(x)$
given by \eqref{eqhar_p} and $\e=\f{1}{1024}$. We use $\tg
t=\f{1}{50000}$, $\tg x=\f{1}{65536}$ for the TS, $\tg
t=\f{1}{10}$, $\tg x=\f{1}{8192}$ for the BD method, and $\tg
t=\f{1}{100000}$, $\tg x=\f{1}{131072}$ for the ``exact''
solution. }\label{fig23}\vspace{-3mm}
\[\Delta^{\rm ts}_\infty(t)=1.61,\
\Delta^{\rm bd}_\infty(t)=9.16\E-2,\
\Delta^{\rm ts}_2(t)=2.63\E-1,\
\Delta^{\rm bd}_2(t)= 1.71\E-2.\]
\end{figure}

\begin{figure} \footnotesize
\begin{center}
\resizebox{1.5in}{!}{\includegraphics{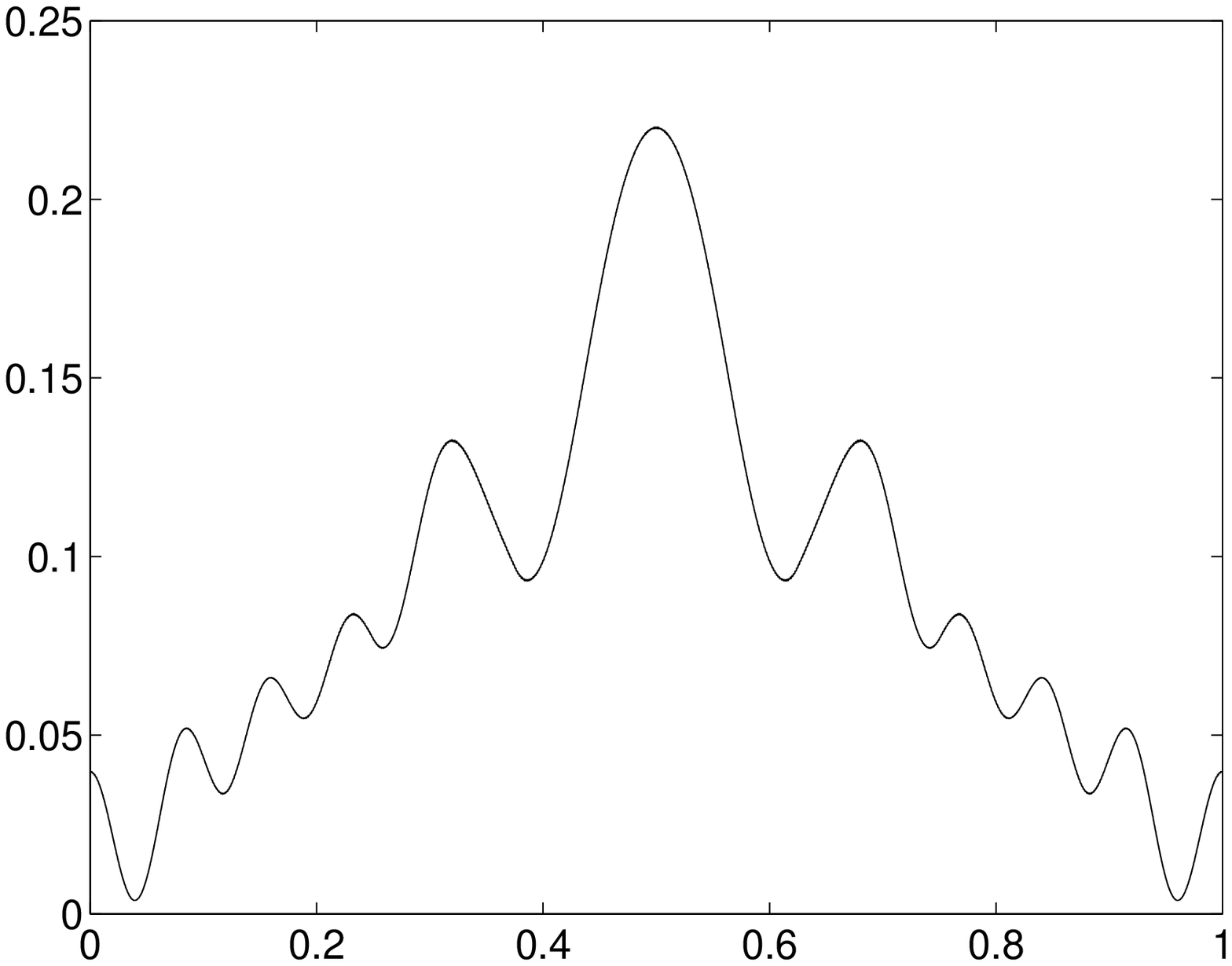}}
\resizebox{1.5in}{!}{\includegraphics{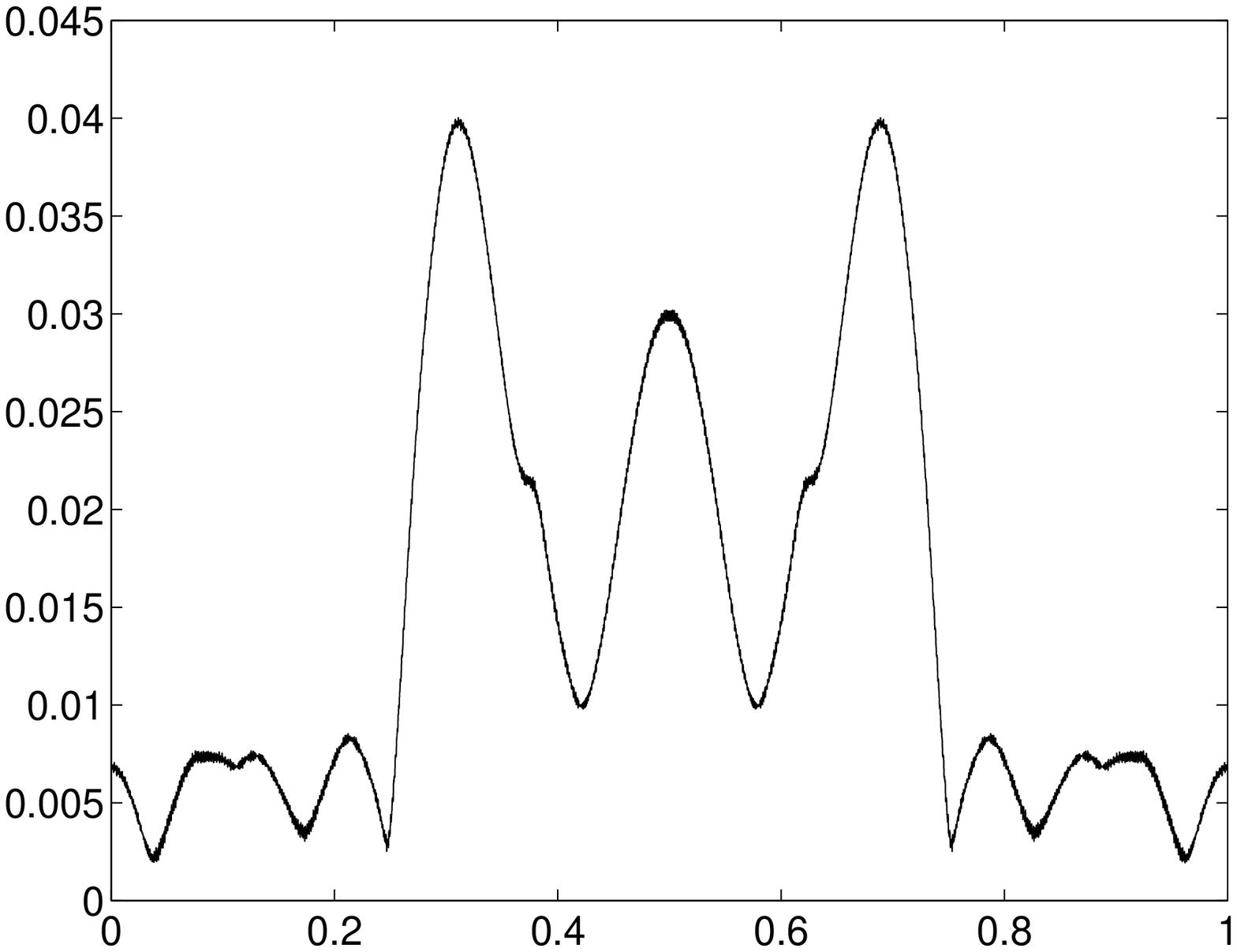}}
\resizebox{1.5in}{!}{\includegraphics{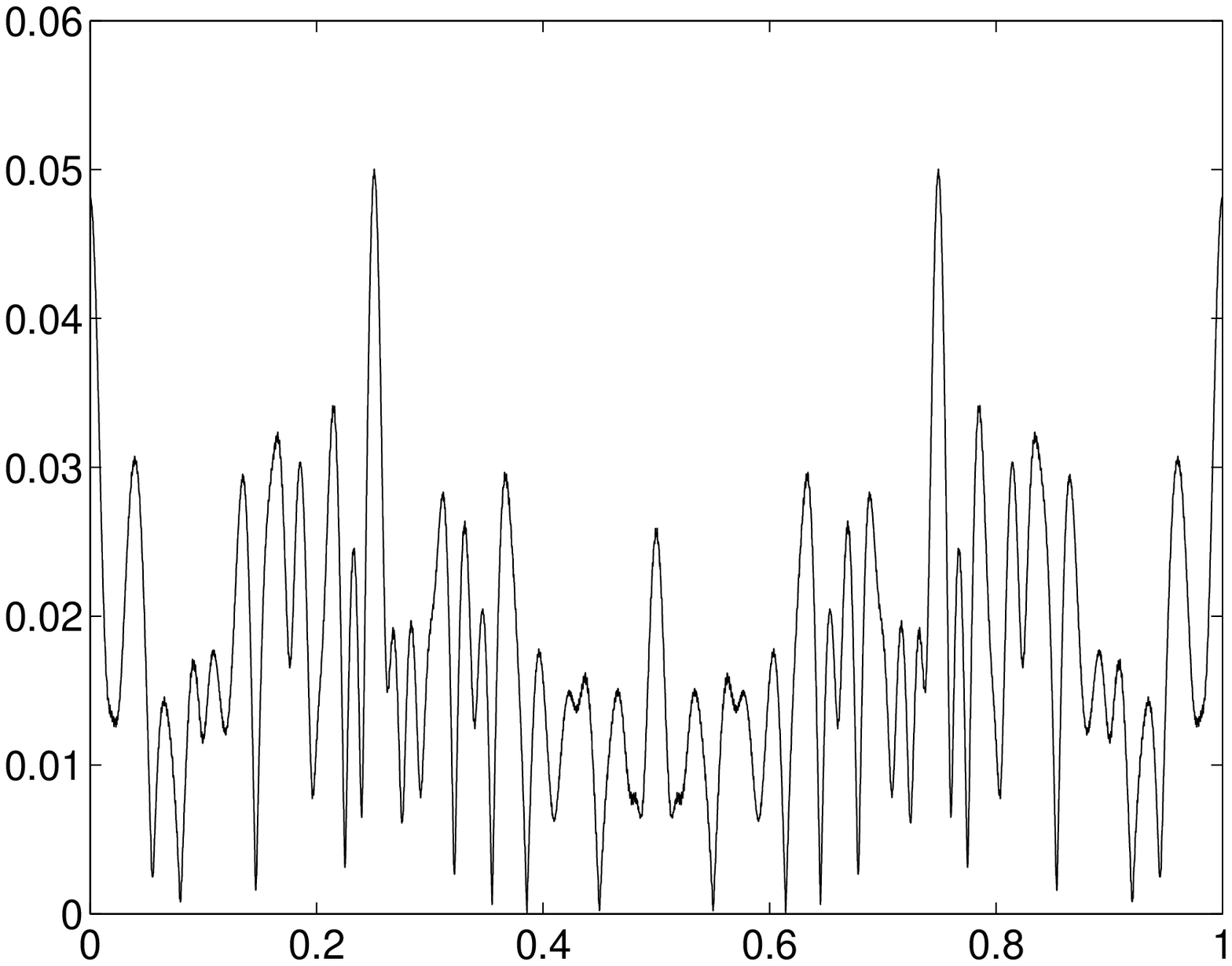}}

$|\psi^{\rm ex}(t,x)|^2$, $|\psi^{\rm ex}(t,x)-\psi^{\rm
ts}(t,x)|$ and $|\psi^{\rm ex}(t,x)-\psi^{\rm bd}(t,x)|$ at
$t=1.0$\vspace{-1mm}
\caption{Numerical results for example
\ref{nuex2} with $U(x)$ given by \eqref{eqstep}, $\e=\f{1}{2}$. We
use $\tg t=\f{1}{100}$, $\tg x=\f{1}{32}$ for the TS, $\tg
t=\f{1}{5}$, $\tg x=\f{1}{32}$ for the Bloch-decomposition method,
and $\tg t=\f{1}{100000}$, $\tg x=\f{1}{8192}$for the ``exact''
solution. \vspace{-3mm}}\label{fig24}
\[\Delta^{\rm ts}_\infty(t)=4.01\E-2,\
\Delta^{\rm bd}_\infty(t)=5.00\E-2,\
\Delta^{\rm ts}_2(t)=1.85\E-2,\
\Delta^{\rm bd}_2(t)=1.98\E-2.\]\vspace{1mm}

\resizebox{1.5in}{!}{\includegraphics{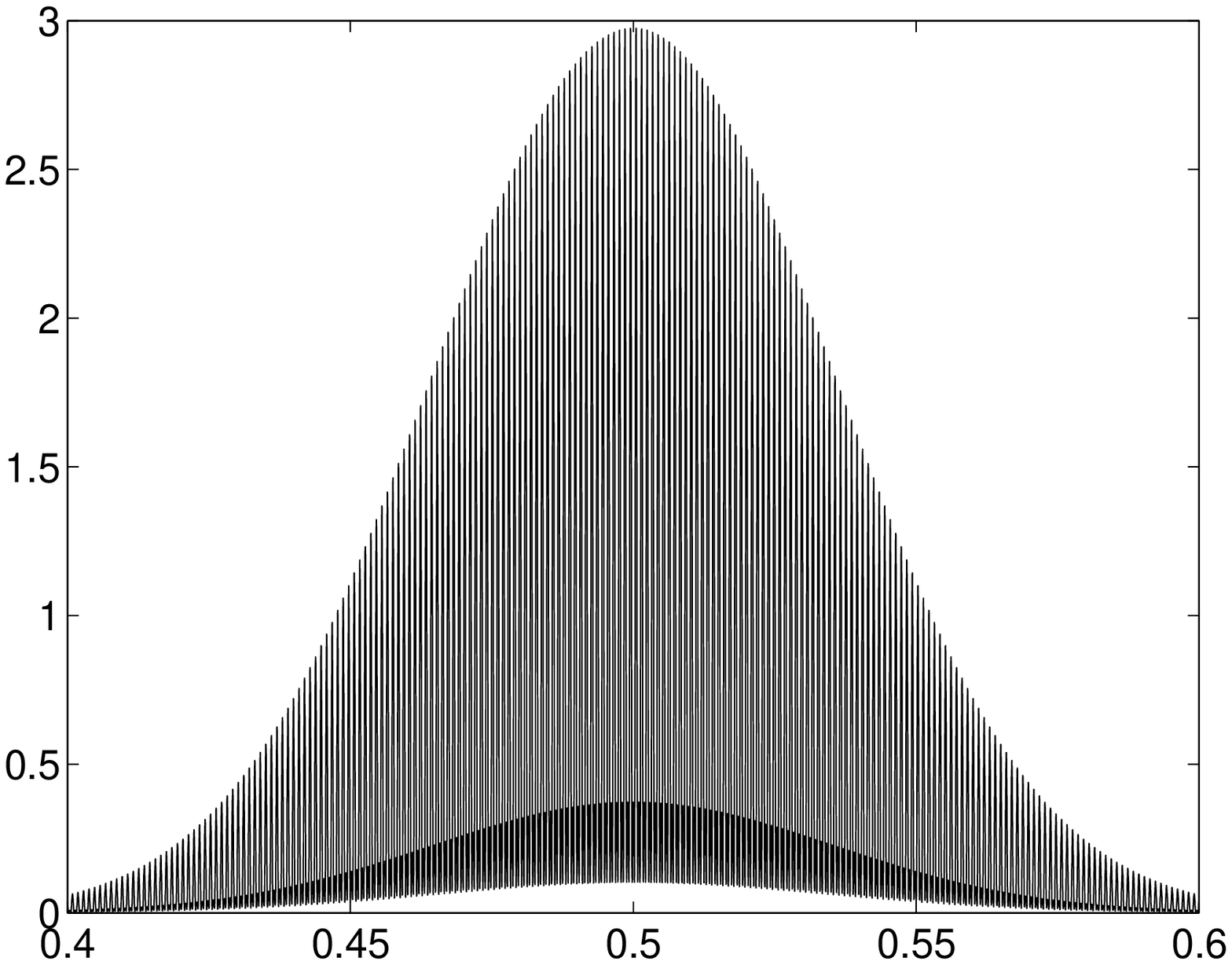}}
\resizebox{1.5in}{!}{\includegraphics{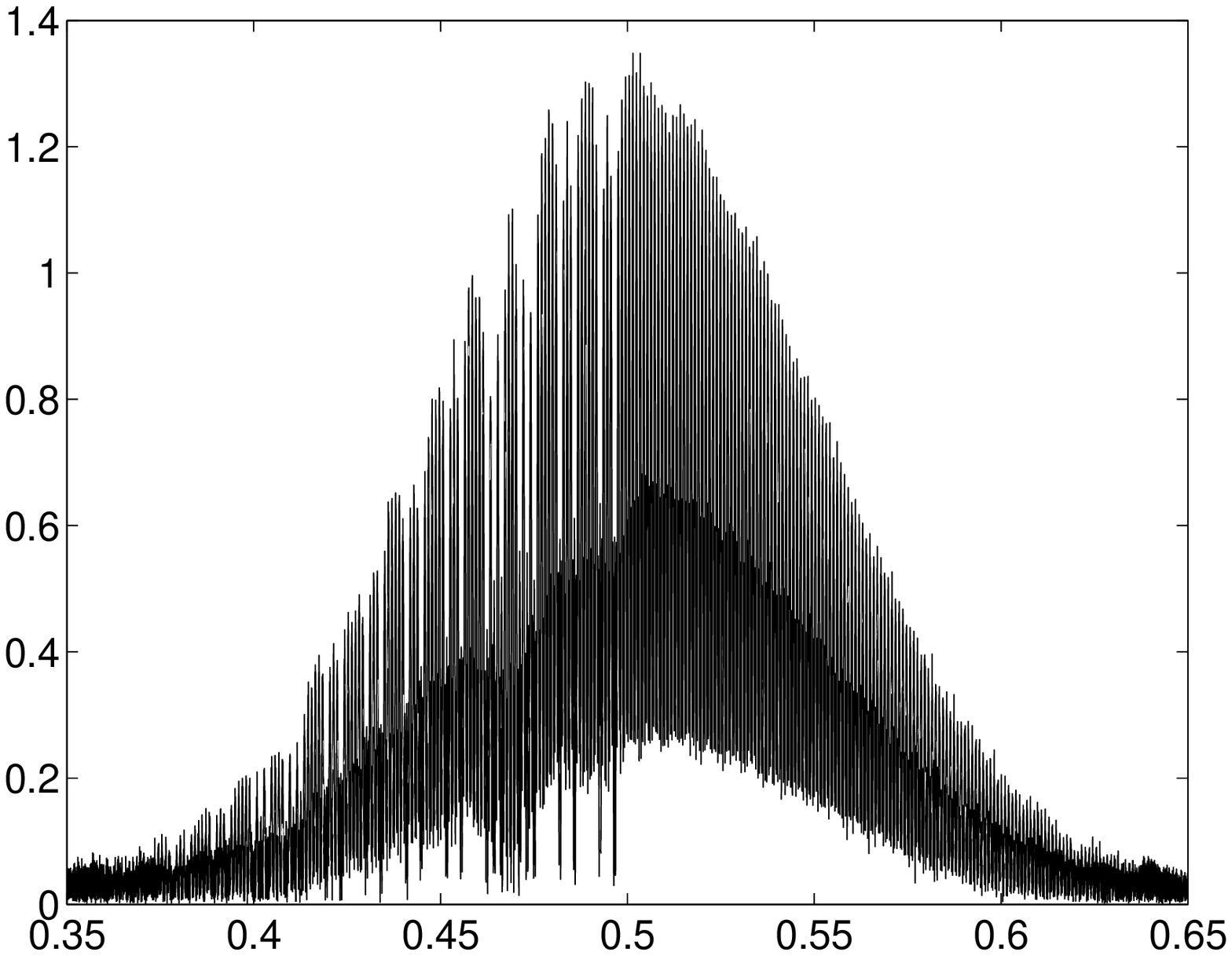}}
\resizebox{1.5in}{!}{\includegraphics{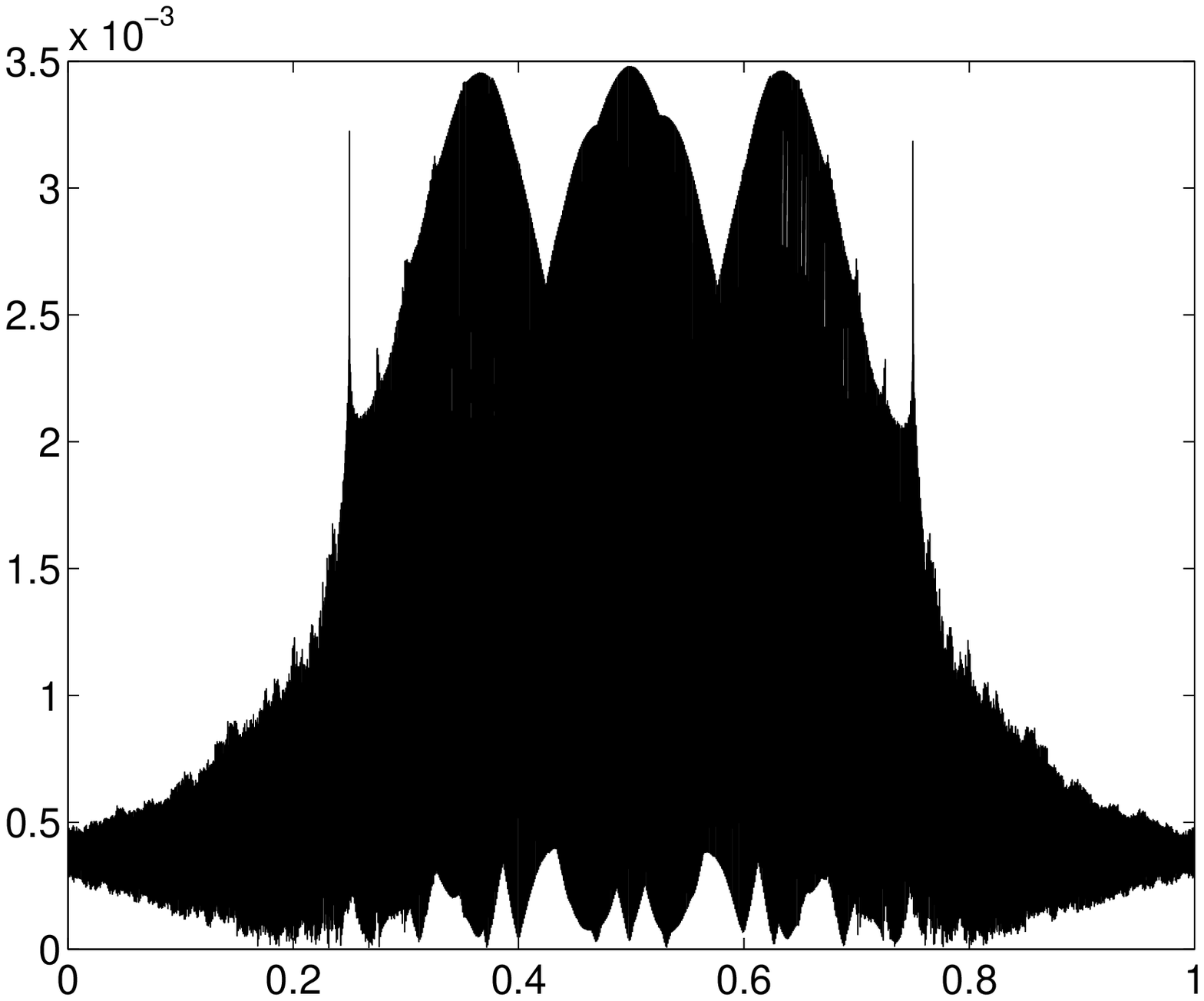}}

$|\psi^{\rm ex}(t,x)|^2$, $|\psi^{\rm ex}(t,x)-\psi^{\rm ts}(t,x)|$
and $|\psi^{\rm ex}(t,x)-\psi^{\rm bd}(t,x)|$ at $t=0.1$\vspace{-1mm}
\end{center}
\caption{Numerical results for example \ref{nuex2}. Here
$U(x)$ is given in \eqref{eqstep}, $\e=\f{1}{1024}$.
We use $\tg t=\f{1}{10000}$, $\tg x=\f{1}{65536}$ for
Time-splitting,
$\tg t=\f{1}{10}$, $\tg x=\f{1}{8192}$ for Bloch-decomposition,
and $\tg t=\f{1}{100000}$, $\tg
x=\f{1}{131072}$ for 'exact' solution.
}\label{fig26}\vspace{-3mm}
\[\Delta^{\rm ts}_\infty(t)=1.35,\
\Delta^{\rm bd}_\infty(t)=3.48\E-3,\
\Delta^{\rm ts}_2(t)=2.23\E-1,\
\Delta^{\rm bd}_2(t)=1.14\E-3.\]
\end{figure}

\begin{figure} \footnotesize
\begin{center}
\resizebox{1.5in}{!}{\includegraphics{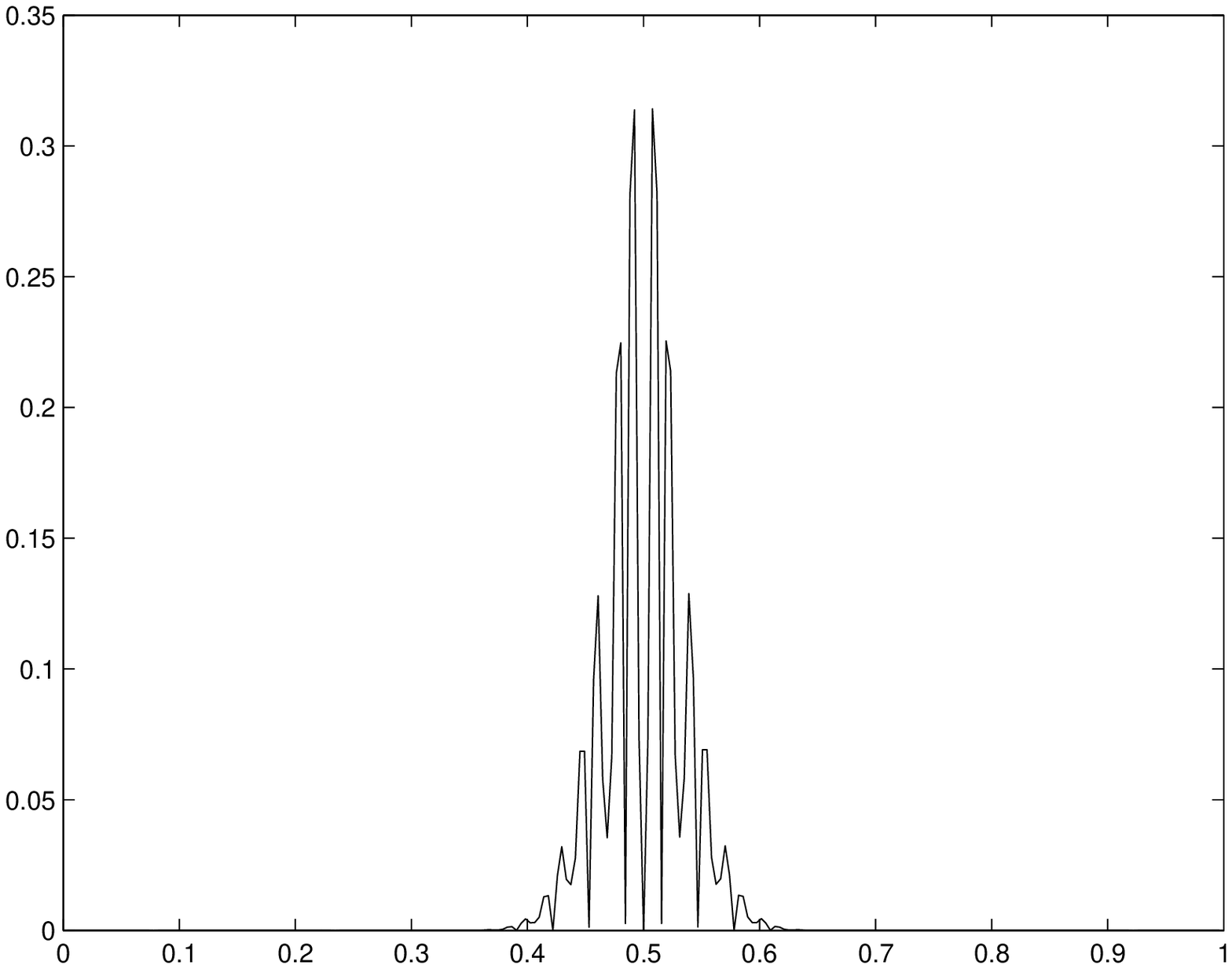}}
\resizebox{1.5in}{!}{\includegraphics{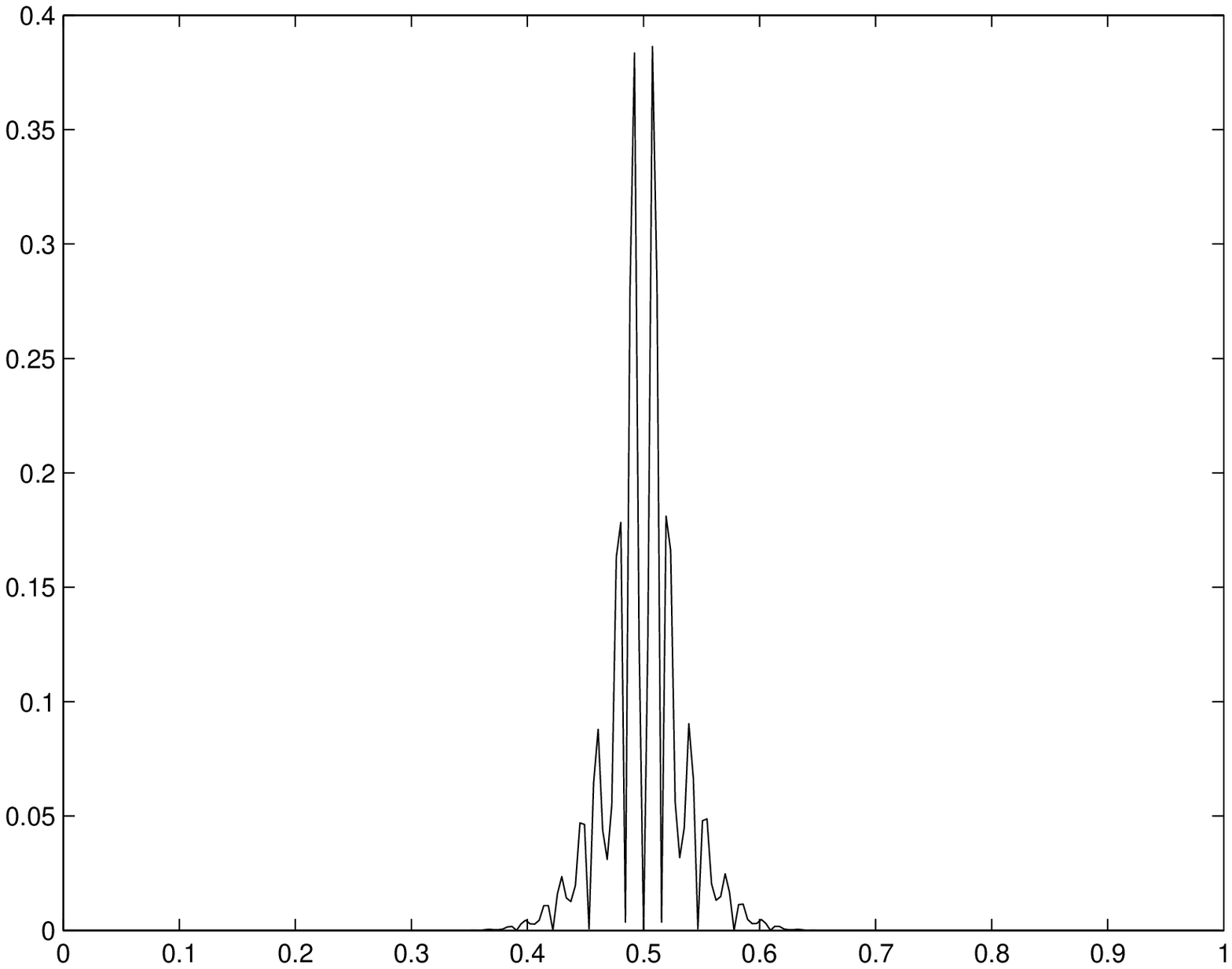}}
\resizebox{1.5in}{!}{\includegraphics{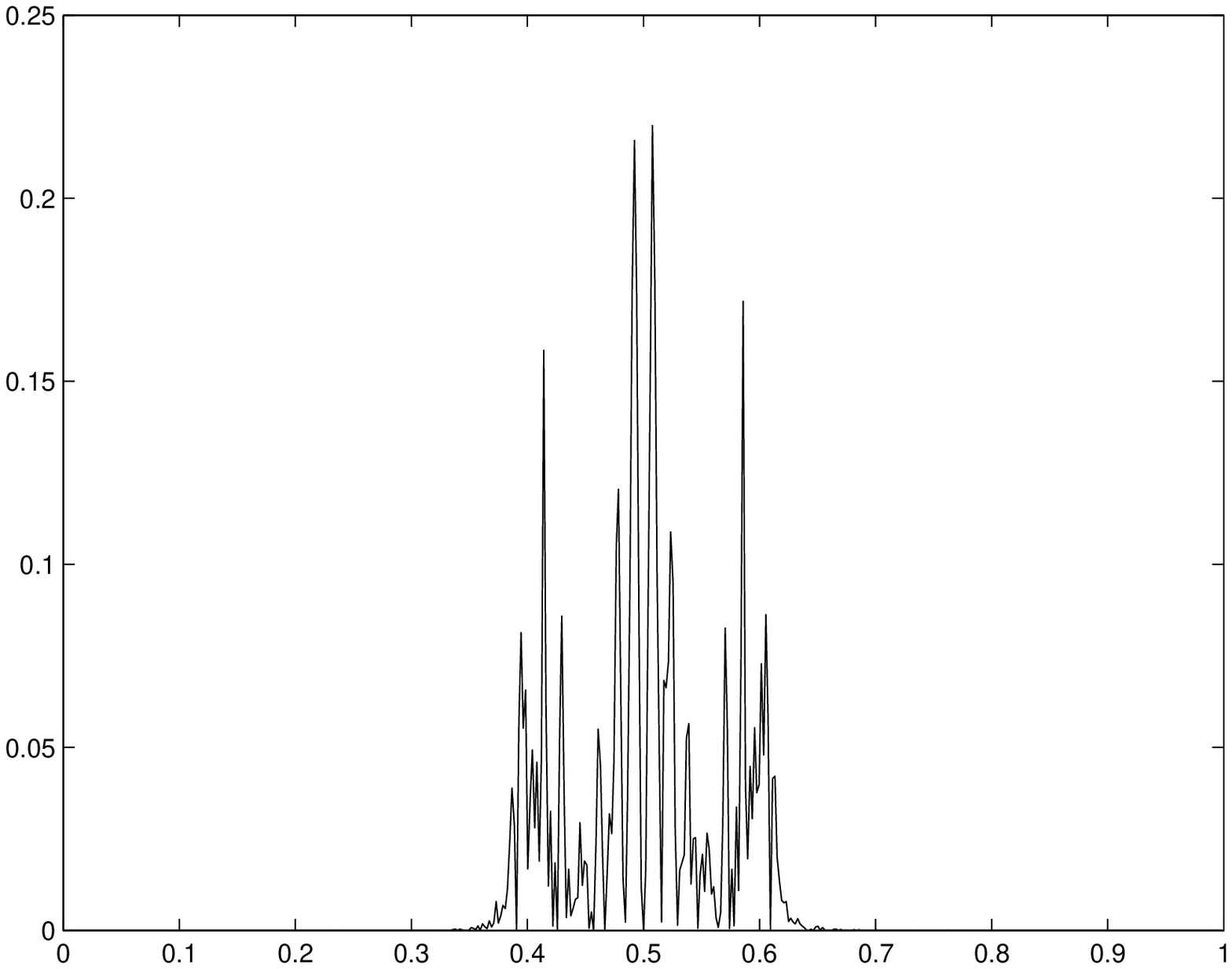}}

$|\psi^{\rm sc}( x,0.1)|^2$, $|\psi^{\rm sc}( x,0.25)|^2$, and
$|\psi^{\rm sc}( x,0.5)|^2$, $\e=\f{1}{32}$.\vspace{4mm}

\resizebox{1.5in}{!}{\includegraphics{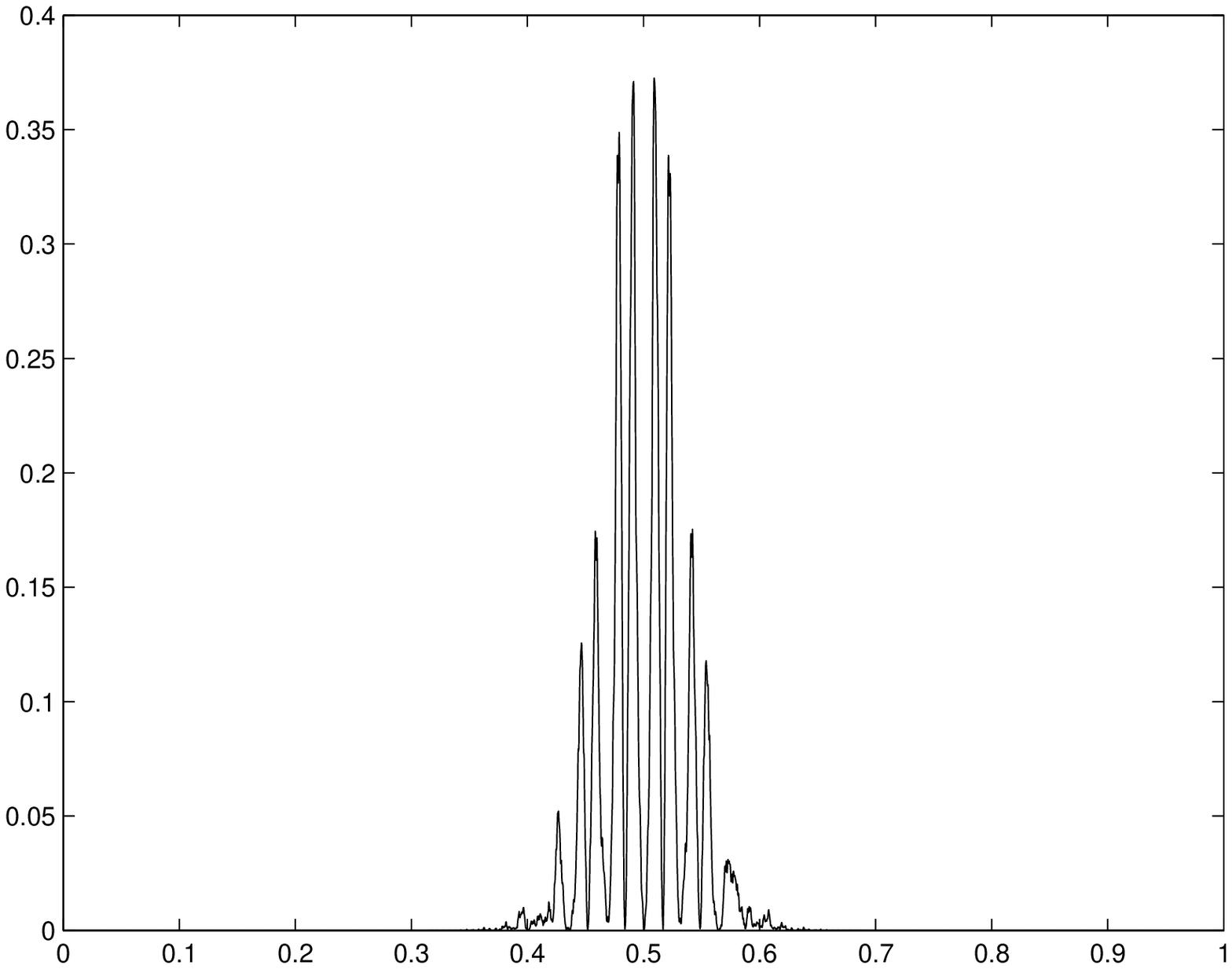}}
\resizebox{1.5in}{!}{\includegraphics{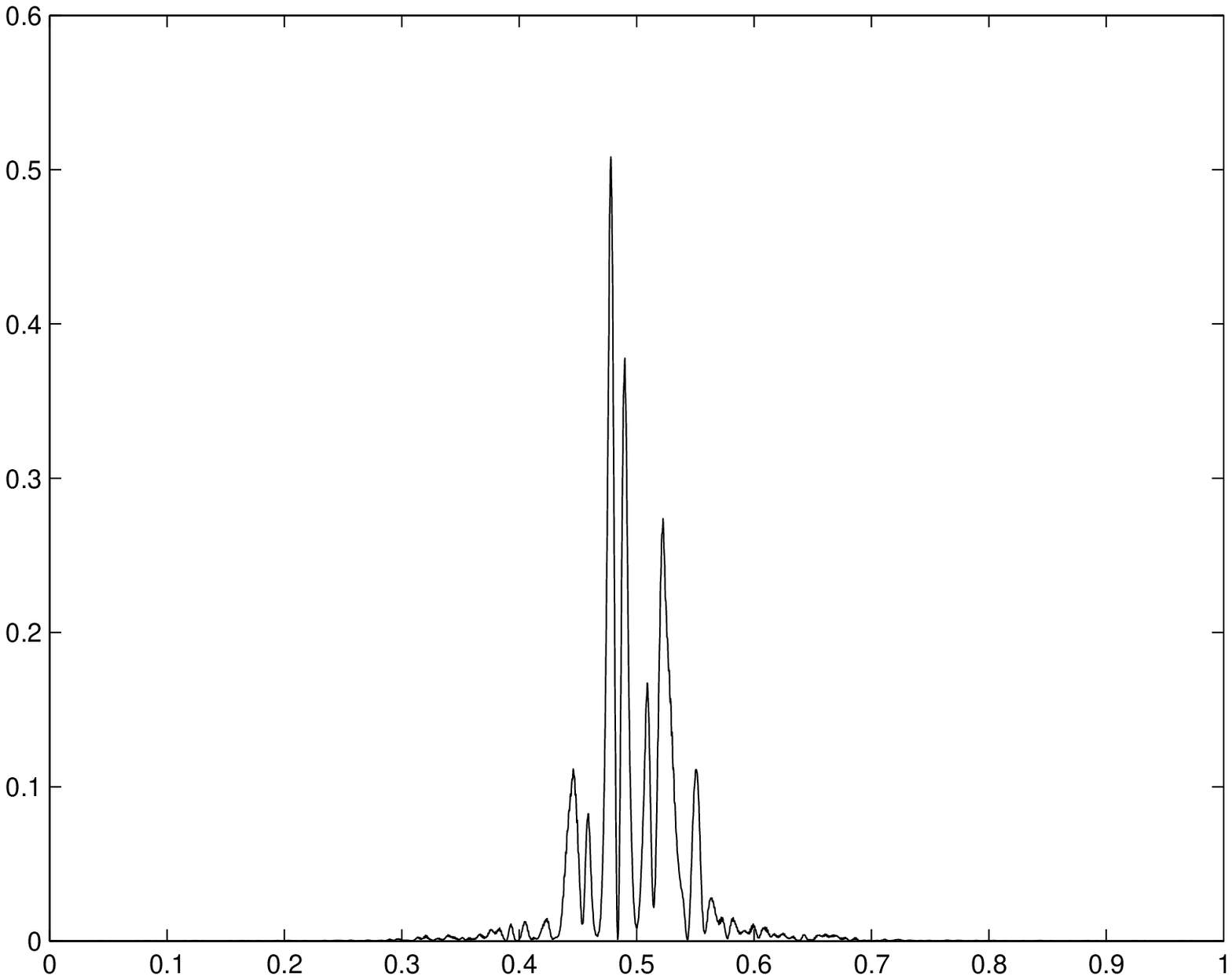}}
\resizebox{1.5in}{!}{\includegraphics{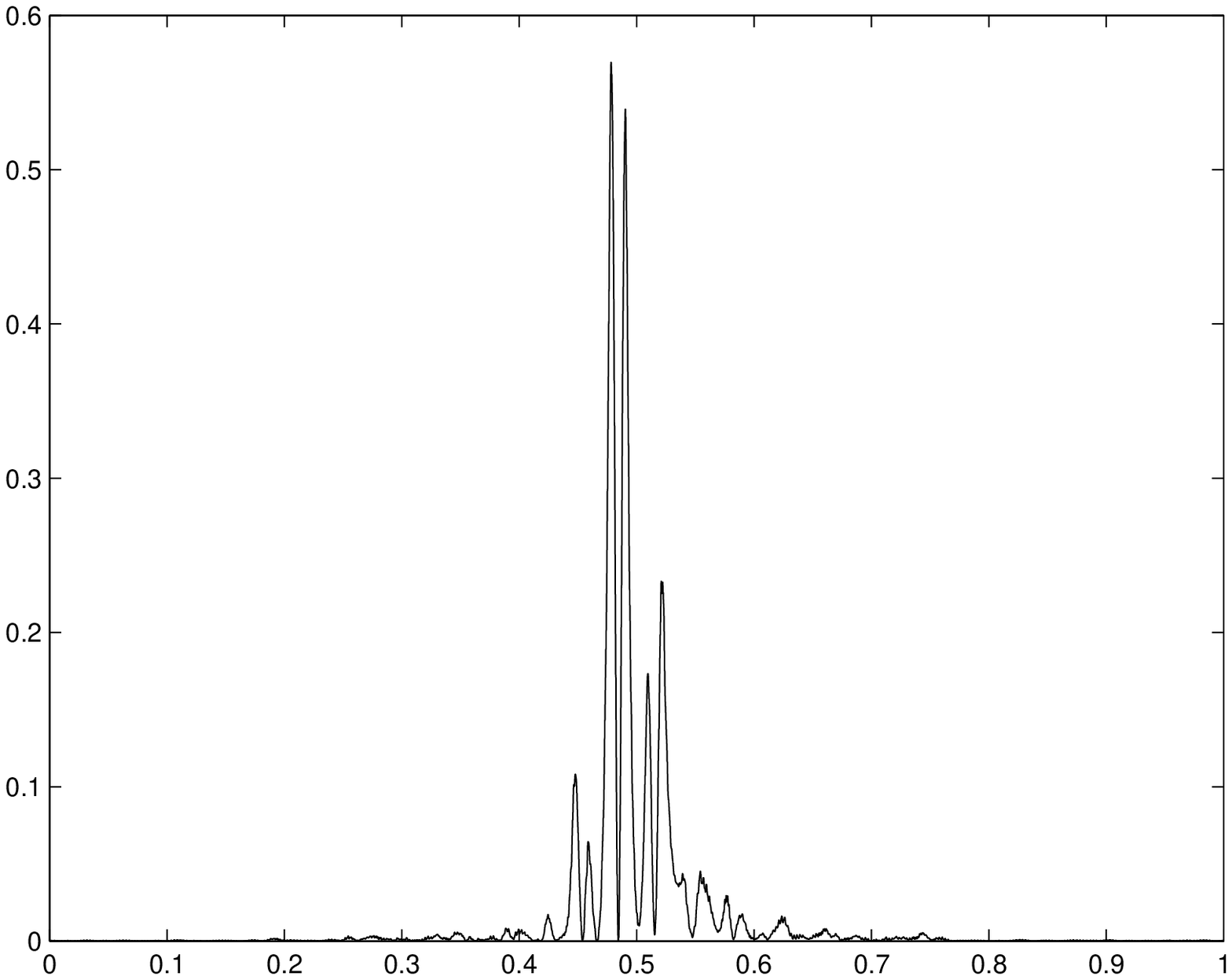}}

$|\psi^{\rm bd}( x,0.1)|^2$, $|\psi^{\rm bd}( x,0.25)|^2$, and
$|\psi^{\rm bd}( x,0.5)|^2$, $\e=\f{1}{32}$.\vspace{4mm}

\resizebox{1.5in}{!}{\includegraphics{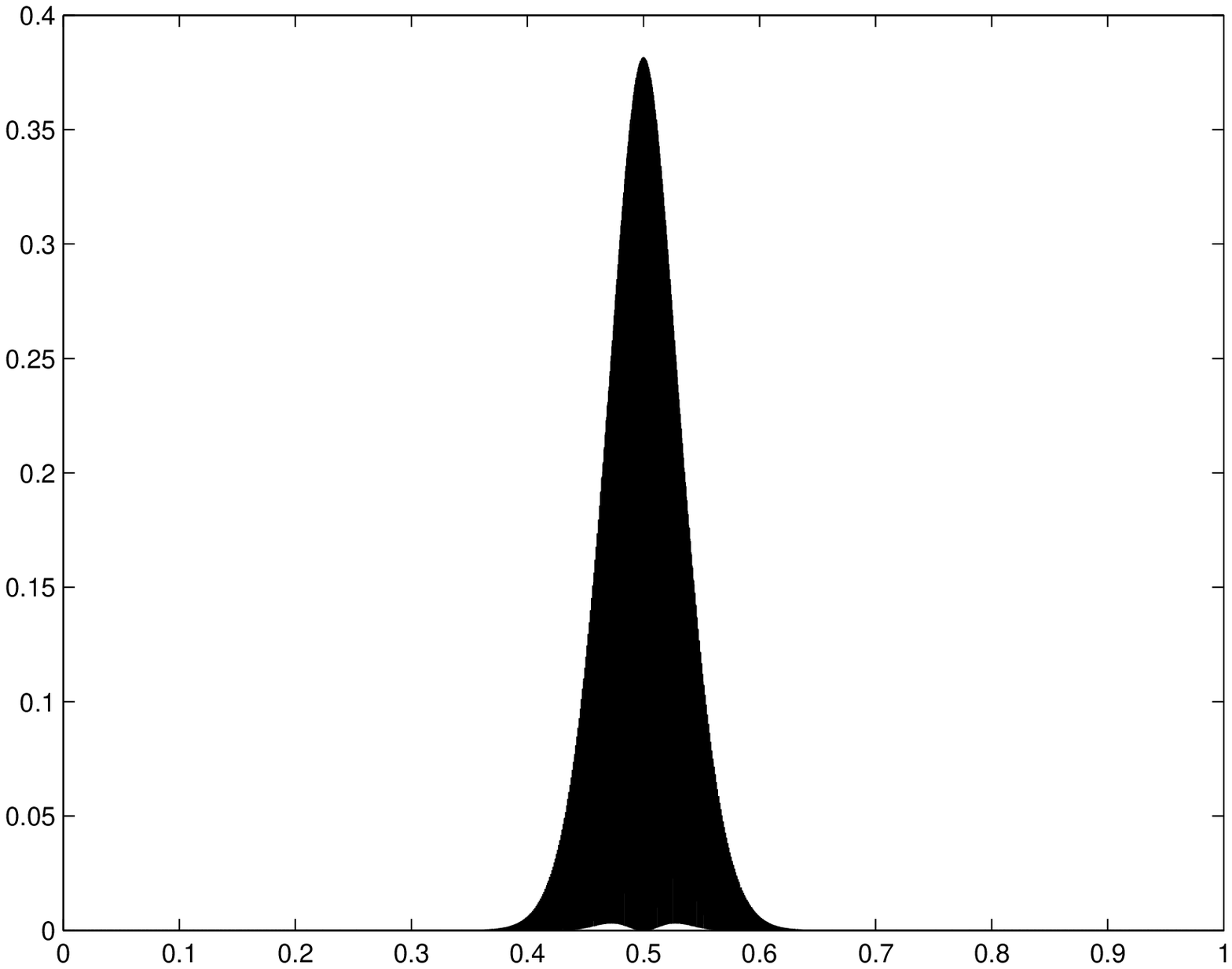}}
\resizebox{1.5in}{!}{\includegraphics{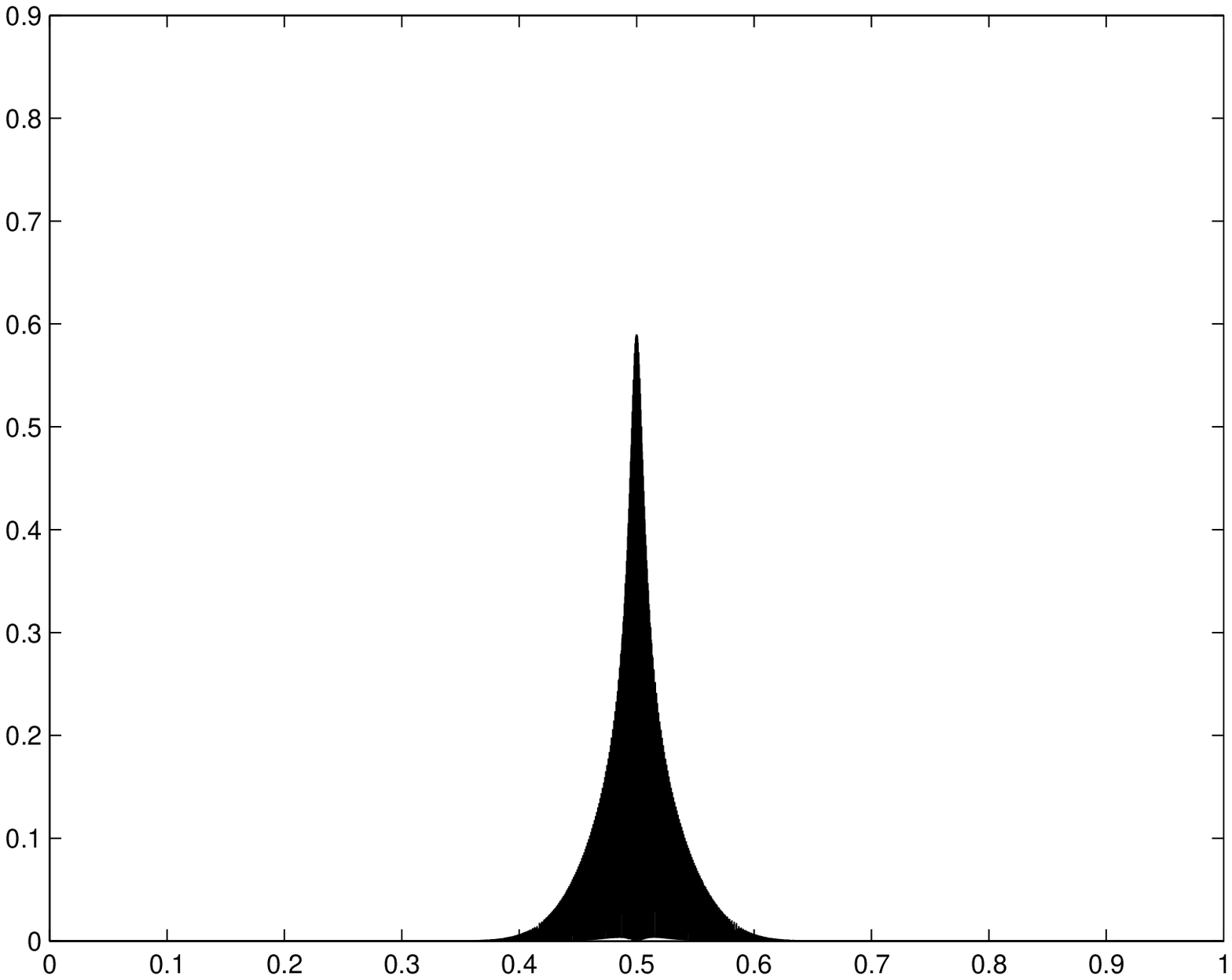}}
\resizebox{1.5in}{!}{\includegraphics{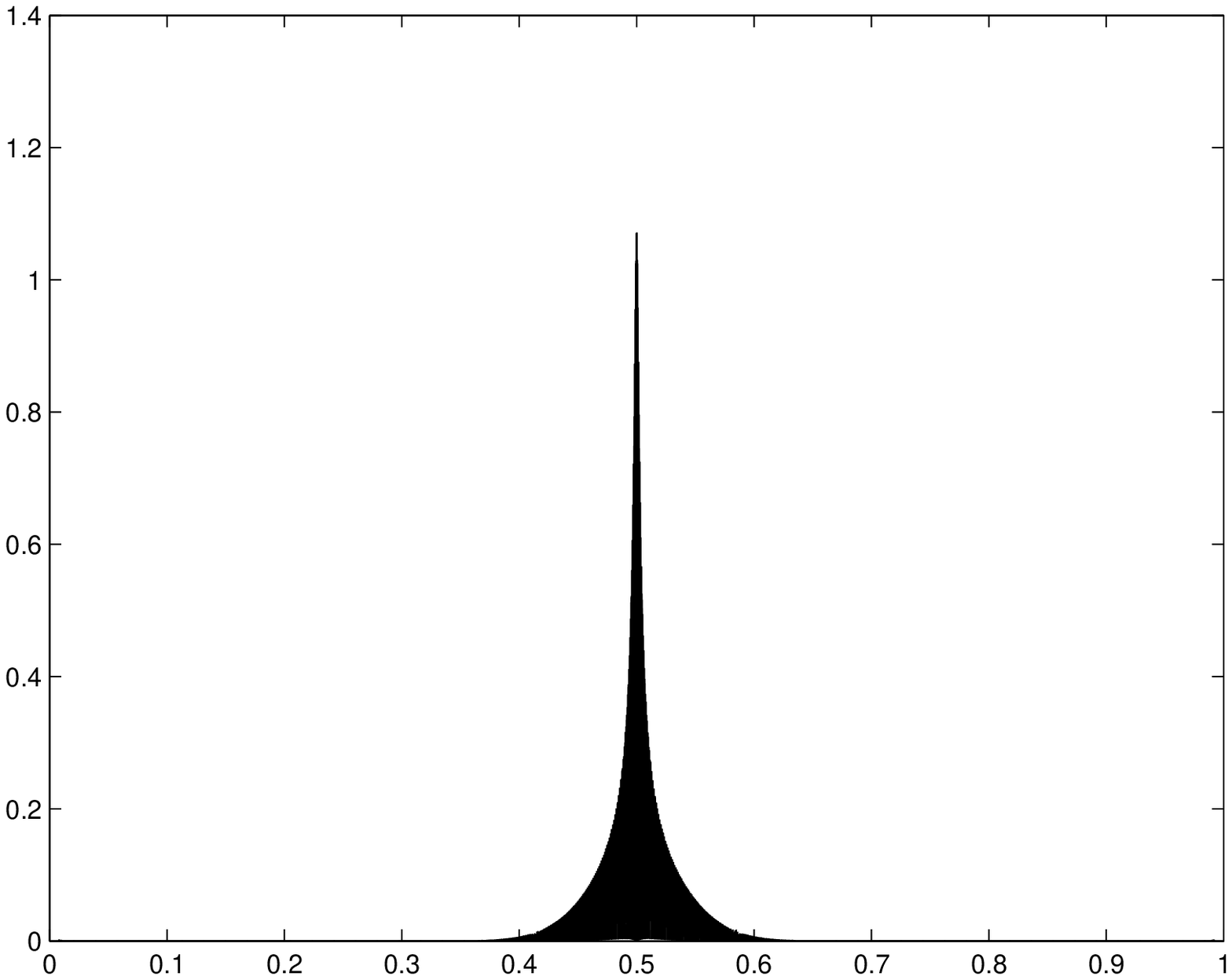}}

$|\psi^{\rm sc}( x,0.1)|^2$, $|\psi^{\rm sc}( x,0.25)|^2$, and
$|\psi^{\rm sc}( x,0.5)|^2$, $\e=\f{1}{1024}$.\vspace{4mm}

\resizebox{1.5in}{!}{\includegraphics{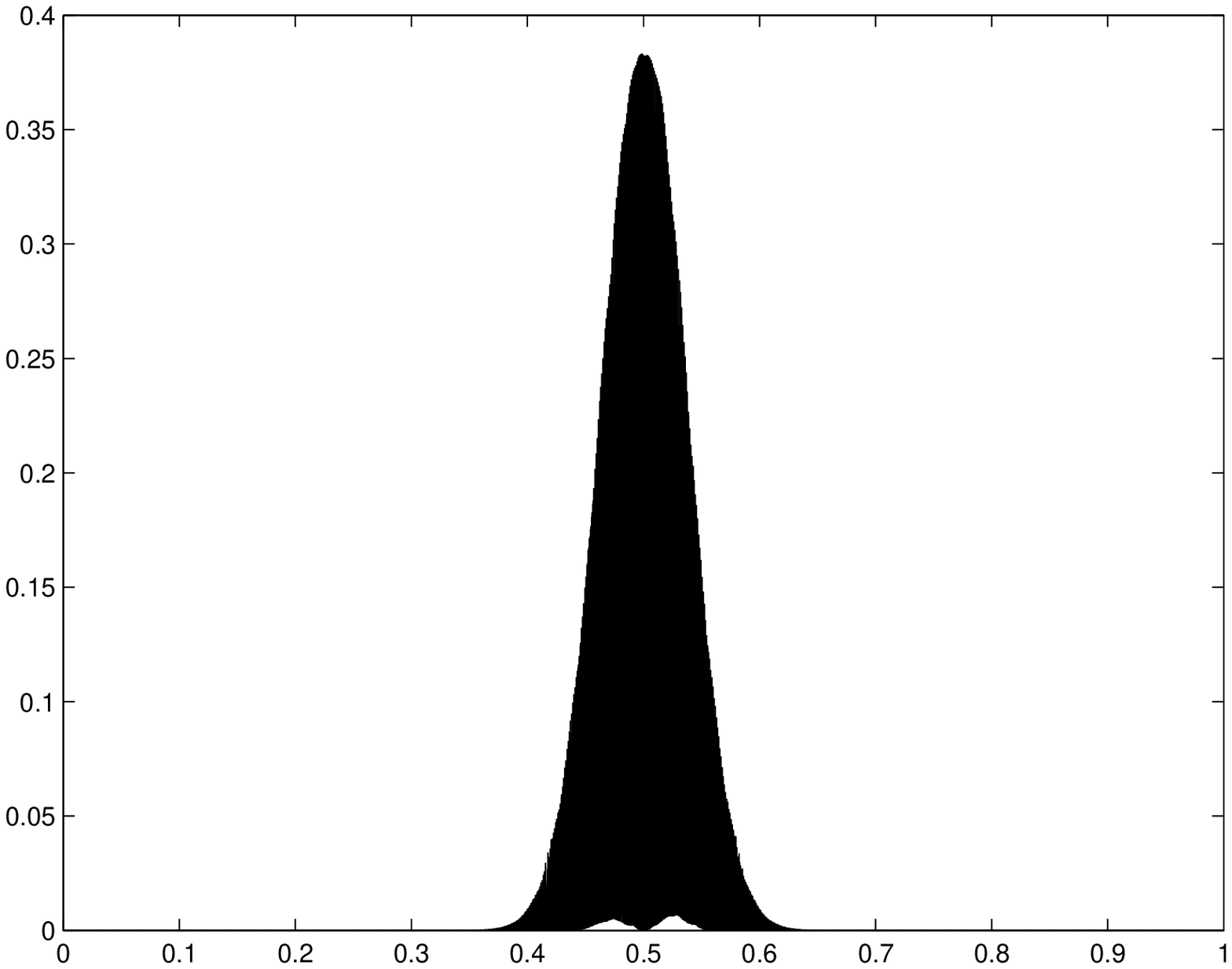}}
\resizebox{1.5in}{!}{\includegraphics{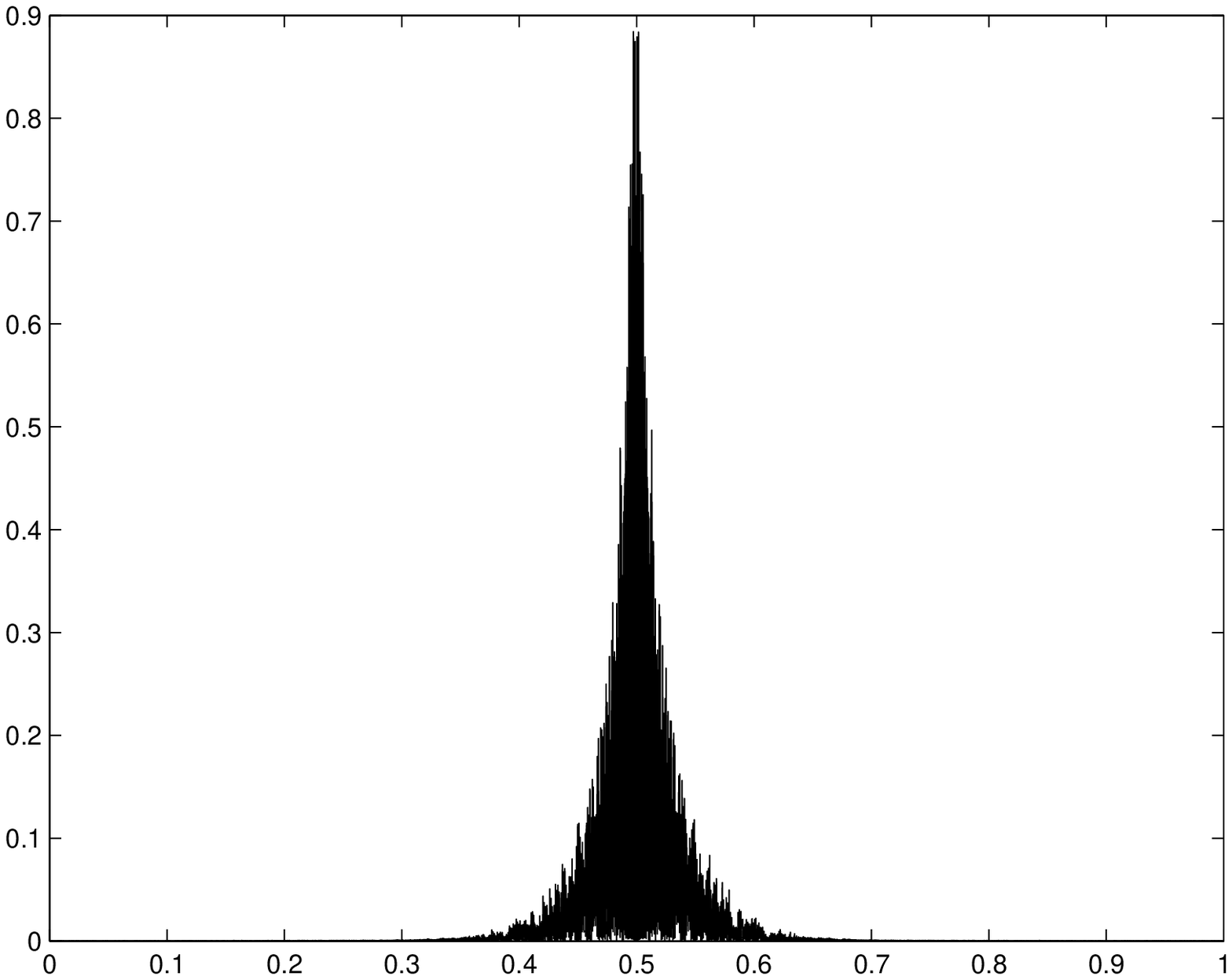}}
\resizebox{1.5in}{!}{\includegraphics{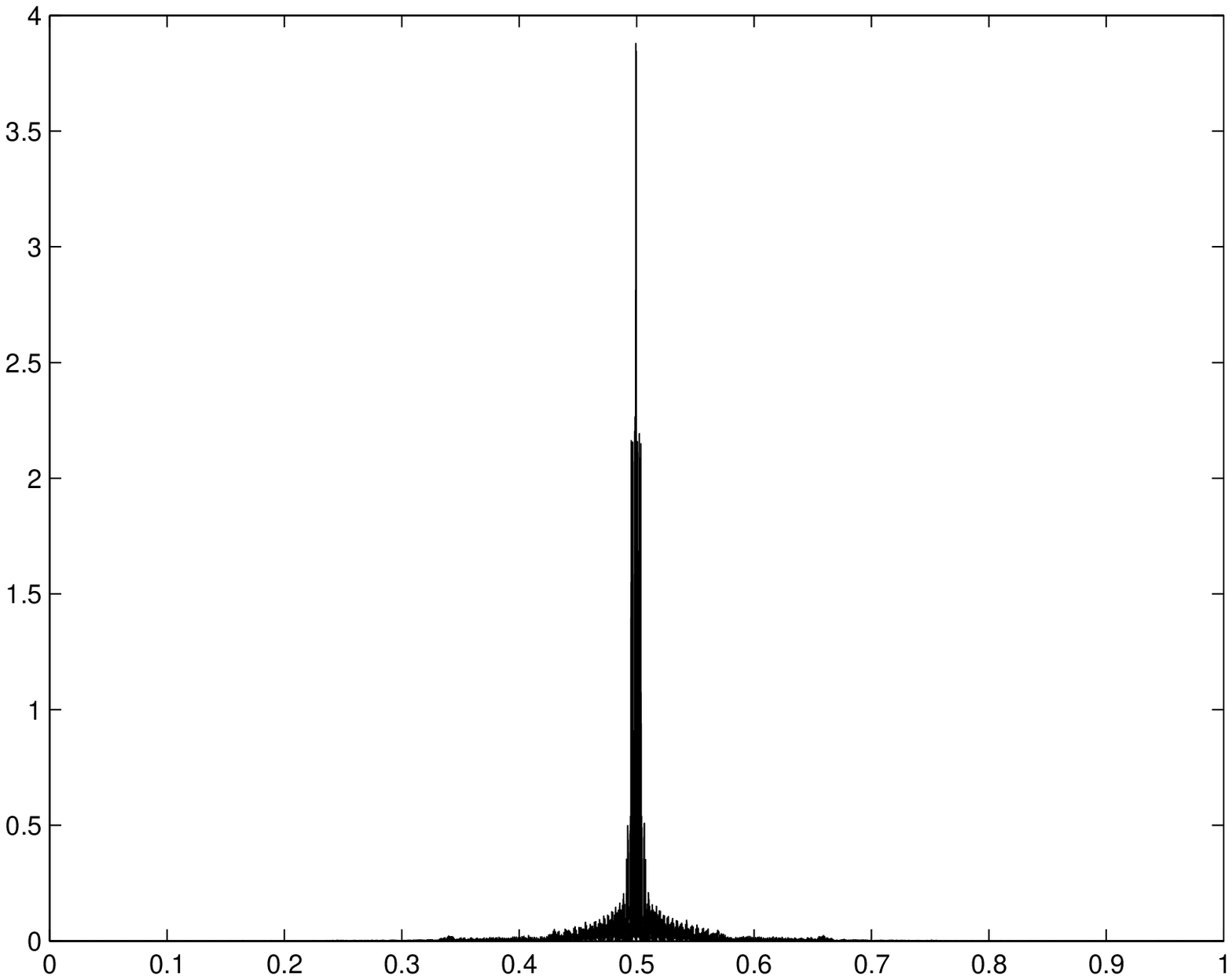}}

$|\psi^{\rm bd}( x,0.1)|^2$, $|\psi^{\rm bd}( x,0.25)|^2$, and
$|\psi^{\rm bd}( x,0.5)|^2$, $\e=\f{1}{1024}$.
\end{center}
\caption{Numerical results for example \ref{nuscex2} with $U(x)$
given by \eqref{eqhar_p}, $\tg t=\f{1}{10000}$, $\tg
x=\f{1}{32768}$. The left column shows the situation before the
caustic, whereas the other two columns respectively present the
numerical results at and after the caustic.}\label{fig41}
\end{figure}

\begin{thebibliography}{99}

\bibitem{AsKn}  {\sc J.\ Asch and A.\ Knauf}, {\em Motion in periodic
potentials}, Nonlinearity {\bf 11} (1998), 175--200.

\bibitem{AsMe} {\sc N.~W. Ashcroft and N.~D. Mermin},
\emph{Solid state physics}, Saunders New York, 1976.

\bibitem {BJM} {\sc W.~Z. Bao, S. Jin, and P. Markowich},
\emph{On time-splitting spectral approximations for the
Schr\"odinger equation in the semiclassical regime},
J. Comp. Phys. {\bf 175} (2002), 487--524.

\bibitem {BJM1} {\sc W.~Z. Bao, S. Jin, and P. Markowich},
\emph{Numerical study of time-splitting spectral discretizations
of nonlinear Schr\"odinger equations in the semi-Classical regime},
SIAM J. Sci. Comp. {\bf 25} (2003), 27--64.

\bibitem {BMP} {\sc P. Bechouche, N. Mauser, and F. Poupaud},
\emph{Semiclassical limit for the
Schr\"odinger-Poisson equation in a crystal}, Comm. Pure Appl. Math.
\textbf{54} (2001),  no. 7, 851--890.

\bibitem {BePo} {\sc P. Bechouche,and F. Poupaud},
\emph{Semi-classical limit of a Schr\"odinger equation for a
stratified material}, Monatsh. Math. {\bf 129} (2000), no. 4,
281--301.

\bibitem {BLP} A. Bensoussan, J.~L. Lions, and G. Papanicolaou,
\emph{Asymptotic Analysis for Periodic Structures}, North-Holland
Pub. Co. (1978).

\bibitem {B} F. Bloch, \emph{\"Uber die Quantenmechanik der Elektronen
in Kristallgittern}, Z. Phys. \textbf{52} (1928), 555--600.

\bibitem{Bl} E.\ I.\ Blount, {\em Formalisms of band theory},
Solid State Physics {\bf 13}, Academic Press, New York, 305--373
(1962).

\bibitem{Bu} K. Busch, \emph{Photonic band structure theory:
assessment and perspectives}, Compte Rendus Physique {\bf 3}
(2002), 53--66.

\bibitem {Ca} R. Carles, \emph{WKB analysis for nonlinear
Schr\"odinger equations with a potential}, Comm. Math. Phys. to appear.

\bibitem {CMS} R. Carles, P. A. Markowich and C. Sparber,
\emph{Semiclassical asymptotics for weakly nonlinear
Bloch waves}, J. Stat. Phys. \textbf{117} (2004), 369--401.

\bibitem{COV} C. Conca, R. Orive, and M. Vanninathan,
\emph{Bloch approximation in homogenization on bounded domains},
Asymptot. Anal. {\bf 41} (2005), no. 1, 71--91.

\bibitem{CSV} C. Conca, N. Srinivasan and M. Vanninathan,
\emph{Numerical solution of elliptic partial differential
equations by Bloch waves method}, in:  Congress on Differential
Equations and Applications/VII CMA (Salamanca, 2001), 63--83,
2001.

\bibitem{CoVa} C. Conca and M. Vanninathan,
\emph{Homogenization of periodic structures via Bloch
decomposition}, SIAM J. Appl. Math. {\bf 57} (1997), no. 6,
1639--1659.

\bibitem{FL} M.~V. Fischetti and S. E. Laux,
\emph{Monte Carlo analysis of electron transport in small semiconductor
devices including band-structure and space-charge effects},
Phys. Rev. B {\bf 38} (1998), 9721--9745.

\bibitem {GMMP} P. G\'erard, P. Markowich, N. Mauser, and F. Poupaud,
\emph{Homogenization Limits and Wigner transforms},
Comm. Pure and Appl. Math {\bf 50} (1997), 323--378.

\bibitem{Go} L. Gosse, \emph{Multiphase semiclassical approximation
of an electron in a one-dimensional crystalline lattice. II.
Impurities, confinement and Bloch oscillations}, J. Comput. Phys.
{\bf 201} (2004), no. 1, 344--375.

\bibitem {GoMa} L. Gosse and P.~A. Markowich, \emph{Multiphase
semiclassical approximation of an electron in a one-dimensional
crystalline lattice - I. Homogeneous problems}, J. Comput Phys.
\textbf{197} (2004), 387--417.

\bibitem {GoMau} L. Gosse and N. Mauser, \emph{Multiphase semiclassical
approximation of an electron in a one-dimensional crystalline
lattice. III. From ab initio models to WKB for
Schr\"odinger-Poisson}, to appear in J. Comput. Phys. {\bf 211}
(2006), no. 1, 326--346.

\bibitem {GRT} J. C. Guillot, J. Ralston, and E. Trubowitz,
\emph{Semiclassical asymptotics in solid-state physics}, Comm.
Math. Phys. \textbf{116} (1998), 401--415.

\bibitem{HFKW} D. Hermann, M. Frank, K. Busch, and P. W\"olfle,
\emph{Photonic band structure computations},
Optics Express {\bf 8} (2001), 167--173.

\bibitem{HJ}R. Horn and C. Johnson, \emph{Matrix analysis},
Cambridge University Press, Cambridge, 1985.

\bibitem{HJMSZ} Z. Huang, S. Jin, P. Markowich, C. Sparber and C. Zheng,
\emph{A Time-splitting spectral scheme for the Maxwell-Dirac
system}, J. Comput. Phys. {\bf 208} (2005), issue 2, 761--789.

\bibitem {JiXi} S. Jin, Z. Xin, \emph{Numerical passage from systems of
conservation laws to Hamilton-Jacobi equations, and a relaxation
scheme}, SIAM J. Num. Anal. \textbf{35} (1998), 2385--2404.

\bibitem{JC} J.~D. Joannopoulos and M.~L. Cohen, \emph{Theory of
Short Range Order and Disorder in Tetrahedrally Bonded
Semiconductors}, Solid State Physics {\bf 31} (1974), 1545.

\bibitem{Ko} H.~J. Korsch and M. Gl\"uck,
\emph{Computing quantum eigenvalues made easy}, Eur. J. Phys.
{\bf 23} (2002), 413--425.

\bibitem{LVF} S.~E. Laux, M.~V. Fischetti, and D.~J. Frank,
\emph{Monte Carlo analysis of semiconductor devices: the DAMOCLES
program}, IBM Journal of Research and Development {\bf 34} (1990),
466--494.

\bibitem{Lu}  J.M.\ Luttinger, {\em The effect of a magnetic field
on electrons in a periodic potential},
Phys.\ Rev.\ {\bf 84} (1951), 814--817 .

\bibitem {PST} G. Panati, H. Spohn, and S. Teufel,
\emph{Effective dynamics for Bloch electrons: Peierls substitution
and beyond}, Comm. Math. Phys. \textbf{242} (2003), 547--578.

\bibitem {ReSi} M. Reed, B. Simon, \emph{Methods of modern mathematical
physics IV. Analysis of operators}, Academic Press (1978).

\bibitem {Te} S. Teufel, \emph{Adiabatic perturbation theory in quantum
dynamics}, Lecture Notes in Mathematics 1821, Springer (2003).

\bibitem{Wi}  C.\ H.\ Wilcox, {\em Theory of bloch waves},
J.\ Anal.\ Math. {\bf 33} (1978), 146--167.

\bibitem{Za} J.\ Zak, {\em Dynamics of electrons in solids in external
fields}, Phys.\ Rev.\ {\bf 168} (1968), 686--695.

\bibitem{Ze} A. Zettel, \emph{Spectral theory and computational methods
for Sturm-Liouville problems}, in D. Hinton and P.~W, Sch\"afer
(eds.). Lecture Notes in Pure and Applied Math. {\bf 191}, Dekker
1997.
\end{thebibliography}
\end{document}